\catcode\endlinechar = 10%
%
\voffset = -5.4mm%
\vsize = 257mm%
\hoffset = -5.4mm%
\hsize = 170mm%
\font\fourteenrm = cmr12 at 14pt%
%
%
%
%
%
%
%
%
%
%
%
%
%
%
%
%
%
\font\sixteensc = cmcsc10 at 16pt%
\font\twelvesc = cmcsc10 at 12pt%
%
%
%
%
\rm%
\parindent = 1cm%
\parskip = 5pt%
\def\monthname{%
	\ifcase\month%
		\or January%
		\or February%
		\or March%
		\or April%
		\or May%
		\or June%
		\or July%
		\or August%
		\or September%
		\or October%
		\or November%
		\or December%
	\fi%
}%
\footline = {\hss\rm\folio\hss}%
\def\Document{article}%
\def\Title#1{{%
	\parskip = 0pt%
	\centerline{\sixteensc #1}%
	\baselineskip = 28.8pt%
	\centerline{\fourteenrm B.R.S.~Roso}%
	\baselineskip = 16.8pt%
	\centerline{\monthname\ \number\year}%
	\vskip 16.8pt%
	\baselineskip = 21.6pt%
}}%
\newcount\SectNum%
\newcount\EltNum%
\SectNum = 0%
\def\Abs{{%
	\vskip\baselineskip%
	\twelvesc%
	\centerline{Abstract}%
}}%
\def\Ack{{%
	\vskip\baselineskip%
	\twelvesc%
	\centerline{Acknowledgements}%
}}%
\def\Sect#1{%
	\advance\SectNum by 1%
	\EltNum = 1%
	{%
		\vskip\baselineskip%
		\goodbreak%
		\twelvesc%
		\centerline{\number\SectNum.\enskip#1}%
		\nobreak%
	}%
}%
\long\def\Elt#1#2{%
	{%
		\parindent = 0cm%
		\goodbreak%
		\par%
		\def\EltName{%
			\number\SectNum%
			.%
			\number\EltNum%
		}%
		\def\EltFullName{%
			#1%
			\ 
			\EltName
		}%
		{\bf \EltFullName:\ }%
		{\ignorespaces#2}%
		\goodbreak%
		\par%
	}%
	\advance\EltNum by 1%
}%
\def\Label#1{%
	\global\edef#1{{\EltFullName}}%
}%
\def\LabelNum#1{%
	\global\edef#1{{\EltName}}%
}%
\long\def\Thm#1{%
	\Elt{Theorem}{#1}%
}%
\long\def\Prop#1{%
	\Elt{Proposition}{#1}%
}%
\long\def\Lem#1{%
	\Elt{Lemma}{#1}%
}%
\long\def\Def#1{%
	\Elt{Definition}{#1}%
}%
\long\def\Rem#1{%
	\Elt{Remark}{#1}%
}%
\long\def\Eg#1{%
	\Elt{Example}{#1}%
}%
\long\def\Cor#1{%
	\Elt{Corollary}{#1}%
}%
\def\qed{{%
	\unskip%
	\nobreak%
	\hfil%
	\penalty50%
	\hskip2em%
	\hbox{}%
	\nobreak%
	\hfil%
	{\bf QED}%
	\parfillskip=0pt%
	\finalhyphendemerits=0%
	\par%
}}%
\long\def\Proof#1{%
	{%
		\goodbreak%
		\par%
		\noindent%
		{\bf Proof:\ }%
		{\ignorespaces #1}%
		\qed%
	}%
	\goodbreak%
}%
%
%
\def\Fig#1{%
	{%
		\goodbreak%
		\par%
		\vskip\parskip%
		\def\EltName{%
			\number\SectNum.\number\EltNum%
		}%
		\def\EltFullName{%
			Figure\ \EltName%
		}%
		\hbox to \hsize {%
			\hfill%
			\vbox{%
				\halign{%
					\hfill%
					##%
					\hfill%
					\nobreak%
					\cr%
					#1%
					\cr%
					{\bf \EltFullName}%
					\cr%
				}%
			}%
			\hfill%
		}%
		\par%
	}%
	\advance\EltNum by 1%
}%

%
%
\font\tenfrak = eufm10%
\font\sevenfrak = eufm7%
\font\fivefrak = eufm5%
\newfam\frakfam%
\textfont\frakfam = \tenfrak%
\scriptfont\frakfam = \sevenfrak%
\scriptscriptfont\frakfam = \fivefrak%
\def\frak{\fam\frakfam\tenfrak}%
\font\tenrunes = frucms at 10pt%
\font\sevenrunes = frucms at 7pt%
\font\fiverunes = frucms at 5pt%
\newfam\runesfam%
\textfont\runesfam = \tenrunes%
\scriptfont\runesfam = \sevenrunes%
\scriptscriptfont\runesfam = \fiverunes%
\def\runes{\fam\runesfam\tenrunes}%
%
%
\def\C{{\bf C}}%
\def\R{{\bf R}}%
\def\Q{{\bf Q}}%
\def\Z{{\bf Z}}%
%
%
\def\defop#1#2{\def#1{\mathop{#2{}}\nolimits}}%
\def\id{{\rm id}}%
\def\d{{\rm d}}%
\defop\D{{\rm D}}%
\def\dbyd#1{{\d\over\d #1}}%
\def\derzero#1{\left.\dbyd{#1}\right|_{#1=0}}%
\def\pdbyd#1{{\partial\over\partial #1}}%
\def\pdbydo#1#2{{\partial^{#2}\over{\partial #1}^{#2}}}%
\def\conj#1{{\overline{#1}}}%
\def\norm#1{\left\|#1\right\|}%
\def\modulus#1{\left|#1\right|}%
\defop\eval{{\rm ev}}%
\def\Hodge{\mathord{\ast}}%
\def\hd{\partial}%
\defop\ahd{\conj{\partial}}%
\defop\ad{{\rm ad}}%
\defop\Ad{{\rm Ad}}%
\def\Dirac{{\frak D}}%
\def\morph{\mathop{\longrightarrow}\limits}%
\def\morphnamed#1{{%
	\setbox0 = \hbox{$\scriptstyle #1$}%
	\dimen0 = \wd0%
	\advance\dimen0 by 10pt%
	\mathop{%
		\;%
		\hbox to \dimen0 {\rightarrowfill}%
		\;%
	}\limits^{\box0}%
}}%
\def\converge#1{{%
	\setbox0 = \hbox{$\scriptstyle #1$}%
	\dimen0 = \wd0%
	\advance\dimen0 by 10pt%
	\mathop{%
		\;%
		\hbox to \dimen0 {\rightarrowfill}%
		\;%
	}\limits_{\box0}%
}}%
\def\isomorph{\morph^\sim}%
\defop\sign{{\rm sign}}%
\defop\index{{\rm ind}}%
\defop\And{\enskip\&\enskip}%
\defop\PD{{\rm PD}}%
\defop\LieDer{{\cal L}}%
\defop\cliff{{c\ell}}%
%
\def\CP{{\bf CP}}%
\def\CPBar{\overline{\bf CP}{}}
\def\SO{{\rm SO}}%
\def\U{{\rm U}}%
%
%
%
\defop\Hom{{\rm Hom}}%
\defop\Ker{{\rm Ker}}%
\defop\Image{{\rm Im}}%
\defop\Graph{{\rm Graph}}%
\defop\Coker{{\rm Coker}}%
\defop\Endo{{\rm End}}%
\defop\Sym{{\rm Sym}}%
\def\Homology{{\rm H}{}}%
\def\RedHomology{\tilde{\rm H}{}}%
\def\GeomCompl{{\rm c}}%
\def\BorelHomology{\GeomCompl\RedHomology}%
\def\HMFrom{\widehat{\rm HM}{}}%
\def\HMTilde{\widetilde{\rm HM}{}}%
\defop\ECH{{\rm ECH}}%
\defop\HFHat{{\widehat{\rm HF}}}%
\defop\HF{{\rm HF}}%
\defop\Harm{{\cal H}}%
\def\Sobolev#1,#2{\mathop{\rm L}\nolimits_{#1}^{#2}}%
\defop\LTwo{{\rm L}^2}%
\defop\LTwoLoc{{\rm L}^2_{\rm loc}}%
\def\DiffClass#1{\mathop{\rm C}\nolimits^{#1}}%
\def\Tan{{\rm T}}%
\def\Cotan{\Tan^\ast}%
\defop\Lie{{\rm Lie}}%
\defop\Cliff{{\rm C\ell}}%
\defop\Aut{{\rm Aut}}%
\def\Exterior{{\Lambda}}%
\defop\Ext{{\rm Ext}}%
\defop\Fred{{\rm Fred}}%
\defop\Stab{{\rm Stab}}%
\defop\Pin{{\rm Pin}}%
\def\SpinC{{{\rm Spin}^\C}}%
\def\Spinor{{\cal S}}%
\defop\Interior{{\rm int}}%
\defop\Obj{{\rm ob}}%
\def\Colim{\mathop{\rm Colim}}%
\defop\Ball{{\rm B}}%
\defop\supp{{\rm supp}}%
%
\defop\incl{\hookrightarrow}%
\def\epsilon{\varepsilon}%
\def\SetComp#1#2{\left\{\,#1\;\;\vrule\;\;#2\,\right\}}
\def\Half{{\textstyle{1\over2}}}
\def\Frac#1#2{{\textstyle{#1\over#2}}}
\def\Diagram#1{{%
	\def\Up{\Big\uparrow}%
	\def\Down{\Big\downarrow}%
	\def\Right{\longrightarrow}%
	\def\Left{\longleftarrow}%
	\advance\baselineskip by 1em%
	\vcenter{
	\halign{%
		\hfil$##$\hfil&&\quad\hfil$##$\hfil\cr%
		\mathstrut\cr%
		\noalign{\kern-\baselineskip}%
		#1\crcr%
		\mathstrut\cr%
		\noalign{\kern-\baselineskip}%
	}}%
}}%

\defop\SW{{\cal X}}%
\defop\FDSW{{\cal V}}%
\defop\SWF{{\rm SWF}}%
\defop\CSD{{\rm CSD}}%
\def\SpinCStruct{{\frak s}}%
\def\SpinStruct{{\cal S}}%
\def\Conf{{\cal C}}%
\def\Gauge{{\cal G}}%
\def\GaugeLie{{\frak G}}%
\def\GaugeLieAct{{\frak L}}%
\def\NormGauge{\Gauge^\circ}%
\def\NormGaugeLie{\GaugeLie^\circ}%
\defop\Inv{{\rm Inv}}%
\defop\Conley{{I}}%
\defop\KohnDirac{{\frak KD}}%
\def\LCSlice{{\cal K}}%
\def\ELCSlice{{\cal K}^{\rm E}}%
\def\AGCSlice{{\cal K}^{\rm AGC}}%
\def\ELCProj{\Pi^{\rm ELC}}%
\def\GCProj{\Pi^{\rm GC}}%
\def\LCProj{\Pi^{\rm LC}}%
\def\GCSW{\SW^{\rm GC}}%
\def\AGCProj{\Pi^{\rm AGC}}%
\def\GCGaugeOrbitSlice{{\cal J}}%
\def\AGCSW{\SW^{\rm AGC}}%
\def\FDAAGCSW#1{\SW^{{\rm AGC},#1}}%
\def\QuotBlowUpFDAAGCSW#1#2{%
	[\SW^{{\rm AGC},#1}_{#2}]^\sigma%
}%
\def\Connex{{\cal A}}%
\def\ContactConf{C}%
\def\ContactCirc{U}%
\def\ContactNormal#1#2{E^{#1}_{#2}}%
\def\ContactNormalSt#1#2{E^{{\rm s},#1}_{#2}}%
\def\ContactNormalUnst#1#2{E^{{\rm u},#1}_{#2}}%
\def\FDASW{\SW}%
\def\FDASWFlow{\varphi}%
\def\Sphere{{\rm S}}%
\def\Disk{{\rm D}}%
\def\Vol{{\rm Vol}}%
\def\ContactNVBF{{\runes e}}%
\def\ContactFNVBF{{\runes f}}%
\def\Thom{\Theta}%
\def\ContactThom{{\cal T}}%
\def\RO{{\rm RO}}%
\def\B{{\rm B}}%
\def\E{{\rm E}}%
\def\Spectra#1#2{{#1{\cal S}#2}}%
\def\HoSpectra#1#2{{{h}#1{\cal S}#2}}%
\def\StabHoSpectra#1#2{{\bar{h}#1{\cal S}#2}}%
\defop\SignOp{{\rm Sign}}%
\def\Universe{{\cal U}}
\def\AnotherUniverse{{\cal V}}
\def\CoulombUniverse{{\cal W}}
\def\MetricDesuspParamRes{{\frak n}}%
\def\ThomClass{{\theta}}%
\def\BorelHMTilde{\GeomCompl\HMTilde}%
\def\HantzscheWendt{Y_{\rm HW}}%
\defop\TB{{\rm tb}}%
\def\Const{K}%
\def\Perturbation{{\frak q}}%

\global\def\AtiyahPatodiSingerSpectralCited{0}%
\def\AtiyahPatodiSingerSpectralAlt{%
	\global\def\AtiyahPatodiSingerSpectralCited{1}%
	Atiyah,\ Patodi\ \&\ Singer\ 1975%
}%
\global\def\CauchyCoursCited{0}%
\def\CauchyCours{%
	\global\def\CauchyCoursCited{1}%
	Cauchy\ (1821)%
}%
\global\def\ChelnokovMednykhHantzscheWendtCited{0}%
\def\ChelnokovMednykhHantzscheWendtAlt{%
	\global\def\ChelnokovMednykhHantzscheWendtCited{1}%
	Chelnokov\ \&\ Mednykh\ 2020%
}%
\global\def\ColinGhigginiHondaEquivalenceICited{0}%
\def\ColinGhigginiHondaEquivalenceIAuthors{%
	\global\def\ColinGhigginiHondaEquivalenceICited{1}%
	Colin,\ Ghiggini\ \&\ Honda%
}%
\def\ColinGhigginiHondaEquivalenceIYear{%
	\global\def\ColinGhigginiHondaEquivalenceICited{1}%
	2012a%
}%
\global\def\ColinGhigginiHondaEquivalenceIICited{0}%
\def\ColinGhigginiHondaEquivalenceIIYear{%
	\global\def\ColinGhigginiHondaEquivalenceIICited{1}%
	2012b%
}%
\global\def\ConleyIsolatedCited{0}%
\def\ConleyIsolated{%
	\global\def\ConleyIsolatedCited{1}%
	Conley\ (1978)%
}%
\def\ConleyIsolatedAlt{%
	\global\def\ConleyIsolatedCited{1}%
	Conley\ 1978%
}%
\global\def\CorneaHomotopicalCited{0}%
\global\def\CostenobleWanerEquivariantPoincareCited{0}%
\global\def\CostenobleWanerEquivariantOrdinaryCited{0}%
\global\def\CristofaroGardinerAbsoluteCited{0}%
\def\CristofaroGardinerAbsoluteAlt{%
	\global\def\CristofaroGardinerAbsoluteCited{1}%
	Cristofaro-Gardiner\ 2013%
}%
\global\def\DonnellyEtaCited{0}%
\global\def\EcheverriaNaturalityCited{0}%
\global\def\FloerRefinementCited{0}%
\def\FloerRefinement{%
	\global\def\FloerRefinementCited{1}%
	Floer\ (1987)%
}%
\def\FloerRefinementAlt{%
	\global\def\FloerRefinementCited{1}%
	Floer\ 1987%
}%
\global\def\FloerWittensComplexCited{0}%
\global\def\FloerZehnderEquivariantConleyCited{0}%
\def\FloerZehnderEquivariantConleyAlt{%
	\global\def\FloerZehnderEquivariantConleyCited{1}%
	Floer\ \&\ Zehnder\ 1988%
}%
\global\def\GeigesIntroductionCited{0}%
\def\GeigesIntroduction{%
	\global\def\GeigesIntroductionCited{1}%
	Geiges\ (2008)%
}%
\global\def\GhigginiOzsvathSzaboCited{0}%
\global\def\GhigginiTightCited{0}%
\def\GhigginiTightAlt{%
	\global\def\GhigginiTightCited{1}%
	Ghiggini\ 2008%
}%
\global\def\GhigginiLiscaStipsiczClassificationCited{0}%
\def\GhigginiLiscaStipsiczClassificationAlt{%
	\global\def\GhigginiLiscaStipsiczClassificationCited{1}%
	Ghiggini,\ Lisca\ \&\ Stipsicz\ 2006%
}%
\global\def\GhigginiHornMorrisTightCited{0}%
\global\def\GirouxStructuresCited{0}%
\def\GirouxStructures{%
	\global\def\GirouxStructuresCited{1}%
	Giroux\ (2000)%
}%
\def\GirouxStructuresAlt{%
	\global\def\GirouxStructuresCited{1}%
	Giroux\ 2000%
}%
\global\def\GompfHandlebodyCited{0}%
\def\GompfHandlebody{%
	\global\def\GompfHandlebodyCited{1}%
	Gompf\ (1998)%
}%
\def\GompfHandlebodyAlt{%
	\global\def\GompfHandlebodyCited{1}%
	Gompf\ 1998%
}%
\global\def\GompfStipsiczFourManifoldsCited{0}%
\def\GompfStipsiczFourManifolds{%
	\global\def\GompfStipsiczFourManifoldsCited{1}%
	Gompf\ \&\ Stipsicz\ (1999)%
}%
\def\GompfStipsiczFourManifoldsAlt{%
	\global\def\GompfStipsiczFourManifoldsCited{1}%
	Gompf\ \&\ Stipsicz\ 1999%
}%
\global\def\GreenleesMayGeneralizedCited{0}%
\global\def\HellConleyCited{0}%
\global\def\HondaClassificationCited{0}%
\def\HondaClassification{%
	\global\def\HondaClassificationCited{1}%
	Honda\ (1999)%
}%
\def\HondaClassificationAlt{%
	\global\def\HondaClassificationCited{1}%
	Honda\ 1999%
}%
\global\def\HondaClassificationIICited{0}%
\global\def\HutchingsTaubesGluingCited{0}%
\def\HutchingsTaubesGluingAlt{%
	\global\def\HutchingsTaubesGluingCited{1}%
	Hutchings\ \&\ Taubes\ 2007%
}%
\global\def\IidaTaniguchiSeibergWittenCited{0}%
\def\IidaTaniguchiSeibergWitten{%
	\global\def\IidaTaniguchiSeibergWittenCited{1}%
	Iida\ \&\ Taniguchi\ (2021)%
}%
\global\def\KhandhawitNewGaugeCited{0}%
\global\def\KhandhawitLinSasahiraUnfoldedIICited{0}%
\global\def\KhuzamLiftingCited{0}%
\def\KhuzamLifting{%
	\global\def\KhuzamLiftingCited{1}%
	Khuzam\ (2012)%
}%
\def\KhuzamLiftingAlt{%
	\global\def\KhuzamLiftingCited{1}%
	Khuzam\ 2012%
}%
\global\def\KronheimerMrowkaMonopolesAndContactStructuresCited{0}%
\def\KronheimerMrowkaMonopolesAndContactStructures{%
	\global\def\KronheimerMrowkaMonopolesAndContactStructuresCited{1}%
	Kronheimer\ \&\ Mrowka\ (1997)%
}%
\global\def\KronheimerMrowkaMonopolesAndThreeManifoldsCited{0}%
\def\KronheimerMrowkaMonopolesAndThreeManifolds{%
	\global\def\KronheimerMrowkaMonopolesAndThreeManifoldsCited{1}%
	Kronheimer\ \&\ Mrowka\ (2007)%
}%
\def\KronheimerMrowkaMonopolesAndThreeManifoldsAlt{%
	\global\def\KronheimerMrowkaMonopolesAndThreeManifoldsCited{1}%
	Kronheimer\ \&\ Mrowka\ 2007%
}%
\global\def\KronheimerMrowkaOzsvathSzaboMonopolesCited{0}%
\def\KronheimerMrowkaOzsvathSzaboMonopoles{%
	\global\def\KronheimerMrowkaOzsvathSzaboMonopolesCited{1}%
	Kronheimer,\ Mrowka,\ Ozsv\'ath\ \&\ Szab\'o\ (2007)%
}%
\def\KronheimerMrowkaOzsvathSzaboMonopolesAlt{%
	\global\def\KronheimerMrowkaOzsvathSzaboMonopolesCited{1}%
	Kronheimer,\ Mrowka,\ Ozsv\'ath\ \&\ Szab\'o\ 2007%
}%
\global\def\KurlandHomotopyCited{0}%
\def\KurlandHomotopyAlt{%
	\global\def\KurlandHomotopyCited{1}%
	Kurland\ 1982%
}%
\global\def\KutluhanLeeTaubesHFEqHMICited{0}%
\def\KutluhanLeeTaubesHFEqHMIAuthors{%
	\global\def\KutluhanLeeTaubesHFEqHMICited{1}%
	Kutluhan,\ Lee\ \&\ Taubes%
}%
\def\KutluhanLeeTaubesHFEqHMIYear{%
	\global\def\KutluhanLeeTaubesHFEqHMICited{1}%
	2020a%
}%
\global\def\KutluhanLeeTaubesHFEqHMIICited{0}%
\def\KutluhanLeeTaubesHFEqHMIIYear{%
	\global\def\KutluhanLeeTaubesHFEqHMIICited{1}%
	2020b%
}%
\global\def\KutluhanLeeTaubesHFEqHMIIICited{0}%
\def\KutluhanLeeTaubesHFEqHMIIIYear{%
	\global\def\KutluhanLeeTaubesHFEqHMIIICited{1}%
	2020c%
}%
\global\def\KutluhanLeeTaubesHFEqHMIVCited{0}%
\def\KutluhanLeeTaubesHFEqHMIVYear{%
	\global\def\KutluhanLeeTaubesHFEqHMIVCited{1}%
	2021a%
}%
\global\def\KutluhanLeeTaubesHFEqHMVCited{0}%
\def\KutluhanLeeTaubesHFEqHMVYear{%
	\global\def\KutluhanLeeTaubesHFEqHMVCited{1}%
	2021b%
}%
\global\def\LewisMaySteinbergerEquivariantCited{0}%
\def\LewisMaySteinbergerEquivariant{%
	\global\def\LewisMaySteinbergerEquivariantCited{1}%
	Lewis Jr.,\ May\ \&\ Steinberger\ (1986)%
}%
\def\LewisMaySteinbergerEquivariantAlt{%
	\global\def\LewisMaySteinbergerEquivariantCited{1}%
	Lewis Jr.,\ May\ \&\ Steinberger\ 1986%
}%
\global\def\LidmanManolescuEquivalenceCited{0}%
\def\LidmanManolescuEquivalence{%
	\global\def\LidmanManolescuEquivalenceCited{1}%
	Lidman\ \&\ Manolescu\ (2018a)%
}%
\def\LidmanManolescuEquivalenceAlt{%
	\global\def\LidmanManolescuEquivalenceCited{1}%
	Lidman\ \&\ Manolescu\ 2018a%
}%
\global\def\LidmanManolescuCoveringCited{0}%
\def\LidmanManolescuCovering{%
	\global\def\LidmanManolescuCoveringCited{1}%
	Lidman\ \&\ Manolescu\ (2018b)%
}%
\def\LidmanManolescuCoveringAlt{%
	\global\def\LidmanManolescuCoveringCited{1}%
	Lidman\ \&\ Manolescu\ 2018b%
}%
\global\def\LinSolvCited{0}%
\def\LinSolv{%
	\global\def\LinSolvCited{1}%
	Lin\ (2020)%
}%
\global\def\LinLipnowskiHyperbolicCited{0}%
\def\LinLipnowskiHyperbolic{%
	\global\def\LinLipnowskiHyperbolicCited{1}%
	Lin\ \&\ Lipnowski\ (2022)%
}%
\global\def\LinRubermanSavelievFroyshovCited{0}%
\def\LinRubermanSavelievFroyshovAlt{%
	\global\def\LinRubermanSavelievFroyshovCited{1}%
	Lin,\ Ruberman\ \&\ Saveliev\ 2018%
}%
\global\def\ManolescuSeibergWittenFloerCited{0}%
\def\ManolescuSeibergWittenFloer{%
	\global\def\ManolescuSeibergWittenFloerCited{1}%
	Manolescu\ (2003)%
}%
\def\ManolescuSeibergWittenFloerAlt{%
	\global\def\ManolescuSeibergWittenFloerCited{1}%
	Manolescu\ 2003%
}%
\global\def\ManolescuGluingCited{0}%
\global\def\MatkovicClassificationCited{0}%
\def\MatkovicClassificationAlt{%
	\global\def\MatkovicClassificationCited{1}%
	Matkovi\v c\ 2018%
}%
\global\def\MayEtAlEquivariantCited{0}%
\def\MayEtAlEquivariant{%
	\global\def\MayEtAlEquivariantCited{1}%
	May\ \&\ al.\ (1996)%
}%
\def\MayEtAlEquivariantAlt{%
	\global\def\MayEtAlEquivariantCited{1}%
	May\ \&\ al.\ 1996%
}%
\global\def\McCordMappingsAndHomologicalCited{0}%
\global\def\McCordMappingsAndMorseCited{0}%
\global\def\PalisMeloGeometricCited{0}%
\global\def\MischaikowConleyCited{0}%
\def\MischaikowConley{%
	\global\def\MischaikowConleyCited{1}%
	Mischaikow\ (1995)%
}%
\global\def\MoserElementaryCited{0}%
\def\MoserElementary{%
	\global\def\MoserElementaryCited{1}%
	Moser\ (1971)%
}%
\global\def\NicolaescuNotesCited{0}%
\def\NicolaescuNotesAlt{%
	\global\def\NicolaescuNotesCited{1}%
	Nicolaescu\ 2000%
}%
\global\def\OzsvathSzaboAbsolutelyCited{0}%
\def\OzsvathSzaboAbsolutely{%
	\global\def\OzsvathSzaboAbsolutelyCited{1}%
	Ozsv\'ath\ \&\ Szab\'o\ (2003)%
}%
\global\def\OzsvathSzaboHolomorphicCited{0}%
\def\OzsvathSzaboHolomorphic{%
	\global\def\OzsvathSzaboHolomorphicCited{1}%
	Ozsv\'ath\ \&\ Szab\'o\ (2004)%
}%
\global\def\OzsvathSzaboContactCited{0}%
\def\OzsvathSzaboContact{%
	\global\def\OzsvathSzaboContactCited{1}%
	Ozsv\'ath\ \&\ Szab\'o\ (2005)%
}%
\global\def\PalisOnMorseSmaleCited{0}%
\global\def\PetitSpinCCited{0}%
\def\PetitSpinC{%
	\global\def\PetitSpinCCited{1}%
	Petit\ (2005)%
}%
\global\def\RamosAbsoluteCited{0}%
\def\RamosAbsoluteAlt{%
	\global\def\RamosAbsoluteCited{1}%
	Ramos\ 2018%
}%
\global\def\RolfsenKnotsCited{0}%
\def\RolfsenKnots{%
	\global\def\RolfsenKnotsCited{1}%
	Rolfsen\ (2003)%
}%
\global\def\RosoThesisCited{0}%
\def\RosoThesis{%
	\global\def\RosoThesisCited{1}%
	Roso\ (2023)%
}%
\def\RosoThesisAlt{%
	\global\def\RosoThesisCited{1}%
	Roso\ 2023%
}%
\global\def\ShubinPseudodifferentialCited{0}%
\global\def\TaubesSeibergWittenSymplecticCited{0}%
\def\TaubesSeibergWittenSymplectic{%
	\global\def\TaubesSeibergWittenSymplecticCited{1}%
	Taubes\ (1994)%
}%
\global\def\TaubesSWToGrCited{0}%
\def\TaubesSWToGr{%
	\global\def\TaubesSWToGrCited{1}%
	Taubes\ (1996)%
}%
\global\def\TaubesWeinsteinICited{0}%
\def\TaubesWeinsteinI{%
	\global\def\TaubesWeinsteinICited{1}%
	Taubes\ (2007)%
}%
\def\TaubesWeinsteinIAlt{%
	\global\def\TaubesWeinsteinICited{1}%
	Taubes\ 2007%
}%
\global\def\TaubesWeinsteinIICited{0}%
\def\TaubesWeinsteinII{%
	\global\def\TaubesWeinsteinIICited{1}%
	Taubes\ (2009)%
}%
\def\TaubesWeinsteinIIAlt{%
	\global\def\TaubesWeinsteinIICited{1}%
	Taubes\ 2009%
}%
\global\def\TaubesEmbeddedICited{0}%
\def\TaubesEmbeddedIAuthors{%
	\global\def\TaubesEmbeddedICited{1}%
	Taubes%
}%
\def\TaubesEmbeddedIYear{%
	\global\def\TaubesEmbeddedICited{1}%
	2010a%
}%
\global\def\TaubesEmbeddedIICited{0}%
\def\TaubesEmbeddedIIYear{%
	\global\def\TaubesEmbeddedIICited{1}%
	2010b%
}%
\global\def\TaubesEmbeddedIIICited{0}%
\def\TaubesEmbeddedIIIYear{%
	\global\def\TaubesEmbeddedIIICited{1}%
	2010c%
}%
\global\def\TaubesEmbeddedIVCited{0}%
\def\TaubesEmbeddedIVYear{%
	\global\def\TaubesEmbeddedIVCited{1}%
	2010d%
}%
\global\def\TaubesEmbeddedVCited{0}%
\def\TaubesEmbeddedV{%
	\global\def\TaubesEmbeddedVCited{1}%
	Taubes\ (2010e)%
}%
\def\TaubesEmbeddedVYear{%
	\global\def\TaubesEmbeddedVCited{1}%
	2010e%
}%
\global\def\TomDieckTransformationCited{0}%
\def\TomDieckTransformationAlt{%
	\global\def\TomDieckTransformationCited{1}%
	tom Dieck\ 1987%
}%
\global\def\WuLegendrianCited{0}%
\def\WuLegendrianAlt{%
	\global\def\WuLegendrianCited{1}%
	Wu\ 2006%
}%
\def\Bib{{%
	\parindent=0cm%
	\if\AtiyahPatodiSingerSpectralCited 1%
		{\bf Atiyah,\ M.F.,\ Patodi,\ V.K.\ \&\ Singer,\ I.M.}
		(1975).
		``Spectral Asymmetry and Riemannian Geometry I.''
		{\it Mathematical Proceedings of the Cambridge Philosophical Society},
		77(1),
		pp. 43--69.
		Cambridge University Press.
		\par%
	\fi%
	\if\CauchyCoursCited 1%
		{\bf Cauchy,\ A.L.}
		(1821).
		{\it Cours d'Analyse de l'\'Ecole Royale Polytechnique}.
		de Bure, Paris.
		\par%
	\fi%
	\if\ChelnokovMednykhHantzscheWendtCited 1%
		{\bf Chelnokov,\ G.\ \&\ Mednykh,\ A.}
		(2020).
		{\it On the coverings of Hantzsche-Wendt manifold}.
		arXiv preprint arXiv:2009.06691.
		\par%
	\fi%
	\if\ColinGhigginiHondaEquivalenceICited 1%
		{\bf Colin,\ V.,\ Ghiggini,\ P.\ \&\ Honda,\ K.}
		(2012a).
		{\it The equivalence of Heegaard Floer and embedded contact homology via open book decompositions I}.
		arXiv preprint arXiv:1208.1074.
		\par%
	\fi%
	\if\ColinGhigginiHondaEquivalenceIICited 1%
		{\bf Colin,\ V.,\ Ghiggini,\ P.\ \&\ Honda,\ K.}
		(2012b).
		{\it The equivalence of Heegaard Floer and embedded contact homology via open book decompositions II}.
		arXiv preprint arXiv:1208.1077.
		\par%
	\fi%
	\if\ConleyIsolatedCited 1%
		{\bf Conley,\ C.}
		(1978).
		{\it Isolated Invariant Sets and the Morse Index},
		(38).
		American Mathematical Society.
		\par%
	\fi%
	\if\CorneaHomotopicalCited 1%
		{\bf Cornea,\ O.}
		(2000).
		``Homotopical dynamics: suspension and duality.''
		{\it Ergodic Theory and Dynamical Systems},
		20(2),
		pp. 379--391.
		Cambridge University Press.
		\par%
	\fi%
	\if\CostenobleWanerEquivariantPoincareCited 1%
		{\bf Costenoble,\ S.R.\ \&\ Waner,\ S.}
		(1992).
		``Equivariant Poincar\'e Duality.''
		{\it Michigan Mathematical Journal},
		39(2),
		pp. 325--351.
		\par%
	\fi%
	\if\CostenobleWanerEquivariantOrdinaryCited 1%
		{\bf Costenoble,\ S.R.\ \&\ Waner,\ S.}
		(2016).
		{\it Equivariant ordinary homology and cohomology}.
		Springer.
		\par%
	\fi%
	\if\CristofaroGardinerAbsoluteCited 1%
		{\bf Cristofaro-Gardiner,\ D.}
		(2013).
		``The absolute gradings on embedded contact homology and Seiberg--Witten Floer cohomology.''
		{\it Algebraic \& Geometric Topology},
		13(4),
		pp. 2239--2260.
		Mathematical Sciences Publishers.
		\par%
	\fi%
	\if\DonnellyEtaCited 1%
		{\bf Donnelly,\ H.}
		(1978).
		``Eta invariants for $G$-spaces.''
		{\it Indiana University Mathematics Journal},
		27(6),
		pp. 889--918.
		\par%
	\fi%
	\if\EcheverriaNaturalityCited 1%
		{\bf Echeverria,\ M.}
		(2020).
		``Naturality of the contact invariant in monopole Floer homology under strong symplectic cobordisms.''
		{\it Algebraic \& Geometric Topology},
		20(4).
		\par%
	\fi%
	\if\FloerRefinementCited 1%
		{\bf Floer,\ A.}
		(1987).
		``A Refinement of the Conley Index and an Application to the Stability of Hyperbolic Invariant Sets.''
		{\it Ergodic Theory and Dynamical Systems},
		7(1).
		\par%
	\fi%
	\if\FloerWittensComplexCited 1%
		{\bf Floer,\ A.}
		(1989).
		``Witten's complex and infinite-dimensional Morse theory.''
		{\it Journal of differential geometry},
		30(1),
		pp. 207--221.
		Lehigh University.
		\par%
	\fi%
	\if\FloerZehnderEquivariantConleyCited 1%
		{\bf Floer,\ A.\ \&\ Zehnder,\ E.}
		(1988).
		``The equivariant Conley index and bifurcations of periodic solutions of Hamiltonian systems.''
		{\it Ergodic Theory and Dynamical Systems},
		8,
		pp. 87--97.
		Cambridge University Press.
		\par%
	\fi%
	\if\GeigesIntroductionCited 1%
		{\bf Geiges,\ H.}
		(2008).
		{\it An introduction to contact topology}.
		Cambridge University Press.
		\par%
	\fi%
	\if\GhigginiOzsvathSzaboCited 1%
		{\bf Ghiggini,\ P.}
		(2006).
		``Ozsv\'ath-Szab\'o invariants and fillability of contact structures.''
		{\it Mathematische Zeitschrift},
		253(1).
		\par%
	\fi%
	\if\GhigginiTightCited 1%
		{\bf Ghiggini,\ P.}
		(2008).
		``On tight contact structures with negative maximal twisting number on small Seifert manifolds.''
		{\it Algebraic \& Geometric Topology},
		8(1),
		pp. 381--396.
		\par%
	\fi%
	\if\GhigginiLiscaStipsiczClassificationCited 1%
		{\bf Ghiggini,\ P.,\ Lisca,\ P.\ \&\ Stipsicz,\ A.I.}
		(2006).
		``Classification of Tight Contact Structures on Small Seifert 3-Manifolds with $e_0\geq0$.''
		{\it Proceedings of the American Mathematical Society},
		pp. 909--916.
		\par%
	\fi%
	\if\GhigginiHornMorrisTightCited 1%
		{\bf Ghiggini,\ P.\ \&\ Van Horn-Morris,\ J.}
		(2016).
		``Tight contact structures on the Brieskorn spheres $-\Sigma(2,3,6n-1)$ and contact invariants.''
		{\it Journal f\"ur die reine und angewandte Mathematik},
		2016(718),
		pp. 1--24.
		\par%
	\fi%
	\if\GirouxStructuresCited 1%
		{\bf Giroux,\ E.}
		(2000).
		``Structures de contact en dimension trois et bifurcations des feuilletages de surfaces.''
		{\it Inventiones Mathematicae},
		141(3),
		pp. 615--689.
		\par%
	\fi%
	\if\GompfHandlebodyCited 1%
		{\bf Gompf,\ R.E.}
		(1998).
		``Handlebody construction of Stein surfaces.''
		{\it Annals of mathematics},
		pp. 619--693.
		\par%
	\fi%
	\if\GompfStipsiczFourManifoldsCited 1%
		{\bf Gompf,\ R.E.\ \&\ Stipsicz,\ A.I.}
		(1999).
		{\it 4-manifolds and Kirby calculus},
		(20).
		American Mathematical Soc..
		\par%
	\fi%
	\if\GreenleesMayGeneralizedCited 1%
		{\bf Greenlees,\ J.P.C.\ \&\ May,\ J.P.}
		(1995).
		{\it Generalized tate cohomology},
		543.
		American Mathematical Soc..
		\par%
	\fi%
	\if\HellConleyCited 1%
		{\bf Hell,\ J.}
		(2009).
		{\it Conley Index at Infinity}.
		Ph.D. Dissertation, Freie Universität Berlin.
		\par%
	\fi%
	\if\HondaClassificationCited 1%
		{\bf Honda,\ K.}
		(1999).
		``On the classification of tight contact structures I.''
		{\it Geometry \& Topology},
		4,
		pp. 309--368.
		\par%
	\fi%
	\if\HondaClassificationIICited 1%
		{\bf Honda,\ K.}
		(2000).
		``On the classification of tight contact structures II.''
		{\it Journal of Differential Geometry},
		55(1),
		pp. 83--143.
		Lehigh University.
		\par%
	\fi%
	\if\HutchingsTaubesGluingCited 1%
		{\bf Hutchings,\ M.\ \&\ Taubes,\ C.H.}
		(2007).
		``Gluing pseudoholomorphic curves along branched covered cylinders I.''
		{\it Journal of Symplectic Geometry},
		5(1),
		pp. 43--137.
		International Press of Boston.
		\par%
	\fi%
	\if\IidaTaniguchiSeibergWittenCited 1%
		{\bf Iida,\ N.\ \&\ Taniguchi,\ M.}
		(2021).
		``Seiberg-Witten Floer Homotopy Contact Invariant.''
		{\it Studia Scientiarum Mathematicarum Hungarica},
		58(4),
		pp. 505--558.
		Akad\'emiai Kiad\'o Budapest.
		\par%
	\fi%
	\if\KhandhawitNewGaugeCited 1%
		{\bf Khandhawit,\ T.}
		(2015).
		``A new gauge slice for the relative Bauer--Furuta invariants.''
		{\it Geometry \& Topology},
		19(3),
		pp. 1631--1655.
		Mathematical Sciences Publishers.
		\par%
	\fi%
	\if\KhandhawitLinSasahiraUnfoldedIICited 1%
		{\bf Khandhawit,\ T.,\ Lin,\ J.\ \&\ Sasahira,\ H.}
		(2018).
		{\it Unfolded Seiberg-Witten Floer spectra, II: Relative invariants and the gluing theorem}.
		arXiv preprint arXiv:1809.09151.
		\par%
	\fi%
	\if\KhuzamLiftingCited 1%
		{\bf Khuzam,\ M.B.}
		(2012).
		``Lifting the $3$-dimensional invariant of $2$-plane fields on $3$-manifolds.''
		{\it Topology and its Applications},
		159(3),
		pp. 704--710.
		\par%
	\fi%
	\if\KronheimerMrowkaMonopolesAndContactStructuresCited 1%
		{\bf Kronheimer,\ P.\ \&\ Mrowka,\ T.}
		(1997).
		``Monopoles and Contact Structures.''
		{\it Inventiones Mathematicae},
		130(2),
		pp. 209--255.
		\par%
	\fi%
	\if\KronheimerMrowkaMonopolesAndThreeManifoldsCited 1%
		{\bf Kronheimer,\ P.\ \&\ Mrowka,\ T.}
		(2007).
		{\it Monopoles and three-manifolds},
		10.
		Cambridge University Press.
		\par%
	\fi%
	\if\KronheimerMrowkaOzsvathSzaboMonopolesCited 1%
		{\bf Kronheimer,\ P.,\ Mrowka,\ T.,\ Ozsv\'ath,\ P.\ \&\ Szab\'o,\ Z.}
		(2007).
		``Monopoles and Lens Space Surgeries.''
		{\it Annals of Mathematics},
		pp. 457--546.
		JSTOR.
		\par%
	\fi%
	\if\KurlandHomotopyCited 1%
		{\bf Kurland,\ H.L.}
		(1982).
		``Homotopy Invariants of Repeller-Attractor Pairs. I. The P\"uppe Sequence of an RA Pair.''
		{\it Journal of Differential Equations},
		46(1),
		pp. 1--31.
		\par%
	\fi%
	\if\KutluhanLeeTaubesHFEqHMICited 1%
		{\bf Kutluhan,\ \c C.,\ Lee,\ Y.J.\ \&\ Taubes,\ C.H.}
		(2020a).
		``HF=HM, I: Heegaard Floer homology and Seiberg--Witten Floer homology.''
		{\it Geometry \& Topology},
		24(6),
		pp. 2829--2854.
		Mathematical Sciences Publishers.
		\par%
	\fi%
	\if\KutluhanLeeTaubesHFEqHMIICited 1%
		{\bf Kutluhan,\ \c C.,\ Lee,\ Y.J.\ \&\ Taubes,\ C.H.}
		(2020b).
		``HF=HM, II: Reeb orbits and holomorphic curves for the ech/Heegaard Floer correspondence.''
		{\it Geometry \& Topology},
		24(6),
		pp. 2855--3012.
		Mathematical Sciences Publishers.
		\par%
	\fi%
	\if\KutluhanLeeTaubesHFEqHMIIICited 1%
		{\bf Kutluhan,\ \c C.,\ Lee,\ Y.J.\ \&\ Taubes,\ C.H.}
		(2020c).
		``HF=HM, III: holomorphic curves and the differential for the ech/Heegaard Floer correspondenc.''
		{\it Geometry \& Topology},
		24(6),
		pp. 3013--3218.
		Mathematical Sciences Publishers.
		\par%
	\fi%
	\if\KutluhanLeeTaubesHFEqHMIVCited 1%
		{\bf Kutluhan,\ \c C.,\ Lee,\ Y.J.\ \&\ Taubes,\ C.H.}
		(2021a).
		``HF=HM, IV: The Seiberg--Witten Floer homology and ech correspondence.''
		{\it Geometry \& Topology},
		24(7),
		pp. 3219--3469.
		Mathematical Sciences Publishers.
		\par%
	\fi%
	\if\KutluhanLeeTaubesHFEqHMVCited 1%
		{\bf Kutluhan,\ \c C.,\ Lee,\ Y.J.\ \&\ Taubes,\ C.H.}
		(2021b).
		``HF=HM, V: Seiberg--Witten Floer homology and handle additions.''
		{\it Geometry \& Topology},
		24(7),
		pp. 3471--3748.
		Mathematical Sciences Publishers.
		\par%
	\fi%
	\if\LewisMaySteinbergerEquivariantCited 1%
		{\bf Lewis Jr.,\ L.G.,\ May,\ J.P.\ \&\ Steinberger,\ M.}
		(1986).
		``Equivariant stable homotopy theory.''
		{\it Lecture notes in mathematics},
		1213.
		Springer-Verlag.
		\par%
	\fi%
	\if\LidmanManolescuEquivalenceCited 1%
		{\bf Lidman,\ T.\ \&\ Manolescu,\ C.}
		(2018a).
		``The equivalence of two Seiberg-Witten Floer homologies.''
		{\it As\-t\'e\-risque},
		399.
		Soci\'et\'e Math\'ematique de France.
		\par%
	\fi%
	\if\LidmanManolescuCoveringCited 1%
		{\bf Lidman,\ T.\ \&\ Manolescu,\ C.}
		(2018b).
		``Floer homology and covering spaces.''
		{\it Geometry \& Topology},
		22(5),
		pp. 2817--2838.
		\par%
	\fi%
	\if\LinSolvCited 1%
		{\bf Lin,\ F.}
		(2020).
		``Monopole Floer homology and SOLV geometry.''
		{\it Annales Henri Lebesgue},
		3,
		pp. 1117--1131.
		\par%
	\fi%
	\if\LinLipnowskiHyperbolicCited 1%
		{\bf Lin,\ F.\ \&\ Lipnowski,\ M.}
		(2022).
		``The Seiberg-Witten equations and the length spectrum of hyperbolic three-manifolds.''
		{\it Journal of the American Mathematical Society},
		35(1),
		pp. 233--293.
		\par%
	\fi%
	\if\LinRubermanSavelievFroyshovCited 1%
		{\bf Lin,\ J.,\ Ruberman,\ D.\ \&\ Saveliev,\ N.}
		(2018).
		{\it On the Fr\o yshov invariant and monopole Lefschetz number}.
		arXiv preprint arXiv:1802.07704.
		\par%
	\fi%
	\if\ManolescuSeibergWittenFloerCited 1%
		{\bf Manolescu,\ C.}
		(2003).
		``Seiberg-Witten-Floer stable homotopy type of three-manifolds with $b_1=0$.''
		{\it Geometry \& Topology},
		7(2),
		pp. 889--932.
		Mathematical Sciences Publishers.
		\par%
	\fi%
	\if\ManolescuGluingCited 1%
		{\bf Manolescu,\ C.}
		(2007).
		``A gluing theorem for the relative Bauer-Furuta invariants.''
		{\it Journal of Differential Geometry},
		76(1),
		pp. 117--153.
		\par%
	\fi%
	\if\MatkovicClassificationCited 1%
		{\bf Matkovi\v c,\ I.}
		(2018).
		``Classification of tight contact structures on small Seifert fibered L–spaces.''
		{\it Algebraic \& Geometric Topology},
		18(1),
		pp. 111--152.
		\par%
	\fi%
	\if\MayEtAlEquivariantCited 1%
		{\bf May,\ J.P.\ \&\ al.}
		(1996).
		{\it Equivariant Homotopy and Cohomology Theory},
		(91).
		American Mathematical Soc..
		\par%
	\fi%
	\if\McCordMappingsAndHomologicalCited 1%
		{\bf McCord,\ C.K.}
		(1986).
		{\it Mappings and homological properties in the Conley index theory}.
		Ph.D. Dissertation, University of Wisconsin--Madison.
		\par%
	\fi%
	\if\McCordMappingsAndMorseCited 1%
		{\bf McCord,\ C.K.}
		(1991).
		``Mappings and Morse decompositions in the Conley index theory.''
		{\it Indiana University Mathematics Journal},
		pp. 1061--1082.
		JSTOR.
		\par%
	\fi%
	\if\PalisMeloGeometricCited 1%
		{\bf Melo,\ W.\ \&\ Palis,\ J.}
		(1982).
		{\it Geometric Theory of Dynamical Systems: An Introduction},
		(1).
		Springer-Verlag NewYork Inc.
		\par%
	\fi%
	\if\MischaikowConleyCited 1%
		{\bf Mischaikow,\ K.}
		(1995).
		``Conley index theory.''
		{\it Dynamical Systems},
		pp. 119--207.
		Springer.
		\par%
	\fi%
	\if\MoserElementaryCited 1%
		{\bf Moser,\ L.}
		(1971).
		``Elementary surgery along a torus knot.''
		{\it Pacific Journal of Mathematics},
		38(3),
		pp. 737--745.
		\par%
	\fi%
	\if\NicolaescuNotesCited 1%
		{\bf Nicolaescu,\ L.I.}
		(2000).
		{\it Notes on Seiberg-Witten Theory},
		28.
		American Mathematical Soc..
		\par%
	\fi%
	\if\OzsvathSzaboAbsolutelyCited 1%
		{\bf Ozsv\'ath,\ P.\ \&\ Szab\'o,\ Z.}
		(2003).
		``Absolutely graded Floer homologies and intersection forms for four-manifolds with boundary.''
		{\it Advances in Mathematics},
		173(2),
		pp. 179--261.
		Elsevier.
		\par%
	\fi%
	\if\OzsvathSzaboHolomorphicCited 1%
		{\bf Ozsv\'ath,\ P.\ \&\ Szab\'o,\ Z.}
		(2004).
		``Holomorphic disks and topological invariants for closed three-manifolds.''
		{\it Annals of Mathematics},
		pp. 1027--1158.
		JSTOR.
		\par%
	\fi%
	\if\OzsvathSzaboContactCited 1%
		{\bf Ozsv\'ath,\ P.\ \&\ Szab\'o,\ Z.}
		(2005).
		``Heegaard Floer homology and contact structures.''
		{\it Duke Mathematical Journal},
		129(1).
		\par%
	\fi%
	\if\PalisOnMorseSmaleCited 1%
		{\bf Palis,\ J.}
		(1969).
		``On morse-smale dynamical systems.''
		{\it Topology},
		8(4),
		pp. 385--404.
		Elsevier.
		\par%
	\fi%
	\if\PetitSpinCCited 1%
		{\bf Petit,\ R.}
		(2005).
		``Spinc-structures and Dirac operators on contact manifolds.''
		{\it Differential Geometry and its Applications},
		22(2),
		pp. 229--252.
		Elsevier.
		\par%
	\fi%
	\if\RamosAbsoluteCited 1%
		{\bf Ramos,\ V.G.B.}
		(2018).
		``Absolute gradings on ECH and Heegaard Floer homology.''
		{\it Quantum Topology},
		9(2),
		pp. 207--228.
		\par%
	\fi%
	\if\RolfsenKnotsCited 1%
		{\bf Rolfsen,\ D.}
		(2003).
		{\it Knots and links},
		346.
		American Mathematical Soc..
		\par%
	\fi%
	\if\RosoThesisCited 1%
		{\bf Roso,\ B.R.S.}
		(2023).
		{\it Seiberg-Witten Floer Spectra and Contact Structures}.
		Ph.D. Dissertation, University College London.
		\par%
	\fi%
	\if\ShubinPseudodifferentialCited 1%
		{\bf Shubin,\ M.A.}
		(2000).
		{\it Pseudodifferential Operators and Spectral Theory},
		2.
		Spring\-er.
		\par%
	\fi%
	\if\TaubesSeibergWittenSymplecticCited 1%
		{\bf Taubes,\ C.H.}
		(1994).
		``The Seiberg-Witten invariants and symplectic forms.''
		{\it Mathematical Research Letters},
		1(6),
		pp. 809--822.
		International Press of Boston.
		\par%
	\fi%
	\if\TaubesSWToGrCited 1%
		{\bf Taubes,\ C.H.}
		(1996).
		``SW $\Rightarrow$ Gr: from the Seiberg-Witten equations to pseudo-holomorphic curves.''
		{\it Journal of the American Mathematical Society},
		pp. 845--918.
		JSTOR.
		\par%
	\fi%
	\if\TaubesWeinsteinICited 1%
		{\bf Taubes,\ C.H.}
		(2007).
		``The Seiberg--Witten equations and the Weinstein conjecture.''
		{\it Geometry \& Topology},
		11(4),
		pp. 2117--2202.
		Mathematical Sciences Publishers.
		\par%
	\fi%
	\if\TaubesWeinsteinIICited 1%
		{\bf Taubes,\ C.H.}
		(2009).
		``The Seiberg--Witten equations and the Weinstein conjecture II: More closed integral curves of the Reeb vector field.''
		{\it Geometry \& Topology},
		13(3),
		pp. 1337--1417.
		Mathematical Sciences Publishers.
		\par%
	\fi%
	\if\TaubesEmbeddedICited 1%
		{\bf Taubes,\ C.H.}
		(2010a).
		``Embedded contact homology and Seiberg--Witten Floer cohomology I.''
		{\it Geometry \& Topology},
		14(5),
		pp. 2497--2581.
		Mathematical Sciences Publishers.
		\par%
	\fi%
	\if\TaubesEmbeddedIICited 1%
		{\bf Taubes,\ C.H.}
		(2010b).
		``Embedded contact homology and Seiberg--Witten Floer cohomology II.''
		{\it Geometry \& Topology},
		14(5),
		pp. 2583--2720.
		Mathematical Sciences Publishers.
		\par%
	\fi%
	\if\TaubesEmbeddedIIICited 1%
		{\bf Taubes,\ C.H.}
		(2010c).
		``Embedded contact homology and Seiberg--Witten Floer cohomology III.''
		{\it Geometry \& Topology},
		14(5),
		pp. 2721--2817.
		Mathematical Sciences Publishers.
		\par%
	\fi%
	\if\TaubesEmbeddedIVCited 1%
		{\bf Taubes,\ C.H.}
		(2010d).
		``Embedded contact homology and Seiberg-Witten Floer cohomology IV.''
		{\it Geometry \& Topology},
		14(5),
		pp. 2819--2960.
		Mathematical Sciences Publishers.
		\par%
	\fi%
	\if\TaubesEmbeddedVCited 1%
		{\bf Taubes,\ C.H.}
		(2010e).
		``Embedded contact homology and Seiberg-Witten Floer cohomology V.''
		{\it Geometry \& Topology},
		14(5),
		pp. 2961--3000.
		Mathematical Sciences Publishers.
		\par%
	\fi%
	\if\TomDieckTransformationCited 1%
		{\bf tom Dieck,\ T.}
		(1987).
		{\it Transformation Groups}.
		W. de Gruyter.
		\par%
	\fi%
	\if\WuLegendrianCited 1%
		{\bf Wu,\ H.}
		(2006).
		``Legendrian vertical circles in small Seifert spaces.''
		{\it Communications in Contemporary Mathematics},
		8(2),
		pp. 219--246.
		\par%
	\fi%
}}%

\Title{Seiberg-Witten Floer Spectra and Contact Structures}
\Abs
In this \Document,
the author defines an invariant
of rational homology $3$-spheres
equipped with a contact structure
as an element of a cohomotopy set
of the Seiberg-Witten Floer spectrum
as defined in \ManolescuSeibergWittenFloer.
Furthermore,
in light of the equivalence
established in {\LidmanManolescuEquivalence}
between the Borel equivariant homology of said spectrum
and the Seiberg-Witten Floer homology
of \KronheimerMrowkaMonopolesAndThreeManifolds,
the author shall show that
this homotopy theoretic invariant
recovers the already well known contact element
in the Seiberg-Witten Floer cohomology
(vid.\ e.g.\ \KronheimerMrowkaOzsvathSzaboMonopolesAlt)
in a natural fashion.
Next,
the behaviour of the cohomotopy invariant is considered
in the presence of a finite covering.
This setting naturally asks for the use of Borel cohomology
equivariant with respect to the group of deck transformations.
Hence,
a new
equivariant contact invariant is defined
and its properties studied.
The invariant is then computed in one concrete example,
wherein the author demonstrates
that it opens the possibility of considering
scenarios hitherto inaccessible.
\Ack
The author is greatly indebted to his advisor S.~Sivek
for the invaluable guidance
throughout his doctoral studies
and the sage commentary
during the writing of the present \Document.
The author would like to thank
F.~Lin,
L.~Nicolaescu
and
C.H.~Taubes
for insightful email correspondences.
He also thanks the referee
for his suggestions and corrections.
This work was supported by the
Engineering and Physical Sciences Research Council.
\vfootnote{1}{
	Keywords:
	Seiberg-Witten Floer spectra,
	contact structures,
	regular coverings,
	attractor-repeller pairs.
}
\vfootnote{2}
{
	MSC2020:
	57R58,
	57K33,
	57M10,
	37B30
}
\par
\Sect{Introduction}
It is well known from the work of
{\TaubesSeibergWittenSymplectic}
that
symplectic $4$-manifolds
with $b_2^+>1$
have non-trivial Seiberg-Witten invariants.
This is accomplished
by perturbing the Seiberg-Witten equations
in a manner dictated by the symplectic form
multiplied by a large positive real number.
The effect of doing so
is that the Seiberg-Witten equations
gain an obvious canonical solution,
which turns out to be unique and non-degenerate.
\par
Perhaps a little less well known
is that one can do something similar
for contact $3$-manifolds.
This was pursued in
{\TaubesWeinsteinI}
and was ultimately crucial for the proof of the
Weinstein conjecture.
In the $3$-dimensional case,
one adds,
as a perturbation,
the contact form multiplied,
again,
by a positive real number.
A canonical solution to the Seiberg-Witten equations
immediately becomes apparent,
and,
if the real number be made large enough,
one finds that this solution is automatically non-degenerate.
Uniqueness does not hold in the $3$-dimensional case
in the same form as it does in the $4$-dimensional case;
if it did,
all contact rational homology $3$-spheres
would be {\it L-spaces.}
Nonetheless,
another sort of uniqueness does hold
but one which concerns
Seiberg-Witten {\it trajectories.}
The distinguished contact monopole
does not admit any non trivial Seiberg-Witten trajectories
coming into it in the {\it forward time limit.}
This implies that this solution
defines a cocycle
and therefore a class
in {\it monopole Floer cohomology.}
As it turns out,
this class is the well known contact invariant
studied in
{\KronheimerMrowkaOzsvathSzaboMonopoles}
and
is equivalent to the contact invariants
in Heegaard Floer and embedded contact homologies.
\par
In \ManolescuSeibergWittenFloer,
a Seiberg-Witten Floer {\it spectrum\/}
was defined,
which was later shown to recover the monopole Floer cohomologies
through its Borel $\U(1)$-equivariant cohomology.
An important detail here is that
the construction of the spectrum avoids altogether
the use of any {\it generic\/} perturbations.
Taubes' approach to defining the contact invariant,
despite requiring a generic perturbation
in order to work in a Morse theoretic setting,
has the property that
the contact monopole
is already non-degenerate
{\it before\/}
the addition of the generic perturbation.
In the present \Document,
the author applies Taubes' approach to the contact invariant
in the context of the Seiberg-Witten Floer spectrum
in order to conveniently avoid
the use of generic perturbations altogether.
\Thm{
	Given a contact rational homology $3$-sphere
	$(Y,\lambda)$
	there exists a {\it cohomotopical contact invariant,}
	$$
		\SWF(Y,\SpinCStruct_\lambda)
		\to\ContactThom(\lambda),
	$$
	where $\SpinCStruct_\lambda$ is the $\SpinC$ structure
	on $Y$ defined by the contact form $\lambda$
	and $\ContactThom(\lambda)$
	is a (de)suspension of a
	$\U(1)$-equivariant Thom space of a vector bundle
	over the $\U(1)$-orbit of the contact monopole
	in global Coulomb gauge.
	Moreover,
	the classical {\it cohomological\/} contact invariant
	is recovered by pulling back via this map
	a class in the Borel $\U(1)$-equivariant cohomology
	of $\ContactThom(\lambda)$.
}
\par
A similar invariant is constructed in
\IidaTaniguchiSeibergWitten;
however,
the work in the present {\Document}
is entirely independent
and differs in significant ways.
Firstly,
the author uses a different set of analytical results
to ensure the existence of his invariant
and also relies more heavily on certain aspects
of {\it Conley theory\/} to define it.
Secondly,
here,
$\U(1)$-equivariance is kept manifest throughout,
which means one can consider
implications in Borel equivariant cohomology;
indeed,
the author was able to prove,
via the techniques developed in
\LidmanManolescuEquivalence,
how the cohomotopical invariant
recovers the well known cohomological invariant
by passing to Borel equivariant cohomology.
It should also be noted that
the construction presented in the present {\Document}
and the one of {\IidaTaniguchiSeibergWitten}
are sufficiently different
that it is not clear if the two invariants
are indeed equivalent or not.
Of course,
one would be inclined to think that
two such invariants
should really be holding the same information.
However,
proving their equivalence might be a difficult task
due to the different analytical foundations used,
so this goal is not pursued
in the present \Document.
\par
The author's main goal,
after arming himself with the cohomotopical invariant,
was to study {\it covering spaces.}
The avoidance of generic perturbations
is important in this context,
as demonstrated in
\LidmanManolescuCovering,
due to the impossibility of
producing sufficiently generic
{\it equivariant\/} perturbations.
The author's cohomotopical contact invariant
can also be made $G$-equivariant
for $G$ the group of {\it deck transformations\/}
of a finite regular covering.
This allows him to consider
Borel $G$-equivariant cohomology
and deduce certain vanishing and non-vanishing results
via the use of the {\it localization theorem.}
\par
The central question one considers
is whether the lift of a tight contact structure
remains tight
or becomes overtwisted in the covering.
The results derived via use of the contact invariant
shed light on this problem.
In the work of \LinLipnowskiHyperbolic,
the term {\it minimal L-space\/}
is introduced
to refer to rational homology $3$-spheres
having a single solution to the Seiberg-Witten equations
for any $\SpinC$ structure.
For such a manifold,
the Seiberg-Witten Floer spectrum is always the {\it sphere spectrum.}
Moreover,
in the case of a finite covering,
the Seiberg-Witten Floer $G$-spectrum
is the {\it sphere $G$-spectrum.}
In this context,
the author shall establish the following theorem.
\Thm{
	Let
	$\pi:Y\to Y/G$
	be a
	{\it regular prime order covering\/}
	of minimal L-spaces
	and suppose that $\lambda$ be a tight contact form on $Y/G$
	with non-vanishing {\it cohomological\/} contact invariant.
	Then,
	the lifted contact form,
	$\pi^*\lambda$,
	has non-vanishing cohomological contact invariant
	(and,
	therefore,
	is tight)
	provided that
	$d_3(\Ker\pi^*\lambda)+1/2=d(Y,\pi^*\SpinCStruct_\lambda)$,
	where $d_3$ denotes
	Gompf's three-dimensional invariant of hyperplane fields
	and $d$ denotes the Ozsv\'ath-Szab\'o invariant.
}
Examples of minimal L-spaces
include all {\it sol\/} rational homology $3$-spheres
due to the work of {\LinSolv}.
Another example is the {\it Hantzsche-Wendt\/} manifold,
the unique flat rational homology $3$-sphere.
A few more examples exist amongst
the {\it hyperbolic\/} manifolds
as shown in {\LinLipnowskiHyperbolic}.
However,
the most evident class of examples of minimal L-spaces
is that of the {\it elliptic manifolds.}
In this case,
increased knowledge of the contact topology,
in particular the fact that the cohomological contact invariant
never vanishes for tight contact structures
and no two distinct contact structures
have the same $\SpinC$ structure,
allows the author to prove the following stronger theorem.
\Thm{
	Let
	$\pi:Y\to Y/G$
	be a
	{\it prime order regular covering\/}
	of elliptic manifolds
	and suppose that $\lambda$ be a tight contact form on $Y/G$.
	Then,
	the lifted contact form,
	$\pi^*\lambda$,
	is {\it isotopic\/} to a tight contact form
	$\lambda'$
	on $Y$
	if and only if it be {\it homotopic\/}
	to $\lambda'$.
}
This result leads to a scheme
for determining tightness of the lift
of a tight contact structure on an elliptic manifold
based purely on the {\it homotopy theoretic obstruction classes\/}
of the contact structures involved.
This reduces significantly the complexity of the problem
and can be used in concrete calculations.
The rationale is to try to determine the obstruction theoretic invariants
of hyperplane fields --
that is,
the $d_3$ invariant
and the $\SpinC$ structure --
for the {\it lifted contact structure,}
$\pi^*\lambda$,
from those of $\lambda$.
If those be seen to match the values
for a known tight contact structure
on $Y$,
then one shall know that
$\pi^*\lambda$
is isotopic to it.
The issue that arises is that
it is non trivial to determine the lifting behaviour
of the obstruction theoretic invariants,
especially of the $d_3$ invariant.
\par
In order to solve this problem,
the author found himself having to develop techniques
which seem not be discussed in the literature
in a particularly well detailed manner.
These shall be detailed in the present {\Document}
and rely mostly on use of the Kirby calculus.
As an example,
the author shall study the case
a certain tight contact structure
on the $(-8)$-surgery on the left-handed trefoil
which shall be shown to
lift to a virtually overtwisted contact structure
on the lens space $L(12,7)$
via the double covering.
\par
The main ingredient needed to perform these calculations
is a form of {\it $G$-equivariant almost-complex filling\/}
for the given covering of contact manifolds,
which consists of an almost complex
$4$-manifold-with-boundary
extending the given $G$-action on its contact boundary
and potentially having a branching surface in its interior.
Such a filling
can often be produced
by appealing to the notion of
{\it equivariant handle attachments\/}
in the context of the Kirby calculus.
With such a filling at hand,
one can apply the $G$-signature theorem,
as was done by \KhuzamLifting,
to deduce the lifting behaviour of $d_3$ invariants.
\par
The other matter that one must understand carefully
is the lifting behaviour of $\SpinC$ structures.
This is more elementary,
albeit still difficult in practice,
and can be tackled in different manners.
The method pursued here
shall follow a similar approach
to the lifting of $d_3$ invariants
by using Kirby calculus to express {\it Spin\/} structures
in terms of obstruction theory
and then studying the lifting behaviour of Spin structures.
The behaviour of $\SpinC$ structures
follows easily thence.
\par
This work was extracted from parts of the author's
Ph.D. dissertation,
\RosoThesis.

\Sect{Seiberg-Witten Equations and Contact Structures}
This section shall introduce the basic
definitions and analytical results
required from Seiberg-Witten theory.
Throughout this \Document,
the author shall use a version
of the Seiberg-Witten equations
adapted to the presence of a contact form
which was first introduced in
{\TaubesWeinsteinI} and {\TaubesWeinsteinII}
and was subsequently used
in Taubes' work in the correspondence between
monopole Floer homology and embedded contact homology.
\par
Consider an oriented $3$-manifold $Y$
satisfying $b_1(Y)=0$.
A contact form $\lambda$
is a $1$-form on $Y$
satisfying $\lambda\wedge\d\lambda>0$.
The subbundle of $\Tan Y$
given by $\Ker\lambda$
is called a {\it coorientable\/} contact structure.
In this \Document,
all contact structures shall be assumed coorientable.
As shall be seen,
the version of the Seiberg-Witten equations
which shall be used
always admits a canonical solution,
$C_\lambda$,
which,
provided a certain parameter $r>0$
be made large enough,
is nicely behaved in two fundamental ways.
It is non-degenerate irrespective of any genericity requirements,
and it is not the forward time limit of any Seiberg-Witten trajectory.
These two properties shall be instrumental later
in the present \Document.
The solution $C_\lambda$
is essentially defined in a manner
that make its spinor component bounded away from zero;
a feature which is unique to this solution
provided $r$ be made large enough.
\par
In what follows,
agree to fix a metric $g$ on $Y$
with the property that $\lambda\wedge\d\lambda=\Vol_g$.
Use $\xi:=\Ker\lambda$.
Fix a complex structure $J\in\Endo(\xi)$
on the bundle $\xi$
compatible with $g$
in the sense that
$g(-,-)=\d\lambda|_\xi(-,J-)$.
Use $R\in\Gamma\Tan Y$
to denote the Reeb vector field;
that is,
the vector field satisfying
$\iota_R\d\lambda=0$,
$\iota_R\lambda=1$.
Write $\xi\otimes\C=\Exterior^{1,0}\xi\oplus\Exterior^{0,1}\xi$,
where $\Exterior^{1,0}\xi$ and $\Exterior^{0,1}\xi$ are,
respectively,
the $(\pm i)$-eigenbundles of $J$.
Likewise,
for the dual,
write $\xi^*\otimes\C=\Exterior^{1,0}\xi^*\oplus\Exterior^{0,1}\xi^*$.
Denote
$
	\Exterior^{p,q}\xi^*:=
	\Exterior^p_\C\Exterior^{1,0}\xi^*
	\otimes
	\Exterior^q_\C\Exterior^{0,1}\xi^*.
$
There is canonical $\SpinC$ structure on $Y$
defined via the specification
of its spinor representation bundle
in the following way.
\Def{\Label\SpinorsAndCliffMult
	Define the spinor bundle
	$$
		\Spinor_\lambda:=\bigoplus_q\Exterior^{0,q}\xi^*
		=\Exterior^{0,0}\xi^*\oplus\Exterior^{0,1}\xi^*
	$$
	with Clifford multiplication
	$\cliff:\Tan Y\to\Endo_\C(\Spinor_\lambda)$
	defined so as to satisfy,
	for $\alpha\in\Exterior^{0,q}\xi^*$
	and $X\in\xi$,
	$$
		\cliff(R)\alpha=(-1)^{q+1}i\alpha,\quad
		\cliff(X)\alpha={\sqrt 2}\left(
			(X^{0,1})^*\wedge\alpha
			-\iota_{(X^{0,1})}\alpha
		\right).
	$$
}
\Rem{
	The notation
	$(X\mapsto X^*):\Tan Y\otimes\C\to\Cotan Y\otimes\C$
	denotes
	the $\C$-{\it anti\/}linear isomorphism induced by the metric.
}
\Rem{
	Note that this fully determines the Clifford action
	due to the fact that
	$\Tan Y=\langle R\rangle_\R\oplus\xi$.
	This map can be checked to indeed define
	an irreducible Clifford module;
	vid.\ {\PetitSpinC} for a proof
	and more details about this matter.
}
\Def{
	Denote the underlying $\SpinC$ structure,
	that is,
	the principal $\SpinC(3)$-bundle,
	by $\SpinCStruct_\lambda\to Y$.
}
\Rem{
	Notice that
	the determinant line bundle is simply
	$
		\det\SpinCStruct_\lambda
		=\Exterior^{0,1}\xi^*
		\cong\xi^*
	$
	where $\xi^*$ is equipped with the complex structure
	induced by $J$.
}
\Def{
	Use $\Connex(E)$ to denote
	the affine space of Hermitian connexions
	on a Hermitian vector bundle $E\to Y$.
}
\Def{
	Define the {\it configuration space\/}
	by
	$
		\Conf(Y,\SpinCStruct_\lambda):=
		\Connex(\det\SpinCStruct_\lambda)
		\times\Gamma(\Spinor_\lambda)
	$.
	Use 
	$\Conf(Y,\SpinCStruct_\lambda)_k$
	to denote its completion
	in the topology induced by the Sobolev norm
	$\Sobolev k,2$.
}
\Def{
	The tangent bundle
	$\Tan\Conf(Y,\SpinCStruct_\lambda)
	\to\Conf(Y,\SpinCStruct_\lambda)$
	is the bundle having as fibre over
	$C\in\Conf(Y,\SpinCStruct_\lambda)$
	the space
	$
		\Tan\Conf(Y,\SpinCStruct_\lambda)|_C
		=\Gamma(i\Cotan Y\oplus\Spinor_\lambda)
	$.
	Use 
	$\Tan\Conf(Y,\SpinCStruct_\lambda)_k$
	for its Sobolev $\Sobolev k,2$ completion.
}
\Def{
	Use $\tau:\Spinor_\lambda\to i\Cotan Y$
	to denote the quadratic map
	defined by sending
	$$
		\psi\mapsto
		\cliff^{-1}\left(
			\psi\otimes\psi^*
			-{1\over2}|\psi|^2\id
		\right).
	$$
}
\Def{
	Denote
	$
		\Dirac_A:
		\Gamma\Spinor_\lambda
		\to\Gamma\Spinor_\lambda
	$
	the {\it Dirac operator\/} defined by a connexion
	$A$ on $\det\Spinor_\lambda$.
}
\Def{
	Define the {\it contact configuration,}
	$
		\ContactConf_\lambda
		\equiv
		(A_\lambda,\psi_\lambda)
		\in \Connex(\det\SpinCStruct_\lambda)
		\oplus\Spinor_\lambda
	$
	by setting
	its spinor component
	to be the constant function
	$\psi_\lambda:=1\in\Gamma\xi^*_{0,0}\cong Y\times\C$
	and
	by requiring its connexion component $A_\lambda$
	to solve the Dirac equation
	$\Dirac_{A_\lambda}\psi_\lambda=0$.
}
\Rem{
	This canonical configuration
	shall play a pivotal r\^ole in this \Document.
}
\Def{\Label\SWDef
	Let $r>0$.
	The
	{\it canonically perturbed
	Seiberg-Witten vector field\/}
	of $(Y,\lambda)$
	is the vector field
	$
		\SW_{\lambda,r}:
		\Conf(Y,\SpinCStruct_\lambda)
		\to\Tan\Conf(Y,\SpinCStruct_\lambda)
	$
	given by
	$$
		\SW_{\lambda,r}(A,\psi)=
		\left(
			\Hodge{1\over2}\left(
				F_A
				-F_{A_\lambda}
			\right)
			+r\tau(\psi)
			-{ir\over2}\lambda
			,\;
			\Dirac_A\psi
		\right).
	$$
}
\Rem{
	Notice that,
	for any value of $r>0$,
	the contact configuration
	solves the Seiberg-Witten equation
	$$
		\SW_{\lambda,r}(\ContactConf_\lambda)=0.
	$$
	That is,
	$\ContactConf_\lambda$ is a {\it fixed point\/}
	of the Seiberg-Witten vector field.
}
\Def{
	Denote the {\it gauge group\/}
	by $\Gauge(Y):=\DiffClass\infty(Y,\U(1))$.
	Use $\Gauge(Y)_k$
	for its completion in the Sobolev $\Sobolev k,2$
	norm.
}
\Def{
	For $(A,\psi)\in\Conf(Y,\SpinCStruct_\lambda)$,
	define the {\it linearized gauge action\/}
	by
	$$
		\GaugeLieAct_{(A,\psi)}:
		\DiffClass\infty(Y,i\R)
		\to\Tan\Conf(Y,\SpinCStruct_\lambda)|_{(A,\psi)},
		\quad
		u\mapsto (-\d u, u\psi).
	$$
	Use $\GaugeLieAct_{(A,\psi)}^*$
	for its formal $\LTwo$-adjoint.
}
\Def{
	\Label\DefLocalCoulomb
	The {\it local Coulomb gauge\/}
	is the subspace
	$$
		\LCSlice_C
		:=\Ker\GaugeLieAct_C^*
		\subset
		\Tan\Conf(Y,\SpinCStruct_\lambda)|_C.
	$$
	Denote by $\LCSlice\to\Conf(Y,\SpinCStruct_\lambda)$
	the vector bundle with fibres
	the local Coulomb gauges.
	Use
	$\LCSlice_{C,k}$
	and $\LCSlice_k$
	to denote the $\Sobolev k,2$ Sobolev completions.
}
\Def{
	Denote by
	$
		\LCProj_C:
		\Tan\Conf(Y,\SpinCStruct_\lambda)|_C
		\to\LCSlice_C
	$
	the $\LTwo$-orthogonal projection.
}
\Def{
	A configuration
	$(A,\psi)\in\Conf(Y,\SpinCStruct_\lambda)$
	is said to be {\it irreducible\/}
	if $\psi$ is not identically zero.
}
\Def{
	\Label\NonDegenDef
	An irreducible solution
	$C\in\SW_{\lambda,r}^{-1}(0)$
	is said to be {\it non-degenerate\/}
	if the derivative
	$$
		\LCProj_C\circ\D_C\SW_{\lambda,r}:
		\LCSlice_C\to\LCSlice_C
	$$
	be surjective.
}
\Rem{
	It is customary in Seiberg-Witten theory
	to ensure non-degeneracy of all solutions
	through the addition of a generic perturbation;
	however,
	often,
	doing so has its downside.
	One of the key advantages
	of the Seiberg-Witten Floer spectrum construction
	(\ManolescuSeibergWittenFloerAlt)
	is that it avoids the need for such a
	perturbation altogether.
	This was crucial in the results of
	{\LidmanManolescuCovering}
	concerning the Smith-type inequality
	of Seiberg-Witten Floer homology.
	Avoiding the use of a generic perturbation
	shall also be exploited in the present \Document;
	however,
	as will be seen,
	it is still necessary to make sure that
	$\ContactConf_\lambda$
	be non-degenerate.
	This is ensured by the following theorem of Taubes.
}
\Rem{
	For ease of reference,
	the author shall make use of the following convention.
	When a proposition state the existence
	of a certain constant
	which shall be of use later in the text,
	that constant shall be labelled
	with the number of the proposition.
}
\Thm{
	(\TaubesWeinsteinIAlt)
	\Label\ContactConfNonDegen
	\LabelNum\ContactConfNonDegenNum
	There exists $r_\ContactConfNonDegenNum>0$
	such that,
	for $r>r_\ContactConfNonDegenNum$,
	the contact configuration
	is non-degenerate.
}
\Proof{
	Vid.\ \TaubesWeinsteinI,
	Lemma 3.3.	
}
Besides the non-degeneracy property,
the configuration $C_\lambda$
enjoys two {\it uniqueness\/} properties
that shall prove important.
\Thm{
	(\TaubesWeinsteinIIAlt)
	\Label\SpinorContactConfBoundedbelow
	\LabelNum\SpinorContactConfBoundedbelowNum
	There exists
	$r_\SpinorContactConfBoundedbelowNum>0$
	and $\delta_\SpinorContactConfBoundedbelowNum>0$
	such that,
	for $r>r_\SpinorContactConfBoundedbelowNum$,
	the only configuration
	$C=(a, \psi)$,
	up to a gauge transformation,
	satisfying
	$$
		\SW_{\lambda,r}(C)=0,
		\quad
		|\psi|\geq 1-\delta_\SpinorContactConfBoundedbelowNum
	$$
	is the contact configuration
	$C=\ContactConf_\lambda$.
}
\Proof{
	Vid.\ \TaubesWeinsteinII,
	Proposition 2.8.
}
\Thm{
	\Label\NoTrajToContactConf
	\LabelNum\NoTrajToContactConfNum
	There exists $r_\NoTrajToContactConfNum>0$
	such that,
	for $r>r_\NoTrajToContactConfNum$,
	any trajectory	
	$\gamma:\R\to\Conf(Y,\SpinCStruct_\lambda)$
	satisfying
	$$
		\dbyd t\gamma(t)=
		-\SW_{\lambda,r}\big(\gamma(t)\big),\quad
		\lim_{t\to-\infty}\gamma(t)=C,\quad
		\lim_{t\to\infty}\gamma(t)=C_\lambda,
	$$
	where the limits are with respect to
	the Sobolev norm $\Sobolev k,2$
	for any $k\geq 5$,
	must satisfy $\gamma(t)=C_\lambda$
	for all $t\in\R$.
}
\Proof{
	This
	is {\it nearly\/} what is stated by
	\TaubesWeinsteinII,
	Proposition 5.15,
	but not quite.
	Please find an adaptation of Taubes' proof
	in the appendix of the present \Document.
}

\Sect{Review of Seiberg-Witten Floer Spectra}
Armed with the analytic results
from the previous section,
the first goal of the present {\Document} shall be
to define a homotopy theoretic invariant
emerging from the contact configuration.
This invariant shall live
in an equivariant cohomotopy set
of the Seiberg-Witten Floer Spectrum
$\SWF(Y,\SpinCStruct_\lambda)$
defined by \ManolescuSeibergWittenFloer.
This section shall review that construction
with a slight adaptation;
the Seiberg-Witten flow used
shall be the one canonically perturbed
by the contact form
as was described in the previous section.
This shall allow for the definition
of the contact invariant
in the next section.
\par
Let $Y$, $\lambda$ and $g$
be as in the previous section.
\Def{
	The {\it unperturbed Seiberg-Witten vector field\/} is
	$$
		\SW:
		\Conf(Y,\SpinCStruct_\lambda)
		\to\Tan\Conf(Y,\SpinCStruct_\lambda),
		\quad
		\SW(A,\psi)=
		\left(
			{1\over2}\Hodge F_A
			-\tau(\psi),\;
			\Dirac_A\psi
		\right).
	$$
}
\Rem{
	The construction in
	{\ManolescuSeibergWittenFloer}
	uses $\SW$
	as the Seiberg-Witten vector field.
	The version of the Seiberg-Witten vector field
	used in the present {\Document}
	is $\SW_{\lambda,r}$
	and it differs from $\SW$ in two ways.
	Firstly,
	the spinor component of
	$\SW_{\lambda,r}$
	is scaled by $r^{1/2}$
	compared to $\SW$;
	this distinction shall be evidently immaterial
	in the construction.
	Secondly,
	$\SW_{\lambda,r}$
	contains the constant term
	$((ir/2)\lambda-\Hodge(1/2)F_{A_\lambda},0)$
	added on.
	In the language of
	\LidmanManolescuEquivalence,
	this amounts to the addition of a
	``very tame'' perturbation
	(vid.\ \LidmanManolescuEquivalenceAlt,
	Definition 4.4.2),
	which,
	as demonstrated there
	(vid.\ \LidmanManolescuEquivalenceAlt,
	Proposition 6.1.6),
	does not affect the construction of the spectrum.
	Hence,
	the spectra defined with $\SW_{\lambda,r}$
	and $\SW$
	shall be the same.
}
\Def{
	The {\it normalized Gauge group\/}
	is the subgroup $\NormGauge(Y)\subset\Gauge(Y)$
	consisting of those $u\in\Gauge(Y)$
	which can be written as
	$u=e^{if}$
	such that $\int_Y \Hodge f = 0$.
}
\Def{
	By the {\it global Coulomb gauge
	with respect to the connexion $A_\lambda$,}
	one means
	$$
		W:=\left(
			A_\lambda+
			\Ker(
				\d^*:i\Omega^1(Y)
				\to i\Omega^0(Y)
			)
		\right)
		\oplus\Gamma(\Spinor_\lambda)
		\subset\Conf(Y,\SpinCStruct_\lambda).
	$$
}
\Rem{
	The affine space
	$W$ shall be thought of as a vector space
	with zero being $(A_\lambda,0)$;
	that is,
	connexions shall be thought
	of as purely imaginary $1$-forms
	by subtracting $A_\lambda$.
}
\Rem{
	Any $\NormGauge(Y)$-equivalance class
	$[(A_\lambda+a,\psi)]\in\Conf(Y,\SpinCStruct_\lambda)/\NormGauge(Y)$
	has a unique representative
	in global Coulomb gauge;
	that is,
	there is a unique $(A_\lambda+a',\psi')\in W$
	such that $[(A_\lambda+a',\psi')]=[(A_\lambda+a,\psi)]
	\in\Conf(Y,\SpinCStruct_\lambda)/\NormGauge(Y)$.
}
\Def{
	The {\it global Coulomb projection},
	$\GCProj:\Conf(Y,\SpinCStruct_\lambda)\to W$,
	is given
	by sending a configuration $(A_\lambda+a,\psi)\in\Conf(Y,\SpinCStruct_\lambda)$
	to its unique $\NormGauge(Y)$-equivalent
	in global Coulomb gauge.
}
\Rem{
	The map $\GCProj$ may be computed as follows.
	Let $G:\Sobolev m,2(Y)\to\Sobolev {m+2},2(Y)$
	denote the Green's operator
	of the Laplacian $\Delta:\Omega^0(Y)\to\Omega^0(Y)$.
	One can show that
	$$
		\GCProj(A_\lambda + a,\psi)
		=(A_\lambda + a - \d G\d^* a, e^{G\d^* a}\psi).
	$$
	As a consequence,
	note that $\GCProj$
	maps bounded sets to bounded sets.
}
\Def{
	The {\it enlarged local Coulomb slice\/}
	is the subspace
	$$
		\ELCSlice_{(A,\psi)}
		\subset
		\Tan\Conf(Y,\SpinCStruct_\lambda)|_{(A,\psi)}
	$$
	defined as the $\LTwo$-orthogonal complement
	to the orbits of $\NormGauge(Y)$.
}
\Def{
	Denote by $\NormGaugeLie(Y)$
	the Lie algebra of $\NormGauge(Y)$.
}
\Rem{
	Any equivalence class
	$[(b,\phi)]\in\Tan\Conf(Y,\SpinCStruct_\lambda)|_{(A,\psi)}/\NormGaugeLie(Y)$
	has a unique representative
	in enlarged local Coulomb gauge.
}
\Def{
	\Label\EnlargedLocalCoulombProjection
	By the {\it enlarged local Coulomb projection},
	one means
	$$
		\ELCProj_{(A,\psi)}:
		\Tan\Conf(Y,\SpinCStruct_\lambda)|_{(A,\psi)}
		\to\ELCSlice_{(A,\psi)}
	$$
	defined by sending
	a vector
	to the unique representative
	in enlarged local Coulomb gauge
	of its equivalence class
	in the quotient
	$\Tan\Conf(Y,\SpinCStruct_\lambda)|_{(A,\psi)}/\NormGaugeLie(Y)$.
}
\Rem{
	Note that $\LCProj$ and $\ELCProj$
	are maps defined on the tangent bundle $\Tan\Conf(Y,\SpinCStruct_\lambda)$,
	whereas $\GCProj$ is defined on $\Conf(Y,\SpinCStruct_\lambda)$.
	Of course,
	$\GCProj$ induces a map $\Tan\Conf(Y,\SpinCStruct_\lambda)\to\Tan W$
	via the pushforward
	$\GCProj_*$.
}
\Def{
	Set $\GCSW_{\lambda,r}:=\GCProj_*\SW_{\lambda,r}$.
}
\Rem{
	Fix some integer $k\geq5$
	and
	consider,
	henceforth,
	$\GCSW_{\lambda,r}$
	as a map $W_k\to W_{k-1}$
	where $W_m$ denotes the completion of $W$
	in the Sobolev norm $\Sobolev m,2$.
}
\Def{
	Define the Fredholm linear operator
	$\ell:W_k\to W_{k-1}$
	by the formula
	$$
		\ell(A_\lambda + a,\psi)
		=\left(
			A_\lambda + {1\over2}\Hodge\d a,
			\;
			\Dirac_{A_\lambda}\psi
		\right).
	$$
}
\Def{
	Define the (non-linear) operator
	$c:W_k\to W_{k-1}$
	by $c:={\GCSW_{\lambda,r}}-\ell$.
}
\Rem{
	Note
	that $\GCSW_{\lambda,r}=\ell+c$
	where
	$c$ is a {\it compact\/} operator
	as explained in
	\ManolescuSeibergWittenFloer, \S4.
}
\Def{
	\Label\SpectralCutOffsNonSmooth
	For $\mu>1$,
	denote by $W^\mu\subset W_k$,
	the subspace consisting of
	the span of the eigenvectors of $\ell$
	with eigenvalues in the interval $(-\mu,\mu)$.
	Use $\tilde p^\mu:W_k\to W^\mu$
	to denote the $\LTwo$-orthogonal projection.
}
The family of operators $\tilde p^\mu$
must now be smoothed out
in a particular way.
For that end,
fix a smooth function
$\beta:\R\to\R$
satisfying $\supp\beta=[0,1]$
and $\int_\R\beta(x)\d x = 1$.
A preliminary version
of the smoothed out family is as follows.
\Def{
	Define a family of operators
	$p^\mu_{\rm prel}:W_k\to W^\mu$
	by
	$$
		p^\mu_{\rm prel}:=
		\int_0^1 \beta(t)\tilde p^{\mu-t}\d t
	$$
}
\Rem{
	\Label\LidmanManolescuTechnicalConditionOnSpectralCutOff
	This preliminary version
	could well be used
	to define the Seiberg-Witten Floer spectrum
	and,
	indeed,
	is essentially the operator family
	which appears in the original definition in
	\ManolescuSeibergWittenFloer.
	However,
	in \LidmanManolescuEquivalence,
	the authors
	use a slightly modified version
	which turns out to be needed
	in proving some technical results.
	Some of those technical results
	shall be used in the present \Document.
	To define the final version of the operator family,
	a few more data need to be fixed.
	Firstly,
	choose an unbounded strictly increasing sequence
	$\{\mu_i\}\subset\R$
	such that,
	for no $i$,
	be $\mu_i$
	an eigenvalue of $\ell$.
	Next,
	fix a sequence of small real numbers
	$\{\epsilon_i\}\subset\R$
	such that
	the intervals
	$[\mu_i-\epsilon_i,\mu_i+\epsilon_i]$
	be disjoint
	and not contain any eigenvalue of $\ell$.
	At last,
	pick smooth bump functions
	$\{\beta_i:\R\to[0,1]\}$
	such that
	$
		\supp\beta_i
		\subset[\mu_i-\epsilon_i,\mu_i+\epsilon_i]
	$.
}
\Def{
	\Label\SpectralCutOffsSmooth
	Define the family of operators
	$p^\mu:W_k\to W^\mu$
	by
	$$
		p^\mu
		:=\sum_i\beta_i(\mu)\tilde p^\mu
		+\left(
			1
			-\sum_i\beta_i(\mu)
		\right)
		p^\mu_{\rm prel}.
	$$
}
\Rem{
	The family of operators $p^\mu$
	is smooth in $\mu$
	but still has the property that,
	for all $i$,
	$p^{\mu_i}=\tilde p^{\mu_i}$.
}
\Def{
	\Label\CSDDef
	By the {\it canonically perturbed
	Chern-Simons-Dirac functional},
	one means
	$$
		\CSD_{\lambda,r}:\Conf(Y,\SpinCStruct_\lambda)\to\R,
	$$
	$$
		\CSD_{\lambda,r}(A_\lambda+a,\psi)
		:=\int_Y \Hodge r\langle\psi,\Dirac_{A_\lambda+a}\psi\rangle
		-\int_Y a\wedge\d a
		-{ri\over2}\int_Y \lambda\wedge\d a.
	$$
}
\Def{
	A {\it finite type curve\/}
	$\gamma:\R\to\Conf(Y,\SpinCStruct_\lambda)$,
	$\gamma=(A,\psi)$,
	is a curve such that
	the maps 
	$t\mapsto\CSD_{\lambda,r}(\gamma(t))$
	and $t\mapsto\norm{\psi(t)}_{{\rm C}^0}$
	be bounded
	as functions $\R\to\R$.
}
\Def{
	A curve $\gamma:\R\to W_k$
	is said to be a
	{\it Seiberg-Witten trajectory
	in global Coulomb gauge\/}
	if
	$$
		\dbyd t\gamma = -\GCSW_{\lambda,r}(\gamma(t)).
	$$
}
\Def{
	For $R>0$,
	and a normed vector space $V$,
	use $\Ball(V,R)\subset\Disk(V,R)\subset V$
	to denote,
	respectively,
	the open and closed balls of radius $R$.
	Use $\Sphere(V,R)=\Disk(V,R)\setminus\Ball(V,R)$
	to denote the sphere of radius $R$.
}
\Thm{\Label\FiniteTypeBounded
	(cf.\ \ManolescuSeibergWittenFloerAlt,
	Proposition 1)
	There exists $R>0$
	such that all finite type trajectories
	of $\GCSW_{\lambda,r}$
	are contained in the ball
	$\Ball(W_k,R)\subset W_k$.
}
\Proof{
	Firstly, note that
	the proof in {\ManolescuSeibergWittenFloer}
	can be easily adapted
	to the present case of the
	{\it perturbed\/}
	Seiberg-Witten flow.
	That result provides a constant $R'>0$
	such that,
	{\it up to a gauge transformation},
	all Seiberg-Witten trajectories
	of finite type
	sit inside the ball
	$\Ball(W_k,R')\subset\Conf(Y,\SpinCStruct_\lambda)_k$.
	Therefore,
	a Seiberg-Witten trajectory
	in global Coulomb gauge
	is,
	locally,
	the global Coulomb projection
	of a Seiberg-Witten trajectory
	residing in the ball
	$\Ball(W_k,R')\subset\Conf(Y,\SpinCStruct_\lambda)_k$.
	But the global Coulomb projection map
	$\GCProj:\Conf(Y,\SpinCStruct_\lambda)_k\to W_k$
	maps bounded sets to bounded sets.
}
\Rem{
	Henceforth,
	assume
	$R>0$
	to be
	such that all
	finite type
	Seiberg-Witten trajectories
	in Coulomb gauge
	fit in $B(W_k, R)$.
}
\Rem{\Label\FixBumpFunction
	Also,
	fix a family of
	$\U(1)$-equivariant bump functions
	$u^\mu:W^\mu\to\R$
	satisfying
	$$
		u^\mu|_{\Disk(W^\mu,2R)} = 1,
		\quad
		u^\mu|_{W^\mu\setminus\Ball(W^\mu,3R)} = 0
	$$
	and,
	$u^\mu$
	constant
	on the sphere
	$\Sphere(W^\mu,t))\subset W^\mu$
	of radius $t$
	for all $t\in[0,\infty)$.
	Note that the norm on $W^\mu\subset W_k$
	is defined by the Sobolev norm $\Sobolev k,2$
	of $W_k$.
}
\Def{
	Define the
	{\it finite dimensional approximation
	to the Seiberg-Witten vector field\/}
	as
	$$
		\FDASW^\mu_{\lambda,r}
		=u^\mu
		\cdot
		(\ell+p^\mu c).
	$$
	Furthermore,
	use
	$
		\FDASWFlow^\mu_{\lambda,r}:
		W^\mu\times\R
		\to W^\mu
	$
	to denote the flow
	given by the O.D.E.
	$$
		\derzero t\gamma(t)
		=-\FDASW^\mu_{\lambda,r}(\gamma(t)).
	$$
	The flow lines of $\FDASWFlow^\mu_{\lambda,r}$
	are called
	{\it approximate Seiberg-Witten trajectories
	in global Coulomb gauge.}
}
\Def{
	\Label\IsolInvSetForSWFDef
	Use
	$S^\mu_{\lambda,r}\subset\Ball(W^\mu,R)$
	to denote the union
	of all flow lines of
	$\FDASWFlow^\mu_{\lambda,r}$
	which remain inside of $B(W^\mu,R)$
	for all time.
}
\Thm{\Label\FDCompactness
	(\ManolescuSeibergWittenFloerAlt,
	Proposition 3)
	For $\mu>0$
	sufficiently large
	compared to $R$,
	any flow line
	of $\FDASWFlow^\mu_{\lambda,r}$
	which be contained in the disk
	$\Disk(W^\mu,2R)$
	is,
	in fact,
	contained in the open ball
	$\Ball(W^\mu,R)$.
}
\Proof{
	The proof in {\ManolescuSeibergWittenFloer}
	of the non-perturbed case of this theorem
	can be trivially adapted to the case at hand.
	Alternatively,
	this result is a special case of
	\LidmanManolescuEquivalence,
	Proposition 6.1.2(i)
	and
	Proposition 6.1.5,
	noting that the perturbation used in this {\Document}
	is ``very tame''.
}
The author shall now recall
the pertinent definitions
from Conley theory.
For a reference,
the reader is directed to
\ConleyIsolated,
{\FloerRefinement}
and \MischaikowConley.
In what follows,
suppose that $G$ be a compact Lie group,
$\Gamma$ be a locally compact Hausdorff space
with a continuous $G$-action
and $\phi:\Gamma\times\R\to\Gamma$
a continuous and equivariant flow.
\Def{
	\Label\MaxInvSetDef
	Let $U\subset\Gamma$
	be a $G$-invariant subset,
	the {\it maximal invariant set\/} of $U$
	is
	$$
		\Inv(U):=
		\{u\in U\mid(\forall t\in\R)(\phi(u,t)\in U)\}.
	$$
}
\Def{
	\Label\IsolInvSetDef
	Let $S\subset\Gamma$
	be a $G$-invariant compact subset.
	$S$ is called an {\it isolated invariant set\/}
	if there be a compact neighbourhood $U\supset S$
	such that $\Inv(U)=S$.
}
\Def{
	\Label\IndexPairDef
	Let $S\subset\Gamma$ be an isolated invariant set.
	A pair of $G$-invariant compact sets $(M,N)$
	satisfying $N\subset M\subset\Gamma$
	is called an {\it index pair\/} for $S$ when:
	{
		\parindent=1cm
		\item{(i)}
			$M\setminus N$ be an isolating neighbourhood
			for $S$;
		\item{(ii)}
			for all $t\geq 0$ and $x\in N$,
			if $\varphi(\{x\}\times[0,t])\subset M$,
			then $\varphi(\{x\}\times[0,t])\subset N$;
		\item{(iii)}
			for all $t\geq 0$ and $x\in M$,
			if $\varphi(x,t)\notin M$,
			then
			$\varphi(\{x\}\times[0,t])\cap N\neq\emptyset$.
	}
}
\Thm{
	\Label\IndexPairExists
	(\ConleyIsolatedAlt, non-equivariant;
	{\FloerRefinementAlt}
	and
	{\FloerZehnderEquivariantConleyAlt}, equivariant)
	For any isolated invariant set $S\subset\Gamma$,
	there exists an index pair $(M,N)$
	and the $G$-equivariant pointed homotopy type $M/N$
	is independent of the choice of $(M,N)$.
}
\Def{
	The $G$-equivariant homotopy type
	of $M/N$
	where $(M,N)$ is an index pair
	for an isolated invariant set $S$
	is called the {\it Conley index\/} of $S$
	and denoted $\Conley_G(S,\phi)$.
}
{\FDCompactness}
can now be reinterpreted in this language.
\Cor{
	$S^\mu_{\lambda,r}$ is a $\U(1)$-invariant isolated invariant set
	with isolating neighbourhood
	$\Disk(W^\mu,2R)$
	for the $\U(1)$-equivariant flow $\FDASWFlow^\mu_{\lambda,r}$.
}
The author shall now introduce
the relevant definitions from
{\it equivariant stable homotopy theory.}
Here,
the author shall deviate slightly
from the route taken in
\ManolescuSeibergWittenFloer;
this is done in the interests
of later sections that shall deal
with Seiberg-Witten Floer spectra
equivariant with respect to the deck transformations
of a covering.
In \ManolescuSeibergWittenFloer,
in order to perform the required {\it desuspensions,}
an {\it ad hod\/} version of the
Spanier-Whitehead category is used.
In the present \Document,
instead,
the author shall use the,
by now,
more standard category of spectra.
This increases slightly the complexity
of the definitions,
but nothing new is gained
as the Spanier-Whitehead category
embeds into the category of spectra
in a simple way.
For more details,
the reader is directed to
\MayEtAlEquivariant.
In what follows,
let $G$ be a compact Lie group.
Whenever the author say {\it $G$-space,}
he means in fact
{\it pointed $G$-space.}
\Def{
	A {\it $G$-universe $\Universe$\/}
	is an orthogonal $G$-representation
	of countable dimension
	having the following two properties:
	{
		\parindent = 1cm
		\item{(i)}
			For each finite dimensional subrepresentation
			$V\subset\Universe$,
			the direct sum of $V$
			with itself countably many times,
			$V^\infty$,
			also occurs as a subrepresentation
			in $\Universe$.
		\item{(ii)}
			The trivial representation $\R$,
			therefore also $\R^\infty$,
			occurs as a subrepresentation in $\Universe$.
		\par
	}
}
\Def{
	For $V$ a $G$-representation,
	denote by $V^+$ its
	{\it one-point compactification;}
	note that $V^+$ is a $G$-space
	and call it a {\it representation sphere of $G$.}
	For $X$ a $G$-space,
	the {\it $V^{\rm th}$ suspension\/}
	of $X$
	is the smash product
	$\Sigma^V X:=V^+\wedge X$.
	The {\it $V^{\rm th}$ loop space\/}
	of $X$
	is the $G$-space
	$\Omega^V X$
	of all maps $V^+\to X$
	with $G$ acting by conjugation.
}
\Rem{
	There is an {\it adjunction\/} between the suspension functor
	$\Sigma^V$ and the loop space functor $\Omega^V$ on $G$-spaces.
	That is to say that,
	for $G$-spaces $X_1$ and $X_2$,
	there is a {\it natural bijection\/}
	between the space of maps $\Sigma^V X_1\to X_2$
	and the space of maps $X_1\to\Omega^V X_2$.
}
\Def{
	A {\it $G$-prespectrum $E$
	indexed on the $G$-universe $\Universe$\/}
	consists of the following data.
	{
		\parindent = 1cm
		\item{(i)}
			A set of $G$-spaces,
			$E_V$,
			one for each finite dimensional
			subrepresentation $V$
			in $\Universe$.
		\item{(ii)}
			A set of
			$G$-equivariant
			{\it structure maps,}
			$\sigma_{V,W}:\Sigma^{W-V}E_V\to E_W$,
			one for each pair of
			nested finite dimensional subrepresentations,
			$V\subset W\subset\Universe$,
			where $W-V$ denotes the orthogonal complement
			of $V$ in $W$.
		\par
	}
}
\Def{
	A {\it map of $G$-prespectra\/}
	$f$
	from a $G$-prespectrum $E$
	to a $G$-prespectrum $F$,
	both indexed on the same universe $\Universe$,
	consists of a set of $G$-equivariant maps of $G$-spaces
	$f_V:E_V\to F_V$,
	one for each finite dimensional subrepresentation
	$V\subset\Universe$,
	such that the evident diagrams
	$$
		\Diagram{
			\Sigma^{W-V}E_V
			&
			\Right
			&
			E_W
			\cr
			\Down
			&
			&
			\Down
			\cr
			\Sigma^{W-V}F_V
			&
			\Right
			&
			F_W
			\cr
		}
			
	$$
	all commute
	for any pair of nested representations
	$V\subset W\subset\Universe$.
	A {\it homotopy\/}
	between two maps of prespectra $f,g:E\to F$
	is a map of prespectra $h:E\wedge I_+\to F$
	where $I_+$ denotes the interval $[0,1]$
	with a disjoint base point added;
	here,
	the smash $E\wedge I_+$
	between a prespectrum and a space
	is simply to be interpreted spacewise.
	A map of prespectra $f$
	is called a {\it weak equivalence\/}
	if all of its constituent maps of spaces,
	$f_V$,
	be weak equivalences.
}
\Def{
	The {\it suspension prespectrum functor,}
	$\Sigma^\infty_\Universe$,
	from $G$-spaces to $G$-prespectra
	indexed on the universe $\Universe$
	is defined by assigning to a $G$-space $X$
	the prespectrum $\Sigma^\infty_\Universe X$
	consisting of
	$(\Sigma^\infty_\Universe X)_V:=\Sigma^V X$
	and structure maps
	the identity maps.
}
\Def{
	Given a $G$-representation $V$
	in the universe $\Universe$
	and $X$ a $G$-space,
	the {\it $V^{\rm th}$ desuspension of $X$,}
	denoted
	$\Sigma^{-V}\Sigma^\infty_\Universe X$,
	is a $G$-prespectrum indexed on $\Universe$
	defined as follows.
	For $W\subset\Universe$,
	the $W^{\rm th}$ space of
	$\Sigma^{-V}\Sigma^\infty_\Universe X$
	is either a single point,
	in the event that $V$ not be contained in $W$,
	or it is the $G$-space
	$\Sigma^{W-V} X$,
	in the event that $V$ be contained in $W$.
	Meanwhile,
	the structure maps
	$\Sigma^{U-W}\Sigma^{W-V}X\to \Sigma^{U-V}X$
	are the evident ones.
}
\Def{
	A $G$-prespectrum $E$ indexed on
	the $G$-universe $\Universe$
	is called a {\it spectrum\/}
	whenever all the adjoints,
	$E_V\to \Omega^{W-V}E_W$,
	to the structure maps
	be {\it homeomorphisms.}
	The {\it category\/} of $G$-spectra
	and maps as defined above
	for prespectra is denoted
	$\Spectra G\Universe$.
	The {\it homotopy category\/} of $G$-spectra
	is defined as the category with the same objects as
	$\Spectra G\Universe$
	but with morphisms being the homotopy classes
	of maps of prespectra;
	this category is denoted $\HoSpectra G\Universe$.
	The {\it stable homotopy category\/} of $G$-spectra
	consists of the category $\HoSpectra G\Universe$
	together with formal inverses
	for all the weak equivalences;
	this category is denoted $\StabHoSpectra G\Universe$.
}
\Thm{
	The forgetful functor from spectra to prespectra
	has a
	{\it right adjoint\/}
	called
	the {\it spectrification functor.}
}
\Proof{
	Vid.\ \MayEtAlEquivariant, \S XII.2.
}
\Def{
	By composing the {\it suspension prespectrum functor,}
	$\Sigma^\infty_\Universe$,
	with the {\it spectrification functor,}
	one obtains a functor from $G$-spaces
	to $\Spectra G\Universe$.
	This functor shall be called the
	{\it suspension spectrum functor\/}
	and shall also be denoted by
	$\Sigma^\infty_\Universe$.
	Likewise,
	the desuspension of a $G$-space $X$,
	$\Sigma^{-V}\Sigma^\infty_\Universe X$
	can be regarded as being in $\StabHoSpectra G\Universe$.
}
\Rem{
	As the notation suggests,
	there is a desuspension functor,
	$\Sigma^{-V}$,
	defined for all spectra,
	not just for suspension spectra.
	However,
	it is somewhat more subtle to define
	and the author shall not require it
	in the present \Document.
}
\Rem{
	The stronger notion of spectra
	as opposed to prespectra
	is not so important in the present {\Document}
	because the main desire is simply
	to be able to perform {\it desuspensions.}
	Nonetheless,
	it has become standard in the literature
	to work with spectra
	because of their ability to classify
	homology and cohomology theories
	and also the superior properties
	that the category of spectra enjoys.
	Therefore,
	the author decided to phrase everything
	in terms of spectra
	for ease of reference.
}
\Rem{
	Notice that the Coulomb gauge $W$
	is a $\U(1)$-universe
	isomorphic to $\R^\infty\oplus\C^\infty$.
	Such an isomorphism can be defined
	by picking a basis of eigenvectors of $\ell$.
	Note that this universe does {\it not\/}
	contain all representations of $\U(1)$.
	Indeed,
	if $\C$ denote the standard representation,
	where $\U(1)\incl\C$
	is the unit circle,
	then the tensor product $\C\otimes_\C\C$
	is not in the universe.
	This choice of universe is compatible
	with what is done in
	\ManolescuSeibergWittenFloer,
	where desuspensions are only allowed
	with respect to $\R$ and $\C$.
}
\Def{
	When thinking of $W$ as a universe,
	denote it by $\CoulombUniverse$.
}
\Def{\Label\WInterval
	Given an interval $I\subset\R$,
	use $W^{I}$
	to denote the span of the eigenvectors
	of $\ell$
	with eigenvalues in the interval $I$.
	Hence,
	e.g.,
	$W^\mu=W^{(-\mu,\mu)}$.
}
\Thm{
	\Label\SWFWellDefWithCoulombUniverse
	(\ManolescuSeibergWittenFloerAlt, \S 7)
	Up to canonical isomorphism in
	$\StabHoSpectra{\U(1)}{\CoulombUniverse}$,
	the spectrum
	$$
		\Sigma^{-W^{(-\mu,0)}}
		\Sigma^\infty_\CoulombUniverse
			\Conley_{\U(1)}(
				S^\mu_{\lambda,r},
				\FDASWFlow^\mu_{\lambda,r}
			)
	$$
	only depends on $Y$,
	the $\SpinC$ structure $\SpinCStruct_\lambda$
	and the metric $g$.
}
A disadvantage of working with the category
$\StabHoSpectra{\U(1)}{\CoulombUniverse}$
is that the universe $\CoulombUniverse$
depends,
at least superficially,
on the metric $g$.
This can be addressed by applying a {\it change of universe\/}
to pass over to a standard choice
of universe isomorphic to $\CoulombUniverse$.
In order to demonstrate that the choices involved
are immaterial to the final result,
the author shall invoke
one more concept
from stable homotopy theory.
\Def{
	Given two $G$-universes
	$\Universe_0$ and $\Universe_1$
	and a linear isometry
	$f:\Universe_0\to\Universe_1$,
	define the associated
	{\it (restrictive) change of universe functor,}
	$$
		f^*:\StabHoSpectra G{\Universe_1}
		\to\StabHoSpectra G{\Universe_0},
	$$
	by defining it on prespectra
	as the functor that sends
	a prespectrum $E$
	with structure maps $\sigma$
	to the prespectrum
	$f^* E$
	with structure maps $f^*\sigma$
	such that
	$$
		(f^*E)_V=E_{f(V)},
		\quad
		(f^*\sigma)_{V,W}=\sigma_{f(V),f(W)}.
	$$
	Likewise,
	define the associated
	{\it (inductive) change of universe functor,}
	$$
		f_*:\StabHoSpectra G{\Universe_0}
		\to\StabHoSpectra G{\Universe_1},
	$$
	by defining it on a prespectrum $E$
	with structure maps $\sigma$
	to be the prespectrum $f_*E$
	with structure maps $f_*\sigma$
	satisfying the following.
	For finite dimensional subrepresentations
	$V\subset W\subset\Universe_1$,
	denote $V':=V\cap f(\Universe_0)$
	and $W':=W\cap f(\Universe_0)$.
	Then,
	set
	$$
		(f_*E)_V:=E_{f^{-1}(V')}\wedge(V-V')^+,
	$$
	and define the structure map
	$(f_*\sigma)_{V,W}$
	to be the composite
	$$\eqalign{
		&
		E_{f^{-1}(V')}
		\wedge
		(V-V')^+
		\wedge
		(W-V)^+
		\cr

		&
		\isomorph
		E_{f^{-1}(V')}
		\wedge
		(W'-V')^+
		\wedge
		(W-W')^+
		\cr

		&
		\morphnamed{\id\wedge f^{-1}\wedge\id}
		E_{f^{-1}(V')}
		\wedge
		(f^{-1}(W')-f^{-1}(V'))^+
		\wedge
		(W-W')^+
		\cr

		&
		\morphnamed{\sigma_{f^{-1}(V'),f^{-1}(W')}\wedge \id}
		E_{f^{-1}(W')}
		\wedge
		(W-W')^+.
	}$$
}
\Prop{(\LewisMaySteinbergerEquivariantAlt, Proposition 1.2)
	$f_*$ is left adjoint to $f^*$.
}
\Thm{
	\Label\ChangesOfUniverseCoherentlyCanoicallyIsomorphic
	(\LewisMaySteinbergerEquivariantAlt,
	Theorem 1.7 and Corollary 1.8)
	The functors
	$
		f^*:\StabHoSpectra G{\Universe_1}
		\to\StabHoSpectra G{\Universe_0}
	$
	defined by different choices of linear isometry $f$
	are {\it canonically and coherently isomorphic.}
	The same holds for the functors $f_*$.
	Moreover,
	if $f$ be an isomorphism,
	$f^*$ is an equivalence of categories
	with its inverse being $f_*$.
}
\Rem{
	The author shall not expand on the precise meaning of
	canonical nor coherent in the present \Document;
	for that,
	the reader is directed to the proof of
	{\LewisMaySteinbergerEquivariant}, Theorem 1.7.
	Suffice it to say that, by canonical, one means that
	the choices which appear in the proof are immaterial,
	and,
	by coherent,
	one means that certain diagrams that should commute
	do indeed commute.
}
\Def{
	Define $\Universe$
	to be the $\U(1)$-universe
	$\R^\infty\oplus\C^\infty$.
}
The essential point to be taken from
{\ChangesOfUniverseCoherentlyCanoicallyIsomorphic}
is that,
since one may pass between the categories
$\StabHoSpectra{\U(1)}\Universe$
and $\StabHoSpectra{\U(1)}\CoulombUniverse$
in a natural fashion,
one can think of the desuspension functor
$\Sigma^{-W^{(-\mu,0)}}$
as being defined on $\StabHoSpectra{\U(1)}\Universe$
by intertwining with the changes of universe.
\Def{
	\Label\DesuspensionIntertwinedWithChangesOfUniverse
	Given a finite dimensional subrepresentation
	$V\subset\CoulombUniverse$,
	define the desuspension functor
	$$
		\Sigma^{-V}:
		\StabHoSpectra{\U(1)}\Universe
		\to
		\StabHoSpectra{\U(1)}\Universe
	$$
	as the composite
	$$
		\StabHoSpectra{\U(1)}\Universe
		\morphnamed{f_*}
		\StabHoSpectra{\U(1)}\CoulombUniverse
		\morphnamed{\Sigma^{-V}}
		\StabHoSpectra{\U(1)}\CoulombUniverse
		\morphnamed{f^*}
		\StabHoSpectra{\U(1)}\Universe,
	$$
	where $f:\Universe\to\CoulombUniverse$
	is any isometric isomorphism
	defined by a choosing a basis of eigenvectors
	for $\ell$.
}
\Prop{
	The endofunctor $\Sigma^{-V}$
	on $\StabHoSpectra{\U(1)}\Universe$
	is well defined up to canonical isomorphism.
}
\Proof{
	Immediate from \ChangesOfUniverseCoherentlyCanoicallyIsomorphic.
}
With this understood,
henceforth,
the author shall often drop $f^*$ and $f_*$ from the notation.
\Def{
	Define the
	{\it metric dependent
	Seiberg-Witten Floer spectrum\/}
	as
	$$\eqalign{
		\SWF(Y,\SpinCStruct_\lambda,g)
		:=
		&
		\Sigma^{-W^{(-\mu,0)}}
		\Sigma^\infty_\Universe
			\Conley_{\U(1)}(
				S^\mu_{\lambda,r},
				\FDASWFlow^\mu_{\lambda,r}
			)
		\cr

		\cong
		&
		f^*
		\Sigma^{-W^{(-\mu,0)}}
		\Sigma^\infty_\CoulombUniverse
			\Conley_{\U(1)}(
				S^\mu_{\lambda,r},
				\FDASWFlow^\mu_{\lambda,r}
			)
		\in\StabHoSpectra{\U(1)}\Universe,
	}$$
	where $f:\Universe\to\CoulombUniverse$
	denotes any isometric isomorphism
	defined by a basis of eigenvectors
	for $\ell$.
}
\Cor{
	The spectrum $\SWF(Y,\SpinCStruct_\lambda,g)$
	is well defined up to canonical isomorphism in
	$\StabHoSpectra{\U(1)}\Universe$.
}
\Proof{
	Corollary of
	{\SWFWellDefWithCoulombUniverse}
	and \ChangesOfUniverseCoherentlyCanoicallyIsomorphic.
}
Next,
the author proceeds to explain how to
(de)suspend away
the metric dependence.
\Def{
	\Label\SquareInSecondCohomologyOfFourManiWithQHSThreeBoundary
	Let $X$ be an oriented $4$-manifold-with-boundary
	such that $\partial X=Y$.
	Given a class $c\in\Homology^2(X,\Z)$
	define its square,
	$c^2\in\Q$,
	as follows.
	Let $c'\in\Homology^2(X;\Q)$
	be the image of $c$
	under the change of coefficients
	$\Homology^2(X;\Z)\to\Homology^2(X;\Q)$.
	Since $b_2(Y)=0$,
	there is an exact sequence
	$$
		\Homology^2(X,Y;\Q)
		\to\Homology^2(X;\Q)
		\to 0.
	$$
	Pick any preimage $\tilde c\in\Homology^2(X,Y;\Q)$
	for $c'$.
	Define
	$$
		c^2:=(\tilde c\smile\tilde c)[X,Y]\in\Q.
	$$
}
\Def{
	\Label\DesuspensionParameter
	Define a number
	$n(Y,\SpinCStruct_\lambda,g)\in\Q$
	as follows.
	Choose some simply-connected
	$4$-manifold-with-boundary $X$
	with $\SpinC$ structure ${\frak t}$
	such that $\partial X=Y$
	and such that ${\frak t}$ agree with
	$\SpinCStruct_\lambda$ on $Y$.
	Assume further that
	$X$ have a neighbourhood of its boundary
	which be isometric to $[0,1]\times Y$.
	Fixing some connexion $B$ on
	$\det{\frak t}$
	extending arbitrarily the connexion
	$A_\lambda$ on
	$\det\SpinCStruct_\lambda$,
	use $\Dirac_{B}^\pm$
	to denote the Dirac operators of
	$(X,{\frak t})$.
	Denote the signature of $X$
	by $\sigma(X)$.
	With this notation in place,
	let
	$$
		n(Y,\SpinCStruct_\lambda,g)
		:=\index_\C\Dirac_{B}^+
		-{1\over8}\big(
			c_1(\det{\frak t})^2
			-\sigma(X)
		\big).
	$$
}
\Prop{
	\Label\AlternativeDescriptionOfDesuspensionParameter
	(\ManolescuSeibergWittenFloerAlt, \S 6)
	The number
	$n(Y,\SpinCStruct_\lambda,g)$
	does not depend on the choices
	involved in its definition.
	Indeed,
	$$
		n(Y,\SpinCStruct_\lambda,g)
		={1\over2}\left(
			\eta(\Dirac_{A_\lambda},0)
			-\dim_\R\Ker\Dirac_{A_\lambda}
			-{1\over4}\eta(\SignOp, 0)
		\right)
	$$
	where $\eta(D,z)$
	denotes the $\eta$ function
	of an operator $D$
	(vid.\ \AtiyahPatodiSingerSpectralAlt),
	and $\SignOp$
	is the operator on $\Omega^1(Y)\oplus\Omega^0(Y)$
	given by
	$$
		\pmatrix{
			\Hodge\d&-\d\cr
			-\d^*&0
		}.
	$$
}
\Prop{
	(\ManolescuSeibergWittenFloerAlt, \S 7)
	If $N$ denote the cardinality of the finite set
	$\Homology_1(Y;\Z)$,
	then
	$8Nn(Y,\SpinCStruct,g)$
	is an integer
	and its residue
	modulo $8N$
	does not depend on $g$.
}
\Def{
	Let	
	$\MetricDesuspParamRes(Y,\SpinCStruct)\in\Q$
	be such that
	$8N\MetricDesuspParamRes(Y,\SpinCStruct)\in\{0,\ldots,8N-1\}$
	be the residue modulo $8N$
	of $8Nn(Y,\SpinCStruct,g)$.
}
\Thm{
	(\ManolescuSeibergWittenFloerAlt, \S 7)
	Up canonical to isomorphism in
	$\StabHoSpectra{\U(1)}{\Universe}$,
	the spectrum
	$$
		\SWF(Y,\SpinCStruct_\lambda)
		:=
		\Sigma^{\C^{
			{\frak n}(Y,\SpinCStruct)
			-n(Y,\SpinCStruct_\lambda,g)}
		}
		\SWF(Y,\SpinCStruct_\lambda,g)
	$$
	only depends on $Y$
	and the $\SpinC$ structure $\SpinCStruct_\lambda$.
}
\Rem{
	\Label\DifferentDesuspensionConvention
	The author is deviating from what is stated in
	{\ManolescuSeibergWittenFloer} slightly
	in using the (de)sus\-pen\-sion
	by
	$
		\C^{
			{\frak n}(Y,\SpinCStruct)
			-n(Y,\SpinCStruct_\lambda,g)
		}
	$
	instead of by
	$\C^{-n(Y,\SpinCStruct_\lambda,g)}$.
	The reason being
	his desire to operate
	in the standard category of spectra,
	$\Spectra {\U(1)}\Universe$,
	which does not allow for a desuspension
	by a ``rational'' dimensional representation
	in any obvious manner.
	In {\ManolescuSeibergWittenFloer},
	a variant of the Spanier-Whitehead category
	is used,
	which can easily be made to formally admit
	such (de)suspensions.
}
The main result of {\LidmanManolescuEquivalence}
is that the $\U(1)$-equivariant Borel cohomology
of the $\SWF$ spectrum
corresponds to the classical monopole Floer cohomology.
In order to state this theorem,
the author should firstly clarify
what is meant by the cohomology of a spectrum
in the present context.
The rightful sort of $G$-equivariant cohomology theory to consider
for $G$-spectra indexed on a universe $\AnotherUniverse$
is that of {\it $\RO(G;\AnotherUniverse)$-graded cohomology\/}
(vid.\ \MayEtAlEquivariantAlt, \S XIII).
In order to avoid having to deal with such complexities,
the author decided to introduce
the following language
which shall simplify considerably the treatment.
\Def{
	Let $h^*_G$ denote a $\Z$-graded
	$G$-equivariant (reduced) cohomology theory
	on $G$-spaces.
	The author shall say that $h^*_G$ satisfies
	the
	{\it suspension axiom
	with respect to the universe $\AnotherUniverse$\/}
	when,
	for any $G$-space $X$
	and any finite dimensional subrepresentation
	$V\subset\AnotherUniverse$,
	there be a natural isomorphism
	$$
		h^n_G(X)\cong h^{n+\dim V}_G(\Sigma^V X).
	$$
}
\Def{
	Given a $G$-prespectrum $E$ indexed on a universe $\AnotherUniverse$
	and $G$-equivariant cohomology theory $h_G^*$ on $G$-spaces
	satisfying the suspension axiom
	with respect to the universe $\AnotherUniverse$,
	define the
	{\it cohomology of $E$\/}
	as
	$$
		h_G^n(E)
		:=
		\Colim_{V\subset\AnotherUniverse}h_G^{n+\dim V}(E_V).
	$$
}
\Rem{
	It is easy to see from this definition
	that,
	for a desuspension spectrum
	$\Sigma^{-V}\Sigma^\infty_\AnotherUniverse X$,
	the cohomology is simply
	given by the cohomology of the space $X$
	but with a {\it grading shift.}
	$$
		h_G^n(\Sigma^{-V}\Sigma^\infty_\AnotherUniverse X)
		=h_G^{n+\dim V}(X).
	$$
}
\Def{
	The
	{\it Borel\/}
	$G$-equivariant cohomology theory for $G$-spaces
	with coefficients in an abelian ring $R$
	is defined as follows.
	Denote by $\E G$
	a free contractible $G$-space
	and by $\B G$
	the quotient $\E G/G$.
	Then,
	for a $G$-space $X$,
	$$
		\BorelHomology^*_G(X;R)
		:=\RedHomology^*((\E G\times X)/G;R).
	$$
}
\Rem{
	Following \MayEtAlEquivariant,
	the notation $c$ indicates the ``geometric completion''
	involved in obtaining the underlying spectrum
	from the equivariant Eilenberg-MacLane spectrum
	and helps one distinguish {\it Borel\/}
	from {\it Bredon\/} cohomology.
}
\Prop{
	For a $G$-universe $\AnotherUniverse$
	containing only finite dimensional subrepresentations
	$V$
	such that the vector bundle
	$(V\times\E G)/G$
	over $\B G$
	be {\it $R$-orientable,}
	the Borel cohomology theory with $R$-coefficients
	{\it satisfies the suspension axiom for $\AnotherUniverse$.}
}
\Proof{
	This follows directly from the Thom isomorphism theorem.
}
\Rem{
	In the case at hand of $G=\U(1)$
	and universe $\Universe=\R^\infty\oplus\C^\infty$,
	it is clear,
	therefore,
	that the Borel cohomology theory
	with $\Z$-coefficients
	satisfies the suspension axiom for $\Universe$
	because a direct sum of a complex representation
	and a trivial real representation
	always define $\Z$-orientable
	vector bundles over $\B G$.
	Hence,
	one can speak of the Borel cohomology of the $\SWF$ spectrum.
}
\Thm{
	(\LidmanManolescuEquivalenceAlt, Theorem 1.2.1)
	Letting
	$\HMFrom^*(Y,\SpinCStruct)$
	denote the $\Q$-graded
	``from'' monopole Floer cohomology
	of {\KronheimerMrowkaMonopolesAndThreeManifolds},
	there is an isomorphism
	$$
		\HMFrom^q(Y,\SpinCStruct)
		\cong
		\BorelHomology_{\U(1)}^{
			q
			-{\frak n}(Y,\SpinCStruct)
		}
		(\SWF(Y,\SpinCStruct);\Z).
	$$
	In particular,
	letting the ``tilde'' monopole Floer cohomology,
	$\HMTilde^*(Y,\SpinCStruct)$,
	be the mapping cone of the $U$-map
	on $\HMFrom^*(Y,\SpinCStruct)$,
	there is an isomorphism
	$$
		\HMTilde^q(Y,\SpinCStruct)
		\cong
		\RedHomology^{
			q
			-{\frak n}(Y,\SpinCStruct)
		}
		(\SWF(Y,\SpinCStruct);\Z).
	$$
}
\Def{
	For future reference,
	define the
	{\it metric dependent tilde monopole Floer cohomology\/}
	to be
	$$
		\HMTilde^*(
			X,
			\SpinCStruct,
			g
		)
		:=\Homology^*(
			\SWF(X,\SpinCStruct,g);
			\Z
		).
	$$
}

\Sect{The Cohomotopical Contact Invariant}
This section shall fulfil the first goal of the present \Document.
The analytic results concerning the contact monopole
$C_\lambda$
shall allow the author to define a
{\it cohomotopical\/} contact invariant
in a manner that can be roughly outlined in the following way.
Given a top dimensional cell in a CW-complex,
if one were to quotient the complex by all other cells,
one obtains a map to a sphere;
this is an element of the cohomotopy of the complex.
In the event that the cell not be top dimensional,
but,
instead,
all the higher dimensional cells
attach {\it null-homotopically\/} onto the given cell,
one can still perform the same quotient
and obtain an element in the cohomotopy set of the complex.
In the case of a $G$-CW-complex,
a similar story can be told about a $G$-cell;
here,
one needs to be more careful
with what is meant by cohomotopy.
In any event,
one obtains a map from the $G$-CW-complex
to the Thom space of a vector bundle
over a $G$-orbit;
analogously to how,
in the non-equivariant setting,
the sphere is the Thom space of a bundle over
the orbit of the trivial group;
that is,
the point.
\par
To achieve this in the present context,
the author shall make use of a fundamental construction
in the Conley theory;
namely,
the notion of {\it attractor-repeller pairs.}
What has already been said
about the contact configuration
shall be summarised as saying
that the orbit of the contact configuration
defines a {\it repeller\/} in the isolated invariant set
defining the Seiberg-Witten Floer spectrum.
Well known results on Conley theory
then provide a {\it cofibration\/}
involving Conley indices.
Due to the non-degeneracy of the contact configuration,
the cofibre map of this cofibration can be interpreted,
in the presence of a choice of $\U(1)$-CW-structure
on the Seiberg-Witten Floer spectrum,
exactly as the quotient of all but one special $\U(1)$-cell
defined by the orbit of the contact monopole.
This cofibre map shall be declared the contact invariant.
\par
Lastly,
one must take care to {\it stabilize\/}
everything so as to make the cofibre map really an invariant
with respect to the spectral cut-off parameter $\mu$
and the metric $g$.
As it turns out,
this does not provide any difficulty
beyond what was already encountered in
\ManolescuSeibergWittenFloer,
and the proofs shall follow closely
what is said there
only with a few extra Conley theoretic inputs
concerning the naturality
of attractor-repeller pair cofibrations.
\par
Let $Y$, $\lambda$ and $g$
be as in the preceding sections.
\Rem{\Label\ContactConfInCoulombBall
	Notice that
	$\ContactConf_\lambda$
	is,
	by definition,
	in the global Coulomb slice
	with respect to $A_\lambda$;
	that is,
	$C_\lambda\in W$.
	Moreover,
	since $\Dirac_{A_\lambda}\psi_\lambda=0$,
	it is also true that
	$C_\lambda\in\Ker\ell$.
	Therefore,
	for $R>0$ sufficiently large,
	$\ContactConf_\lambda\in \Ball(W^\mu,R)$
	for all $\mu>0$.
	Henceforth,
	agree to set
	$R>0$
	large enough so that this last inclusion hold.
}
\Rem{
	It is worth emphasising that the global Coulomb gauge
	does not fix a gauge with respect to
	the entire gauge group $\Gauge(Y)$
	but,
	rather,
	with respect to the normed gauge group
	$\NormGauge(Y)$.
	Hence,
	there is a circle's worth of
	fixed points of $\GCSW_{\lambda,r}$
	in $W_k$
	which are gauge equivalent to
	$\ContactConf_\lambda$.
}
\Def{
	Denote by 
	$\ContactCirc_\lambda\subset W_k$
	the circle of configurations
	gauge equivalent to $\ContactConf_\lambda$
	and call it the
	{\it contact circle.}
}
\Def{
	Use $\GCGaugeOrbitSlice_C$
	to denote the tangent space to the
	$\U(1)$-orbit
	at $C\in W_k$.
}
\Def{
	\Label\GTildeMetric
	Let $\tilde g$
	denote the metric on $W_k$
	defined by
	assigning to tangent vectors
	$(a,\psi),(b,\phi)\in \Tan_C W_k\cong W_k$
	the value
	$
		\Re\left\langle
			\ELCProj_C(a,\psi),
			\ELCProj_C(b,\phi)
		\right\rangle
	$,
	where $\ELCProj_C$
	is the projection to the enlarged local Coulomb gauge
	defined in
	\EnlargedLocalCoulombProjection.
}
\Rem{
	The metric $\tilde g$
	is the one used in many of the technical results
	of {\LidmanManolescuEquivalence}.
	It is notable because it turns the
	Seiberg-Witten vector field $\GCSW$
	in global Coulomb gauge
	into the $\tilde g$-gradient
	of the $\CSD_{\lambda,r}$ functional
	restricted to $W_k$.
	In the present situation,
	it shall be necessary to invoke
	some of those results of {\LidmanManolescuEquivalence}
	which make reference to this metric.
	The $\tilde g$ metric leads to the following
	definition.
}
\Def{
	\Label\AGCSliceDef
	Define the
	{\it local anticircular slice 
	in the global Coulomb gauge\/}
	at $C\in W_k$,
	denoted $\AGCSlice_C$,
	as the $\tilde g$-orthogonal complement
	to $\GCGaugeOrbitSlice_C$
	in $W_k$.
	Use $\AGCSlice_{j,C}$
	for its Sobolev completion of regularity $j\in\Z$
	and use $\AGCProj_C:W_k\to\AGCSlice_C$
	to denote the $\tilde g$-orthogonal projection.
}
\Prop{\Label\ContactCircNonDegen
	For sufficiently large $r>0$,
	at any $C\in\ContactCirc_\lambda$,
	the derivative
	$$
		\AGCProj_C
		\circ
		\D_C\SW_{\lambda,r}
		\circ
		\ELCProj_C:
		\AGCSlice_{k,C}
		\to \AGCSlice_{k-1,C}
	$$
	is surjective.
}
\Proof{
	By gauge equivariance,
	it suffices to prove the result
	for $C=\ContactConf_\lambda$.
	{\ContactConfNonDegen}
	has established that,
	for sufficiently large $r>0$,
	the map
	$$
		\LCProj_{\ContactConf_\lambda}
		\circ\D_{\ContactConf_\lambda}\SW:
		\LCSlice_{\ContactConf_\lambda}
		\to\LCSlice_{\SW_{\lambda,r}
		(\ContactConf_\lambda)}
	$$
	is surjective.
	Hence,
	the result follows directly from
	\LidmanManolescuEquivalence,
	Lemma 5.6.1,
	``$\hbox{(ii)}\Rightarrow\hbox{(iv)}$'';
	cf.\ \LidmanManolescuEquivalence,
	Forumulae (97) and (94).
}
\Rem{
	Henceforth,
	fix $r>0$
	so as to make
	{\ContactCircNonDegen}
	hold.
}
\Rem{
	Bear in mind that the vector space
	$W^\mu$
	is the direct sum of a real vector space
	and a complex vector space,
	whence comes its $\U(1)$-action.
	Note that $\ContactCirc_\lambda\incl W^\mu$
	is the $\U(1)$-equivariant
	embedding of a $\U(1)$-manifold.
}
\Def{
	Let
	$\ContactNormal\mu\lambda\to\ContactCirc_\lambda$
	denote the $\U(1)$-equivariant
	normal bundle of
	$\ContactCirc_\lambda$
	as a submanifold of $W^\mu$.
}
\Prop{\Label\FDAContactCircNonDegen
	For sufficiently large $\mu>0$,
	$\ContactCirc_\lambda$
	is a hyperbolic fixed set
	of the flow $\FDASWFlow^\mu_{\lambda,r}$
	in $W^\mu$.
	In other words,
	for any $C\in\ContactCirc_\lambda$,
	the derivative
	$$
		\D_C\FDASW_{\lambda,r}^\mu:
		\ContactNormal\mu\lambda|_C
		\to
		\ContactNormal\mu\lambda|_C
	$$
	has no eigenvalue with vanishing real part.
	In particular,
	$\ContactCirc_\lambda$
	is a non-degenerate fixed set.
}
\Proof{
	By \ContactCircNonDegen,
	$\ContactCirc_\lambda$
	consists of non-degenerate irreducible fixed points
	of the flow $\GCSW_{\lambda,r}$
	on $W_k$.
	Hence,
	apply the same argument
	used in the proof of
	{\LidmanManolescuEquivalence},
	Proposition 7.3.1,
	to find that the same remains true
	when passing to a finite dimensional approximation
	provided one choose a sufficiently large $\mu$.
}
\Rem{
	Fix $\mu>0$ large enough
	so as to make
	{\FDAContactCircNonDegen}
	hold.
}
\Rem{\Label\IntroUnstableBundle
	Identify
	$\ContactNormal\mu\lambda$
	with a sufficiently small tubular neighbourhood
	of $\ContactCirc_\lambda$
	so as to not contain
	any other fixed points
	of the flow $\FDASWFlow_{\lambda,r}^\mu$.
	This is possible due to the non-degeneracy
	ensured by {\FDAContactCircNonDegen}.
	Now,
	as a vector bundle,
	one can split
	$\ContactNormal\mu\lambda$
	into stable and unstable subbundles as
	$
		\ContactNormalSt\mu\lambda
		\oplus
		\ContactNormalUnst\mu\lambda$,
	where
	$$
		\eqalign{
			\forall v\in
			\ContactNormalSt\mu\lambda,
			&\quad
			\left\langle
				v,
				\D_C\FDASW_{\lambda,r}^\mu (v)
			\right\rangle
			\geq m|v|^2,
			\cr
			\forall v\in
			\ContactNormalUnst\mu\lambda,
			&\quad
			\left\langle
				v,
				\D_C\FDASW_{\lambda,r}^\mu (v)
			\right\rangle
			\leq -m|v|^2
		}
	$$
	for some constant $m>0$.
	Moreover,
	this splitting is preserved by
	$\D_C\FDASW_{\lambda,r}^\mu$
	and,
	hence,
	$\ContactCirc_\lambda$
	is an isolated invariant set
	with index pair
	$
		(\Disk(\ContactNormal\mu\lambda),
		\Sphere(\ContactNormalUnst\mu\lambda))
	$
	where $\Disk$ and $\Sphere$
	denote the unit disk and unit sphere bundles.
}
\Def{
	For $E$ a vector $G$-bundle
	over a compact Hausdorff space,
	denote its $G$-equivariant Thom space by
	$\Thom_G(E):=\Disk E/\Sphere E$.
}
\Cor{\Label\ConleyContactCircIsThomSpace
	The Conley index
	$
		\Conley_{\U(1)}(
			\ContactCirc_\lambda,
			\FDASWFlow^\mu_{\lambda,r}
		)
	$
	is $\U(1)$-equivariantly homotopy equivalent 
	to the $\U(1)$-equivariant Thom space
	$\Thom_{\U(1)}(\ContactNormalUnst\mu\lambda)$
	of the bundle
	$\ContactNormalUnst\mu\lambda\to\ContactCirc_\lambda$.
}
As a direct consequence,
by using a Morse decomposition,
one finds that
the contact circle $\ContactCirc_\lambda$
defines a $\U(1)$-cell
in some $\U(1)$-CW-complex decomposition
of the space
$
	\Conley_{\U(1)}(
		S^\mu_{\lambda,r},
	\FDASWFlow^\mu_{\lambda,r})
$.
If one were to add a {\it generic\/} perturbation
on top of the canonical perturbation in use,
one would find that a possible choice of
$\U(1)$-CW-decomposition
would have a one to one correspondence
between its $\U(1)$-cells
and the set of fixed points
and fixed circles
of the Seiberg-Witten flow.
It is preferable,
of course,
to avoid using such a generic perturbation.
In order to derive a cohomotopical invariant
from a given cell,
it is necessary to establish something about the
attaching maps of the higher dimensional cells.
In particular,
if the given cell be itself top dimensional
or,
more generally,
if one find that
all higher dimensional cells attach null-homotopically
onto the given cell,
then it follows that
one can quotient all but the given cell
in the $\U(1)$-CW-complex
and obtain a map to the quotient
of the contact cell by its boundary.
This is morally the strategy 
which shall be pursued next.
\Thm{
	\Label\ContactCircRepeller
	For sufficiently large $r>0$
	and $\mu>0$,
	there are no non-constant
	approximate Seiberg-Witten trajectories
	in the set $S^\mu_{\lambda,r}$
	which have plus infinity limit
	in the contact circle $\ContactCirc_\lambda$.
}
\Proof{
	The author starts with the argument
	in Step 1 of the proof of Proposition 3 in
	{\ManolescuSeibergWittenFloer}.
	Suppose the result not hold.
	Then,
	there exists an increasing sequence
	$\mu_n\to\infty$
	and a sequence of non-constant trajectories
	$
		\gamma_n:
		\R
		\to
		\Disk(W^{\mu_n},2R)
	$
	satisfying
	$$
		\pdbyd t\gamma_n(t)
		=-\FDASW^{\mu_n}_{\lambda,r}(\gamma_n(t)),
		\quad
		\lim_{t\to\infty}\gamma_n(t)
		=\ContactConf_\lambda.
	$$
	Notice that there is a bound
	$$
		\norm{\pdbyd t\gamma_n(t)}_{\Sobolev k-1,2}
		\leq
		\norm{\ell\gamma_n(t)}_{\Sobolev k-1,2}
		+\norm{p^\mu c \gamma_n(t)}_{\Sobolev k-1,2}
		\leq KR.
	$$
	where $K>0$ is a constant
	independent of both $n$ and $t$.
	This implies that the set of functions
	$$
		\left\{
			\gamma_n:
			\R
			\to
			\Disk(W_{k-1},2R)
		\right\}
	$$
	is equicontinuous.
	By use of the Arzel\`a-Ascoli theorem,
	one can replace this sequence
	by a subsequence which converge
	to some
	$\gamma:\R\to\Disk(W_{k-1},2R)$
	in the compact-open topology.
	Now,
	due to compactness of $c$,
	the sequence of operators
	$(1-p^{\mu_n})c:W_k\to W_{k-1}$
	converges to zero weakly.
	Given this,
	observe that
	$$
		\pdbyd t\gamma_n
		=-(\ell+c)\gamma_n
		+(1-p^{\mu_n})c\gamma_n
		\converge{n\to\infty}
		-(\ell+c)\gamma.
	$$
	uniformly as functions
	from compact subsets of $\R$
	to $W_{k-1}$.
	Hence,
	$$
		\gamma_n(t)
		-\gamma_n(0)
		=\int_0^t \pdbyd s\gamma_n(s)\d s
		\converge{n\to\infty}
		-\int_0^t (\ell+c)\gamma(s)\d s.
	$$
	Together,
	these two last assertions imply that
	$$
		\pdbyd t\gamma(t)
		=-(\ell+c)\gamma(t).
	$$
	Moreover,
	observe that
	$\gamma(t)\converge{t\to\infty}U_\lambda$.
	Therefore,
	$\gamma$
	is the Coulomb projection
	of a Seiberg-Witten trajectory
	with positive infinity limit
	gauge equivalent to the contact configuration.
	By {\NoTrajToContactConf},
	such a trajectory must be constant.
	Let $C_n:=\lim_{t\to-\infty}\gamma_n(t)$.
	By assumption,
	these cannot be equal to $C_\lambda$.
	Note also that
	$\lim_{n\to\infty} C_n =\lim_{t\to-\infty}\gamma(t)=C_\lambda$.
	However,
	{\ContactCircNonDegen}
	combined with
	\LidmanManolescuEquivalence,
	Proposition 7.2.2,
	guarantee that,
	for sufficiently large $n$
	and some small neighbourhood
	$N\supset W_k$
	of $C_\lambda$,
	there cannot be any solution
	to $\FDASW_{\lambda,r}^\mu=0$
	inside of $N$
	other than $C_\lambda$ itself,
	thereby
	contradicting
	convergence of the sequence
	$\{C_n\}$ to $C_\lambda$.
}
A few more concepts from Conley theory
shall be needed next.
In what follows,
suppose that $G$ be a compact Lie group,
$\Gamma$ be a locally compact Hausdorff space
with a continuous $G$-action
and $\phi:\Gamma\times\R\to\Gamma$
a continuous and equivariant flow.
\Def{
	\Label\OmegaAndOmegaStarLimits
	For a $G$-invariant subset
	$T\subset\Gamma$,
	define its
	{\it $\omega$-limit,}
	$\omega(T)$,
	as the maximal invariant set
	of the closure of
	$\varphi(T\times[0,\infty))$.
	Likewise,
	define its 
	{\it $\omega^*$-limit,}
	$\omega^*(T)$,
	as the maximal invariant set
	of the closure of
	$\varphi(T\times(-\infty,0])$.
}
\Def{
	A $G$-invariant subset
	$A\subset S$
	is an {\it attractor\/}
	when there is a neighbourhood
	$U\subset S$ of $A$ as a subspace of $S$
	such that $A=\omega(U)$.
	Likewise,
	a $G$-invariant subset
	$R\subset S$
	of an isolated invariant set $S\subset\Gamma$
	is called an {\it repeller\/}
	when there is a neighbourhood
	$U\subset S$ of $R$ as a subspace of $S$
	such that $R=\omega^*(U)$.
}
\Def{
	Given an attractor $A\subset S$
	in the isolated invariant set $S\subset\Gamma$,
	one defines its {\it complementary repeller,}
	$A^*\subset S$,
	as the set
	$\{x\in S\mid \omega(\{x\})\cap A=\emptyset\}$.
	Given a repeller $R\subset S$,
	one defines its {\it complementary attractor,}
	$R^*\subset S$,
	similarly.
}
\Def{
	\Label\AttractorRepellerPairDef
	An {\it attractor-repeller pair\/} $(A,R)$
	of an isolated invariant set $S\subset\Gamma$
	consists of an attractor $A$
	and a repeller $R$ in $S$
	such that $A=R^*$,
	or,
	equivalently,
	$R=A^*$.
}
\Def{
	\Label\IndexTripleDef
	Let $(A,R)$
	be an attractor-repeller pair for
	of an isolated invariant set $S\subset\Gamma$.
	A triple of $G$-invariant compact sets $(L,M,N)$,
	$N\subset M\subset L\subset\Gamma$,
	is called an {\it index triple\/}
	for $(A,R)$
	when
	$(L,M)$ is an index pair for $R$,
	$(L,N)$ is an index pair for $S$
	and $(M,N)$ is an index pair for $A$.
}
\Thm{
	\Label\IndexTripleExists
	(\ConleyIsolatedAlt, \S I.7, albeit non-equivariantly)
	\Label\IndexTripleExistsCofibWellDef
	For any attractor-repeller pair $(A,R)$
	of an isolated invariant set $S\subset\Gamma$,
	there exists an index triple $(L,M,N)$ for it
	and the induced {\it cofibration,}
	$$
		\Conley_G(A,\phi)
		\to\Conley_G(S,\phi)\to
		\Conley_G(R,\phi),
	$$
	is independent of the choice of $(L,M,N)$
	up to $G$-equivariant homotopy.
}
Now,
one can reinterpret the analytic results
concerning the contact circle
in the Conley theoretic language.
\Thm{
	\Label\ContactCellCofibration
	The contact circle
	$\ContactCirc_\lambda\subset S^\mu_{\lambda,r}$
	is a
	{\it repeller\/}
	in $S^\mu_{\lambda,r}$.
	Hence,
	{\it there exists a cofibration\/}
	$$
		\Conley_{\U(1)}(
			\ContactCirc_\lambda^*,
			\FDASWFlow_{\lambda,r}^\mu
		)
		\to\Conley_{\U(1)}(
			S^\mu_{\lambda,r},
			\FDASWFlow_{\lambda,r}^\mu
		)
		\to\Conley_{\U(1)}(
			\ContactCirc_\lambda,
			\FDASWFlow_{\lambda,r}^\mu
		).
	$$
}
\Proof{
	Follows from
	{\ContactCircRepeller}
	and
	\IndexTripleExistsCofibWellDef.
}
\Def{
	Define the
	{\it spectral cut-off
	and metric dependent
	cohomotopical contact invariant\/}
	as the cofibre map
	$$
		\Psi(\lambda,g,\mu):
		\Conley_{\U(1)}(
			S^\mu_{\lambda,r},
			\FDASWFlow_{\lambda,r}^\mu
		)
		\to\Conley_{\U(1)}(
			\ContactCirc_\lambda,
			\FDASWFlow_{\lambda,r}^\mu
		)
	$$
	of \ContactCellCofibration.
}
\Prop{
	$\Psi(\lambda,g,\mu)$
	does not depend on the choice of $r>0$.
}
\Proof{
	Suppose one chose two
	sufficiently large values
	for $r$,
	call them $r_0<r_1$.
	Without loss of generality,
	assume $|r_0-r_1|$
	to be small.
	Pick the value of $R>0$
	large enough so that
	it satisfy
	{\FiniteTypeBounded}
	for both $r=r_0$ and $r=r_1$.
	Then,
	$\Disk(W^\mu,2R)$
	serves as an isolating neighbourhood
	for all the isolated invariant sets
	$S^\mu_{\lambda,r}$
	under the parametrised family of flows
	$\FDASWFlow^\mu_{\lambda,r}$
	as $r$ varies in $[r_0,r_1]$.
	The continuation properties of Conley theory
	then provide
	a diagram of the form
	$$
		\Diagram{
			\Conley_{\U(1)}(
				S^\mu_{\lambda,r_0},
				\FDASWFlow_{\lambda,r_0}^\mu
			)
			&
			\Right
			&
			\Conley_{\U(1)}(
				\ContactCirc_\lambda,
				\FDASWFlow_{\lambda,r_0}^\mu
			)
			\cr
			\Down
			&
			&
			\Down
			\cr
			\Conley_{\U(1)}(
				S^\mu_{\lambda,r_1},
				\FDASWFlow_{\lambda,r_1}^\mu
			)
			&
			\Right
			&
			\Conley_{\U(1)}(
				\ContactCirc_\lambda,
				\FDASWFlow_{\lambda,r_1}^\mu
			)
			\cr
		}
	$$
	which commutes up to homotopy
	and where the vertical arrows
	are homotopy equivalences
	(cf.\ \ConleyIsolatedAlt, \S III.3.1
	and {\KurlandHomotopyAlt}
	for the non-equivariant case).
}
\Def{
	\Label\MetricDepContactThomIntro
	The
	{\it metric dependent contact Thom space\/}
	is the desuspension
	$$
		\ContactThom(\lambda,g)
		:=
		\Sigma^{-W^{(-\mu,0)}}
		\Sigma^\infty_\Universe
		\Conley_{\U(1)}(
			\ContactCirc_\lambda,
			\FDASWFlow_{\lambda,r}^\mu
		)
		\in\StabHoSpectra{\U(1)}\Universe.
	$$
}
\Def{
	The
	{\it metric dependent
	cohomotopical contact invariant,}
	$$
		\Psi(\lambda,g):
		\SWF(Y,\SpinCStruct_\lambda,g)
		\to
		\ContactThom(\lambda,g).
	$$
	is the desuspension
	$\Sigma^{-W^{(-\mu,0)}}\Psi(\lambda,g,\mu)$
	as a morphism in
	$\StabHoSpectra{\U(1)}\Universe$.
}
\Prop{
	$\Psi(\lambda,g)$
	does not depend on the choice of $\mu>0$.
}
\Proof{
	This proof shall start
	by recalling the setup
	of the first half of the proof of
	Theorem 1 of \ManolescuSeibergWittenFloer.
	Assume two values for $\mu$ be given;
	call them $0<\mu_0<\mu_1$
	and assume both be large enough
	so as to satisfy
	{\FDCompactness}.
	Consider $\mu\in[\mu_0,\mu_1]$.
	Denote $\tilde\FDASWFlow_{\lambda,r}^\mu$
	the flow of the vector field
	$-u^{\mu_1}\cdot(\ell+p^\mu cp^\mu)$
	on $W^{\mu_1}$.
	It is easy to check that,
	for any $\mu\in[\mu_0,\mu_1]$,
	all finite type trajectories of
	$-(\ell+p^\mu cp^\mu)$
	on $W^{\mu_1}$
	are,
	in fact,
	contained in $W^\mu\subset W^{\mu_1}$.
	As a consequence,
	notice that the set
	$S^\mu_{\lambda,r}\subset W^\mu$
	can be identified with
	the union of the
	finite type trajectories of
	$\tilde\FDASWFlow_{\lambda,r}^\mu$
	contained in
	$\Ball(W^{\mu_1},R)$.
	Hence,
	abuse notation and write
	$S^\mu_{\lambda,r}\subset W^{\mu_1}$
	for all $\mu\in[\mu_0,\mu_1]$.
	Moreover,
	in this description,
	$\Disk(W^{\mu_1},2R)$
	is an isolating neighbourhood for
	$S^\mu_{\lambda,r}$
	for all $\mu\in[\mu_0,\mu_1]$.
	By the continuation properties
	of the Conley index,
	$$
		\Conley_{\U(1)}(
			S^{\mu_1},
			\FDASWFlow^\mu_{\lambda,r}
		)\simeq
		\Conley_{\U(1)}(
			S^{\mu_0},
			\tilde\FDASWFlow^\mu_{\lambda,r}
		).
	$$
	Write $W^{\mu_1}=W^{\mu_0}\oplus W'$,
	where $W'$
	is $\LTwo$-orthogonal to $W^{\mu_0}$.
	Of course,
	$W'$ is simply the span
	of the eigenvectors of $\ell$
	with eigenvalues in $(-\mu_1,-\mu_0]\cup[\mu_0,\mu_1)$.
	Use $D\subset W'$
	to denote a small disk around the origin.
	Then,
	if one care to check,
	one finds that
	$$
		\Disk(W^{\mu_0},3R/2)
		\times D
		\subset W^{\mu_1}
	$$
	is also an isolating neighbourhood
	for $S^{\mu_0}_{\lambda,r}$.
	Furthermore,
	with respect to this product,
	one can show that the flow
	$\tilde\FDASWFlow^\mu_{\lambda,r}$
	is homotopic to a product flow
	$\FDASWFlow^{\mu_0}\times{\runes f}$,
	where ${\runes f}$
	is the flow on $D$
	induced by $-\ell$.
	Notice that
	\def\DesuspSphere{\left(W^{(-\mu_1,-\mu_0]}\right)^+}

	$$
		\Conley_{\U(1)}(\{0\},{\runes f})
		=\DesuspSphere.
	$$
	By the behaviour of the Conley index
	under product flows,
	it follows that
	$$
		\Conley_{\U(1)}(
			S^{\mu_1}_{\lambda,r},
			\FDASWFlow_{\lambda,r}^{\mu_1}
		)
		\simeq
		\Conley_{\U(1)}(
			S^{\mu_0}_{\lambda,r},
			\FDASWFlow_{\lambda,r}^{\mu_0}
		)
		\wedge
		\DesuspSphere.
	$$
	Now,
	focus is turned to the contact circle.
	Recall
	the assumption
	from
	{\ContactConfInCoulombBall}
	that
	$\mu_0,\mu_1$
	be large enough
	so that
	$
		\ContactCirc_\lambda
		\subset\Ball(W^{\mu_0},R)
		\subset\Ball(W^{\mu_1},R)
	$.
	As above,
	observe that
	$$
		\Conley_{\U(1)}(
			\ContactCirc_\lambda,
			\FDASWFlow^{\mu_1}_{\lambda,r}
		)
		\simeq
		\Conley_{\U(1)}(
			\ContactCirc_\lambda,
			\tilde\FDASWFlow^{\mu_0}_{\lambda,r}
		).
	$$
	On the other hand,
	considering again the product flow
	$\FDASWFlow^{\mu_0}_{\lambda,r}\times{\runes f}$,
	one finds
	$$
		\Conley_{\U(1)}(
			\ContactCirc_\lambda,
			\FDASWFlow^{\mu_1}_{\lambda,r}
		)
		\simeq
		\Conley_{\U(1)}(
			\ContactCirc_\lambda,
			\FDASWFlow^{\mu_0}_{\lambda,r}
		)
		\wedge
		\DesuspSphere.
	$$
	This leads to a diagram of the form
	$$
		\Diagram{
			\Conley_{\U(1)}(
				S^{\mu_0}_{\lambda,r},
				\FDASWFlow_{\lambda,r}^{\mu_0}
			)
			\wedge
			\DesuspSphere
			&
			\Right
			&
			\Conley_{\U(1)}(
				\ContactCirc_\lambda,
				\FDASWFlow_{\lambda,r}^{\mu_0}
			)
			\wedge
			\DesuspSphere
			\cr
			\Down
			&
			&
			\Down
			\cr
			\Conley_{\U(1)}(
				S^{\mu_1}_{\lambda,r},
				\FDASWFlow_{\lambda,r}^{\mu_1}
			)
			&
			\Right
			&
			\Conley_{\U(1)}(
				\ContactCirc_\lambda,
				\FDASWFlow_{\lambda,r}^{\mu_1}
			)
			\cr
		}
	$$
	which,
	due to the naturality
	of the attractor-repeller cofibration
	under continuation,
	commutes up to homotopy.
	Desuspending everything as needed
	provides the required invariance.
}
\Def{
	\Label\ContactThomIntro
	Define the
	{\it contact Thom space\/}
	as the (de)suspension
	$$
		\ContactThom(\lambda)
		:=\Sigma^{\C^{
			\MetricDesuspParamRes(Y,\SpinCStruct_\lambda)
			-n(Y,\SpinCStruct_\lambda,g)
		}}
		\ContactThom(\lambda,g)
		\in\StabHoSpectra{\U(1)}{\Universe}.
	$$
}
\Def{
	\Label\ContactCohomotopyInvariant
	Define the
	{\it cohomotopical contact invariant\/}
	of the contact rational homology sphere
	$(Y,\lambda)$,
	$$
		\Psi(\lambda):
		\SWF(Y,\SpinCStruct_\lambda)
		\to
		\ContactThom(\lambda),
	$$
	to be the (de)suspension
	$
		\Sigma^{\C^{
			\MetricDesuspParamRes(Y,\SpinCStruct_\lambda)
			-n(Y,\SpinCStruct_\lambda,g)
		}}
		\Psi(\lambda,g)
	$
	as a morphism of
	$\StabHoSpectra{\U(1)}{\Universe}$.
}
\Prop{
	\Label\CohomotopicalInvIndependentOfMetric
	$\Psi(\lambda)$
	does not depend on the choice of metric $g$.
}
\Proof{
	Again,
	the proof shall start by recalling
	what is said by
	{\ManolescuSeibergWittenFloer}
	and then extending the argument
	to deal with the contact circle.
	Since the space of compatible metrics
	is connected,
	it suffices to prove the result 
	for nearby metrics,
	so consider two such metrics
	$g_0,g_1$
	and a smooth path $t\mapsto g_t$
	interpolating them.
	One can choose $\mu>0$ and $R>0$
	large enough
	so as to satisfy the usual requirements
	for all metrics along the path $g_t$.
	The author shall use a subscripted $t$
	to denote the versions of objects
	constructed with the metric $g_t$;
	for example,
	$W^\mu_t$ denotes the version of $W^\mu$
	constructed with $g_t$.
	Notice that,
	perhaps after increasing $\mu$ slightly,
	one can assume that
	$\mu$ not be an eigenvalue of
	$\ell_t$ for any $t\in[0,1]$.
	As a consequence,
	the spaces $W^\mu_t$
	have constant dimension
	as $t$ varies.
	This means the spaces $W^\mu_t$
	form a vector bundle over $[0,1]$
	and,
	therefore,
	the vector spaces
	$W^\mu_t$
	for different values of $t$
	may be identified
	via a trivialisation of this bundle;
	hence,
	use $W^\mu$
	to denote any of the spaces
	$W^\mu_t$.
	For different values of $t$,
	consider the balls $\Ball(W^\mu,R)_t$
	all as subsets of this same
	$W^\mu$.
	By assuming the metrics $g_0,g_1$
	sufficiently close to one another,
	one also finds that,
	for any $t_1,t_2\in[0,1]$,
	$\Ball(W^\mu,R)_{t_1}\subset\Ball(W^\mu,2R)_{t_2}$.
	From this,
	it follows that
	$$
		\bigcap_{t\in[0,1]}\Disk(W^\mu,2R)_t
	$$
	is an isolating neighbourhood
	for $(S^\mu_{\lambda,r})_t$
	with respect to the flow
	$(\FDASWFlow^\mu_{\lambda,r})_t$
	for all $t\in[0,1]$.
	The flow 
	$(\FDASWFlow^\mu_{\lambda,r})_t$
	varies continuously with $t\in[0,1]$;
	hence,
	by Conley theory,
	$$
		\Conley_{\U(1)}\big(
			(S^\mu_{\lambda,r})_0,
			(\FDASWFlow^\mu_{\lambda,r})_0
		\big)
		\simeq
		\Conley_{\U(1)}\big(
			(S^\mu_{\lambda,r})_1,
			(\FDASWFlow^\mu_{\lambda,r})_1
		\big).
		\eqno{(*)}
	$$
	Consider now the contact circle.
	One easily checks that,
	under the change of metric,
	the contact circle,
	$(\ContactCirc_\lambda)_t$,
	moves smoothly in $W^\mu$
	(it is not fixed under changes of metric
	because it depends on the global Coulomb projection).
	By assuming the metrics to be sufficiently close,
	one also sees that
	$$
		\bigcap_{t\in[0,1]}
		\Disk((\ContactNormal\mu\lambda)_t)
	$$
	is an isolating neighbourhood for all
	$(\ContactCirc_\lambda)_t$,
	where
	$(\ContactNormal\mu\lambda)_t\to(\ContactCirc_\lambda)_t$
	is the $g_t$ version
	of the tubular neighbourhood
	introduced in {\IntroUnstableBundle}.
	Therefore,
	$$
		\Conley_{\U(1)}\big(
			(\ContactCirc_\lambda)_0,
			(\FDASWFlow^\mu_{\lambda,r})_0
		\big)
		\simeq
		\Conley_{\U(1)}\big(
			(\ContactCirc_\lambda)_1,
			(\FDASWFlow^\mu_{\lambda,r})_1
		\big),
	$$
	which is not much to say
	due to the characterisation
	of these as Thom spaces
	in {\ConleyContactCircIsThomSpace};
	however,
	the fact that both this homotopy equivalence
	and the one of $(*)$
	come from the same
	deformation of the flow
	allows one to use
	the naturality of the
	attractor-repeller cofibration sequence
	so as to have the following diagram
	commute up to homotopy.
	$$
		\Diagram{
			\Conley_{\U(1)}\big(
				(S^\mu_{\lambda,r})_0,
				(\FDASWFlow_{\lambda,r}^\mu)_0
			\big)
			&
			\Right
			&
			\Conley_{\U(1)}\big(
				(\ContactCirc_\lambda)_0,
				(\FDASWFlow_{\lambda,r_0}^\mu)_0
			\big)
			\cr
			\Down
			&
			&
			\Down
			\cr
			\Conley_{\U(1)}\big(
				(S^\mu_{\lambda,r})_1,
				(\FDASWFlow_{\lambda,r}^\mu)_1
			\big)
			&
			\Right
			&
			\Conley_{\U(1)}\big(
				(\ContactCirc_\lambda)_1,
				(\FDASWFlow_{\lambda,r_1}^\mu)_1
			\big).
			\cr
		}
	$$
	Now,
	the effect of the desuspensions shall be addressed.
	Consider the subspace
	$(W^{(-\mu,0)})_t\subset W^\mu$.
	Note that,
	despite all the $(W^\mu)_t\cong W^\mu$
	being identified,
	the subspaces
	$(W^{(-\mu,0)})_t$
	shall still vary with $t\in[0,1]$.
	Recall that $W$
	is defined as the direct sum
	of a real and a complex space;
	use $n_{\mu,t}$
	to denote the complex dimension
	of this complex summand 
	appearing in $(W^{(-\mu,0)})_t$.
	Notice that
	$n_{\mu,1}-n_{\mu,0}$
	is the spectral flow
	of the family of Dirac operators
	$(\Dirac_{A_\lambda})_t$
	defined as the metric $g_t$ varies.
	One can then check
	with {\AlternativeDescriptionOfDesuspensionParameter}
	that
	$$
		n(Y,\SpinCStruct_\lambda,g_0)
		-n(Y,\SpinCStruct_\lambda,g_1)
		=n_{\lambda,1}
		-n_{\lambda,0}.
	$$
	Without loss of generality,
	suppose
	$
		n(Y,\SpinCStruct_\lambda,g_0)
		\leq
		n(Y,\SpinCStruct_\lambda,g_1).
	$
	Together with the fact that
	the operator family
	$\Hodge_t\d:\Omega^1(Y)\to\Omega^1(Y)$
	has zero spectral flow
	due to $\Homology_1(Y;\R)=0$,
	the above
	implies that
	$$
		(W^{(-\mu,0)})_0
		\cong
		(W^{(-\mu,0)})_1
		\oplus
		\C^{
			n(Y,\SpinCStruct_\lambda,g_1)
			-n(Y,\SpinCStruct_\lambda,g_0)
		}.
	$$
	The result follows
	by combining this with the commuting diagram above.
}
\Rem{
	In the same vein as in
	\ManolescuSeibergWittenFloer,
	the metric invariance
	can be strengthened
	to invariance up to {\it canonical\/} isomorphism,
	which is to say,
	in this context,
	that the isomorphism does not depend
	on the path of metrics interpolating
	the given two metrics,
	but the details shall be left out.
}

\Sect{Recovery of the Cohomological Invariant}
In light of the equivalence,
proved by \LidmanManolescuEquivalence,
between the Borel $\U(1)$-equivariant cohomology
of the Seiberg-Witten Floer spectrum
and the monopole Floer ``from'' cohomology,
this section shall discuss
the relation between
the cohomotopical contact invariant
defined in the previous section
and the well known contact invariants
in Floer cohomology.
\par
In \KronheimerMrowkaOzsvathSzaboMonopoles, \S 6.3,
a distinguished element
of the monopole Floer cohomology group
$\HMFrom^*(Y,\SpinCStruct_\lambda)$,
therein denoted $\check\psi$,
is defined (up to sign)
from the datum of a contact structure
$\Ker\lambda$.
Here,
this class shall be denoted
$\psi(\lambda)$.
In fact,
most of the groundwork
for the definition of this invariant
was done a decade earlier in
\KronheimerMrowkaMonopolesAndContactStructures,
except that the machinery of
monopole Floer homology had not yet been developed.
This same invariant was studied in
\TaubesWeinsteinII,
\S 4,
therein denoted ${\frak t}_r$,
and was shown
(\TaubesWeinsteinIIAlt,
Proposition 4.3)
to be generated by
a single generator of the
monopole Floer cochain complex.
This generator is essentially
the contact configuration,
herein denoted $C_\lambda$,
with the caveat that
a generic perturbation must be used
in that context,
else the monopole Floer cohomology
groups cannot be defined.
\par
By work of
{\TaubesEmbeddedIAuthors}
(\TaubesEmbeddedIYear,
\TaubesEmbeddedIIYear,
\TaubesEmbeddedIIIYear,
{\TaubesEmbeddedIVYear}
\&
\TaubesEmbeddedVYear),
it was established that
there is a natural equivalence
between the monopole Floer cohomology
$\HMFrom^*(Y,\SpinCStruct_\lambda)$
and the embedded contact homology
$\ECH_*(Y,\lambda;0)$
of M.~Hutchings
(vid.\ \HutchingsTaubesGluingAlt).
Furthermore,
in $\ECH$,
there is a very simply defined contact invariant,
which is the class generated
by the empty set of Reeb orbits.
In \TaubesEmbeddedV,
Taubes established that,
under his isomorphism,
the $\ECH$ contact invariant
corresponds to $\psi(\lambda)$.
\par
There is yet another guise
under which $\psi(\lambda)$ appears.
In
{\KutluhanLeeTaubesHFEqHMIAuthors}
(\KutluhanLeeTaubesHFEqHMIYear,
\KutluhanLeeTaubesHFEqHMIIYear,
\KutluhanLeeTaubesHFEqHMIIIYear,
{\KutluhanLeeTaubesHFEqHMIVYear}
\&
\KutluhanLeeTaubesHFEqHMVYear),
the authors construct
isomorphisms between a variant of
the monopole Floer homologies,
called the
{\it balanced\/} monopole Floer homologies,
and the {\it Heegaard\/} Floer homologies
of \OzsvathSzaboHolomorphic;
in case $b_1=0$,
these balanced monopole Floer homologies
agree with the usual monopole Floer homologies.
In \OzsvathSzaboContact,
an invariant of contact structures
is defined (up to sign)
which lives in the group
$\HF^+(-Y,\SpinCStruct_\lambda)$;
this is usually denoted
$c^+(Y,\lambda)$.
This invariant was identified
with the $\ECH$ version
of the contact invariant
in
{\ColinGhigginiHondaEquivalenceIAuthors}
({\ColinGhigginiHondaEquivalenceIYear}
\&
\ColinGhigginiHondaEquivalenceIIYear).
\par
The manner via which
the cohomotopical contact invariant
$\Psi(\lambda)$
recovers the cohomological invariant
$\psi(\lambda)$
is fundamentally rather simple.
Returning to the motif
of the collapse of all but a single cell
in a CW-complex,
one can obtain a class in the cohomology of that complex
by pulling back the generator of the cohomology
of the sphere.
In the case of a $\U(1)$-CW-complex,
as has been said in the preceding section,
the same construction leads to a map
from the complex
to the Thom space of a vector bundle
over an orbit space.
One then proceeds by understanding
the Borel cohomology of the Thom space
as being generated by an equivariant Thom class
and the pullback of this Thom class
provides a cohomological invariant.
Since
the classical cohomological contact invariant
is the class in monopole Floer cohomology
defined by the cochain consisting
only of the generator associated to the contact monopole,
it is not at all surprising
that this approach shall work.
\par
A quick review of the pertinent definitions
is in line.
In what follows,
let $G$ denote a compact Lie group.
\Def{
	\Label\EquivariantThomClass
	Suppose a $G$-equivariant
	(reduced) cohomology theory $h^*_G$
	be given.
	By a $G$-equivariant
	Thom class,
	one means,
	for a $G$-vector bundle $V\to X$,
	a class
	$\ThomClass_G(V)\in h^*_G(\Thom_G(V))$,
	where $\Thom_G(V)$
	denotes the Thom space,
	subject to the requirement that,
	given the inclusion of any $G$-orbit
	$i:{\cal O}\to X$,
	the restriction
	$i^*\ThomClass_G(V)$
	generate
	$h^*_G(\Thom_G(i^*V))$
	as a free
	$h^*_G({\cal O}_+)$-module,
	where ${\cal O}_+$
	denotes the orbit space ${\cal O}$
	with a disjoint base point added.
}
\Rem{
	Specialising this definition
	to when $G$ act trivially
	on the base space $X$,
	the only type of orbit is
	$i:G/G\to X$,
	so the requirement is that
	$i^*\ThomClass_G(V)$
	generate
	$h^*_G(\Thom_G(V|_p))$
	as a free
	$h^*_G(S^0)$-module.
}
\Rem{
	In the other extreme,
	specialising to when $G$ act freely on $X$,
	the only type of orbit is
	$i:G/1\to X$,
	so the requirement is that
	$\ThomClass_G(V)$
	generate
	$h^*_G(\Thom_G(V))$
	as a free
	$h^*_G(X_+)$-module.
}
\Rem{
	It is not at all certain if,
	given a bundle,
	such an equivariant Thom class
	exists,
	or even if the prerequisite that 
	$h^*_G(\Thom_G(i^*V))$
	be a free
	$h^*_G((G/H)_+)$-module
	with a single generator
	is satisfied.
	In the event that it exist,
	one shall say that $V$ is
	{\it $h^*_G$-orientable.}
}
\Rem{
	\Label\ThomClassExistsForFreeBase
	Note that,
	for $G=\U(1)$
	acting on $X\cong\U(1)$ freely,
	there always exists a Thom class
	for any $\U(1)$-bundle $E\to X$
}
\Def{
	Use
	$
		\ThomClass_\lambda
		\in\BorelHomology_{\U(1)}^*(
			\ContactThom(\lambda)
		)
	$
	to denote the
	desuspension of the Thom class
	$\ThomClass_{\U(1)}(\ContactNormalUnst\mu\lambda)$
	suitably suspended or desuspended according to
	{\ContactCohomotopyInvariant},
	where
	$\ContactNormalUnst\mu\lambda\to\ContactCirc_\lambda$
	is the unstable bundle
	of $\ContactCirc_\lambda$
	in $W^\mu$.
}
\Rem{
	Of course,
	this is only defined {\it up to sign}
	which shall be consistent
	with the familiar sign ambiguity
	of the contact invariant in monopole Floer cohomology
	(cf.\ \LinRubermanSavelievFroyshovAlt).
}
\def\Perturbation{{\frak p}}
\Rem{
	In order to draw the comparison
	between the cohomological
	and cohomotopical invariants,
	one must work in the presence
	of an appropriate generic perturbation
	as that is how
	the cohomological invariant is defined
	in
	{\KronheimerMrowkaOzsvathSzaboMonopoles}
	and
	\TaubesWeinsteinII.
	It is not difficult to see that
	the construction of the invariant
	$\Psi(\lambda)$
	is not affected by the addition
	of a {\it small\/} generic perturbation.
	Indeed,
	in the case of a small generic perturbation
	$\Perturbation$,
	there still is a distinguished solution
	nearby $\ContactConf_\lambda$
	and
	unique in that it satisfies
	the lower bound in its spinor component
	as in {\SpinorContactConfBoundedbelow}.
	Call this solution
	$\ContactConf_\lambda(\Perturbation)$.
	Furthermore,
	$\ContactConf_\lambda(\Perturbation)$
	still satisfies the conclusion
	of {\NoTrajToContactConf};
	that is,
	there is no
	perturbed Seiberg-Witten trajectory
	whose forward limit be 
	$\ContactConf_\lambda(\Perturbation)$.
	This follows from the compactness results
	in \KronheimerMrowkaMonopolesAndThreeManifolds,
	Chapter 16.
	The precise details shall be left out,
	but the idea is as follows.
	If this were not true
	for generic perturbations
	of arbitrarily small norm,
	one could form a sequence
	of such perturbations
	converging to zero,
	together with respective
	perturbed Seiberg-Witten trajectories;
	however,
	the compactness properties
	of the Seiberg-Witten map
	lead to a convergent subsequence;
	the limit of this subsequence
	provides a non-trivial
	Seiberg-Witten trajectory
	into the contact configuration
	thereby contradicting
	\NoTrajToContactConf.
	For a detailed discussion
	of the generically perturbed case,
	vid.\ \RosoThesis, \S5.
}
\Thm{
	\Label\CohomotopicalRecoversCohomological
	The cohomological contact invariant
	is recovered by the cohomotopical invariant
	via
	$$
		\pm\psi(\lambda)
		=\Psi(\lambda)^*\big(
			\ThomClass_\lambda
		\big).
	$$
}
\Proof{
	The author shall begin by recalling
	the setup used in
	{\LidmanManolescuEquivalence}
	to equate the cohomology of the
	$\SWF$ spectrum
	and the monopole Floer cohomology.
	To avoid potential confusion,
	the reader is cautioned that,
	in {\LidmanManolescuEquivalence},
	the letter $\lambda$
	is used in the place where $\mu$
	is in the present \Document.
	According to \LidmanManolescuEquivalence,
	Chapter 14,
	one can find a perturbation
	$\Perturbation$
	with the properties that it be
	``very tame''
	(vid.\ \LidmanManolescuEquivalenceAlt,
	Definition 4.4.2),
	``admissible''
	(vid.\ \LidmanManolescuEquivalenceAlt,
	Definition 4.5.8
	or \KronheimerMrowkaMonopolesAndThreeManifoldsAlt,
	Definition 22.1.1)
	and that it satisfy the conclusions
	of \LidmanManolescuEquivalence,
	Proposition 7.4.1,
	Proposition 8.0.1,
	and Proposition 10.0.2.
	Denote by $\SWF_\Perturbation(Y,\SpinCStruct)$
	the version of the Seiberg-Witten Floer spectrum
	constructed with the
	finite dimensional approximations
	to the
	{\it generically perturbed\/}
	Seiberg-Witten vector field
	$$
		\GCSW_{\lambda,r;\Perturbation}
		:=\GCSW_{\lambda,r}+\GCProj_*\Perturbation.
	$$
	Then,
	according to \LidmanManolescuEquivalence,
	Proposition 6.1.6,
	the spectra
	$\SWF_\Perturbation(Y,\SpinCStruct_\lambda)$
	and $\SWF(Y,\SpinCStruct_\lambda)$
	are $\U(1)$-equivariantly
	stably homotopy equivalent.
	The point of using $\Perturbation$
	is that the cohomology of
	$\SWF_\Perturbation(Y,\SpinCStruct_\lambda)$
	may be computed via a Morse complex technique
	which allows one to equate it to
	$\HMFrom^*(Y,\SpinCStruct_\lambda)$.
	\par
	To be more precise,
	define the
	{\it anti-circular
	global Coulomb projection
	of the Seiberg-Witten vector field\/}
	as
	$$
		\AGCSW_{\lambda,r;\Perturbation}
		:=\Pi_{\AGCSlice_k}\GCSW_{\lambda,r;\Perturbation}
	$$
	where $\Pi_{\AGCSlice_j}$
	denotes the pointwise $\tilde g$-orthogonal projection
	to the local slices $\AGCSlice_{k,C}\subset W_k$
	(vid.\ {\GTildeMetric} and {\AGCSliceDef}).
	One can consider finite dimensional approximations
	of this vector field by setting
	$$
		\FDAAGCSW\mu_{\lambda,r;\Perturbation}
		:=p^\mu\AGCSW_{\lambda,r;\Perturbation}.
	$$
	Then,
	one can extend this to a vector field
	on
	$(W^\mu)^\sigma$;
	that is,
	the blow-up
	of $W^\mu$
	at the reducibles.
	Since $\U(1)$-acts freely
	on the blow-up,
	one can quotient it
	and consider
	$(W^\mu)^\sigma/\U(1)$.
	The vector field
	$\FDAAGCSW\mu_{\lambda,r;\Perturbation}$
	then uniquely defines a vector field on
	$(W^\mu)^\sigma/\U(1)$
	which shall be denoted here by
	$\QuotBlowUpFDAAGCSW\mu{\lambda,r;\Perturbation}$.
	One can form a Morse complex
	$(\hat C_\mu,\hat\partial_\mu)$,
	where the abelian group
	$\hat C_\mu$ is
	generated by the stationary points of 
	$\QuotBlowUpFDAAGCSW\mu{\lambda,r;\Perturbation}$
	inside
	$\Ball(W^\mu,2R)^\sigma/\U(1)$,
	and
	the boundary map
	$\hat\partial_\mu$
	is defined by a signed count
	of trajectories of
	$\QuotBlowUpFDAAGCSW\mu{\lambda,r;\Perturbation}$
	having index $1$.
	\par
	According to
	\LidmanManolescuEquivalence,
	Equation (272),
	the cohomology
	$\BorelHomology^i_{\U(1)}(
		\SWF_\Perturbation(Y,\SpinCStruct_\lambda))
	$
	is identified with
	the cohomology of the Morse complex
	$(\hat C_\mu,\hat\partial_\mu)$
	within a grading range
	$i\in\{-M_\mu,\ldots,M_\mu\}$,
	where $M_\mu\in\Q_{>0}$
	is a certain number which is shown
	to go to infinity as the spectral
	cut-off parameter $\mu$ goes to infinity.
	On the other hand,
	\LidmanManolescuEquivalence,
	Equation (277),
	asserts that
	the monopole Floer cohomology
	$\HMFrom^i(Y,\SpinCStruct_\lambda)$
	is also computed from
	the same complex
	$(\hat C_\mu,\hat\partial_\mu)$
	within some other grading range
	$i\in\{-N+1,\ldots,N-1\}$.
	This is done by directly identifying
	the monopoles in the grading range
	$\{-N+1,\ldots,N-1\}$
	with the fixed points of the vector field
	$\QuotBlowUpFDAAGCSW\mu{\lambda,r;\Perturbation}$
	(vid.\ \LidmanManolescuEquivalenceAlt,
	Proposition 9.3.1)
	and by identifying
	the trajectories of spectral flow $1$
	counted in the differential
	of the monopole Floer cochain complex
	with those counted by $\hat\partial_\mu$
	(vid.\ \LidmanManolescuEquivalenceAlt,
	Chapter 13).
	One then makes sure to pick $\mu$
	so that $M_\mu\geq N-1$
	and so that the grading range
	$\{-N+1,\ldots,N-1\}$
	cover at least all the {\it irreducible\/} monopoles.
	\par
	Consider now how the contact invariant
	$\psi(\lambda)$
	manifests itself in this picture.
	Firstly,
	note that
	the gauge equivalence class
	$[\ContactConf_\lambda(\Perturbation)]$
	is an irreducible monopole;
	therefore,
	it is covered by the chain level
	identification above.
	It is established in
	\TaubesWeinsteinII,
	Proposition 4.3,
	that,
	perhaps after yet another increase
	to the parameter $r$,
	the cohomological invariant
	may be described as the class
	in $\HMFrom^*(Y,\SpinCStruct_\lambda)$
	given by the cochain
	consisting only of the generator
	$[\ContactConf_\lambda(\Perturbation)]$.
	Hence,
	in the cohomology of the Morse complex
	$(\hat C_\mu,\hat\partial_\mu)$,
	the element $\psi(\lambda)$
	is represented by the gauge equivalence class
	$[\ContactConf_\lambda(\Perturbation)]$.
	On the other hand,
	the cohomology class
	$
		\Psi(\lambda)^*(
			\ThomClass_\lambda
		)
	$
	is precisely defined
	so that it be generated
	by the class
	$[\ContactConf_\lambda(\Perturbation)]$
	as well;
	this is because
	$\Psi(\lambda)$
	is the quotient map
	of the quotient of
	$\SWF_\Perturbation(Y,\SpinCStruct_\lambda)$
	by every cell other than the
	distinguished
	$\U(1)$-cell associated to
	$[\ContactConf_\lambda(\Perturbation)]$
	(vid.\ \RosoThesisAlt, \S5).
}
Given 
{\CohomotopicalRecoversCohomological},
the author intends to use
the newly constructed invariant
to deduce results about $\psi(\lambda)$
in contexts in which it has proven difficult
to do so
while relying solely on the machinery of
monopole Floer homology,
Heegaard Floer homology
and embedded contact homology.
Note that it is
{\it not\/}
clear,
at the moment,
if there is any case
in which $\Psi(\lambda)$
may hold any more information
than $\psi(\lambda)$ does.
The key advantage of $\Psi(\lambda)$
that the author wishes to emphasise
is that it does not require
a generic perturbation for its definition.
This shall be exploited in the next section.

\Sect{Finite Coverings}
In \LidmanManolescuCovering,
the authors studied
the Seiberg-Witten Floer spectrum
in the presence of a finite regular covering.
The key to their results
was the observation that
the spectrum of the manifold upstairs in the covering
acquires an action of the deck transformation group $G$
and,
upon taking appropriate fixed points of this action,
one obtains the downstairs spectrum.
A Smith-type inequality is then derived
through actual application of Smith theory.
\par
In this section,
the author's goal is
to formulate the contact invariant
in this same scenario of a finite regular covering.
This shall involve studying
the attractor-repeller pair cofibration
used to define the contact invariant
in the $G$-equivariant setting.
In doing so,
one encounters no real difficulty
and it is straightforward to obtain a
$G$-equivariant cohomotopical contact invariant.
\par
Consider a finite group $G$
and a rational homology sphere
$Y$ equipped with a free $G$-action.
Use $\pi:Y\to Y/G$
to denote the quotient map.
Agree to fix a metric $g$ on $Y/G$
and use,
on $Y$,
the induced $G$-invariant metric $\pi^*g$.
Suppose $\lambda$ be a contact form on $Y/G$
so that $\pi^*\lambda$ be a
{\it $G$-equivariant contact form\/}
on $Y$.
Notice that the canonical $\SpinC$ structure
defined by $\lambda$ on $Y/G$
naturally lifts to $Y$
as the canonical $\SpinC$ structure
defined by $\pi^*\lambda$;
that is,
$\pi^*\SpinCStruct_\lambda=\SpinCStruct_{\pi^*\lambda}$.
Likewise,
the connexion $A_{\pi^*\lambda}$
on $\det\SpinCStruct_{\pi^*\lambda}$
is the lift
$\pi^*A_\lambda$
of the connexion
$A_\lambda$ on $\det\SpinCStruct_\lambda$.
Denote by $W$ the global Coulomb slice
with respect to $A_{\pi^*\lambda}$ on $Y$
and by $W'$ the global Coulomb slice
with respect to $A_\lambda$
on $Y/G$.
\par
It is difficult
to study this scenario
in the classical
setting of monopole Floer homology
due to the need for generic perturbations
in order to achieve the required
Morse-Smale condition of Morse theory.
Indeed,
there is no guarantee that a 
sufficiently generic perturbation
chosen for $Y$
in order to satisfy the conditions
for the construction of the group
$\HMFrom^*(Y,\pi^*\SpinCStruct)$
can be made $G$-equivariant
so that it define a valid perturbation
for the construction of
$\HMFrom^*(Y/G,\SpinCStruct)$.
As a consequence,
the behaviour of the monopole Floer homology groups
under coverings
has proven elusive to study
via the classical Morse theoretic approach.
As the construction of the $\SWF$ spectrum
avoids the addition of a generic perturbation,
one can say something significant
using this machinery.
The author starts by recalling
the main observations
of \LidmanManolescuCovering.
\Rem{
	Note that $G$ acts linearly
	on the Coulomb slice $W$.
	Moreover,
	the quotient map $\pi:Y\to Y/G$
	induces an inclusion
	$W'\incl W$
	which identifies $W'$
	with the fixed point space $W^G$.
	This inclusion is
	{\it not,}
	however,
	an isometry
	in the $\LTwo$-norm.
	Nonetheless,
	the $\LTwo$ ball
	$\Ball(W',R')$
	is identified with
	the $\LTwo$ ball
	$\Ball(W, |G|R')^G\subset W$,
	which allows one to construct the Sobolev norm
	on $W$
	so as to have it be $G$-invariant.
	As a consequence,
	in the Sobolev completions,
	one still has
	$W_k'= W_k^G$.
}
\Rem{
	The Fredholm operator $\ell$
	on $Y$ is $G$-equivariant.
	Hence,
	its restriction to $W^G$
	agrees with the analogous operator
	defined on $Y/G$.
	Therefore,
	use $\ell$ to denote {\it both\/} these operators.
	Furthermore,
	note that
	$(W')^\mu=(W^\mu)^G$.
	The map $c$ is also $G$-equivariant,
	so a finite-type Seiberg-Witten trajectory
	in $W_k'$
	is the same thing as
	a $G$-fixed finite-type Seiberg-Witten trajectory
	in $W_k$.
	The same identification can be made
	between the Seiberg-Witten trajectories
	in the {\it finite dimensional approximations.}
}
Use $R>0$
to denote the constant provided by
{\FiniteTypeBounded}
for the case of the manifold $Y$
and $R'>0$
to denote this constant
for the quotient manifold $Y/G$.
By perhaps increasing $R$ or $R'$,
one can ensure that
$R'=R/|G|$.
This means that $R$ will also satisfy
the conclusions of
{\FiniteTypeBounded}
for $Y/G$.
\par
Next,
recall from {\FixBumpFunction}
that the choice of bump functions
$u^\mu$
for $Y$
was made so as to have it constant
on the spheres centred at zero
in the Sobolev norm;
therefore,
the $u^\mu$
are automatically $G$-invariant.
Therefore,
their restrictions to $(W^\mu)^G$
can be used to define the bump functions
required for $Y/G$.
With these conditions,
the
finite dimensional Seiberg-Witten flow
$\FDASWFlow^\mu_{\pi^*\lambda,r}$
of $Y$
is $G$-equivariant
and restricts to $(W^\mu)^G$
as the finite dimensional Seiberg-Witten flow
of $Y/G$.
Now,
if the reader agree to fix 
$\mu>0$
large enough
so as to have
{\FDCompactness}
hold
for both $Y$ and $Y/G$,
then,
it follows that
the isolated invariant set
$S^\mu_{\pi^*\lambda,r}$
is $G$-invariant
and its $G$-fixed subset,
$(S^\mu_{\pi^*\lambda,r})^G$,
is precisely the isolated invariant set
used to define the $\SWF$ spectrum
of $Y/G$.
Moreover,
since the $G$-action
is linear on $W^\mu$,
the isolating neighbourhood
$\Disk(W^\mu,2R)$
of 
$S^\mu_{\pi^*\lambda,r}$
is $G$-invariant
and
$\Disk(W^\mu,R)^G=\Disk\big((W^\mu)^G,R\big)$
serves as an isolating neighbourhood
for the construction of the $\SWF$ spectrum
on $Y/G$.
\par
A few more notions from equivariant stable homotopy theory
shall be required.
In what follows,
use $H$ to denote an arbitrary compact Lie group.
As before,
the reader is directed to {\MayEtAlEquivariant}
for further details.
\par
\Thm{
	Given a closed normal subgroup $K\subset H$,
	there exists a functor,
	called the
	{\it geometric fixed points functor,}
	from $H$-spectra indexed over a universe $\AnotherUniverse$
	to $(H/K)$-spectra indexed over the universe $\AnotherUniverse^K$,
	$$
		\Phi^K:
		\Spectra H\AnotherUniverse
		\to
		\Spectra {(H/K)}\AnotherUniverse^K,
	$$
	satisfying the properties that
	{
		\parindent = 1cm
		\item{(i)}
			$\Phi^K\Sigma^\infty_\AnotherUniverse X
			\cong\Sigma^\infty_{\AnotherUniverse^K} X^K$,
		\item{(ii)}
			$\Phi^K E\wedge \Phi^K F\cong\Phi^K(E\wedge F)$.
	}
}
\Proof{
	Vid.\ \MayEtAlEquivariant, \S XVI.3.
}
\Rem{
	Since the author is simply dealing with suspension spectra,
	these two properties suffice
	in understanding the geometric fixed points
	of the spectra at hand.
	In particular,
	note that,
	for an $H$-space $X$
	and an $H$-representation $V$,
	it follows that
	$$
		\Phi^K\Sigma^{-V}\Sigma^\infty_\AnotherUniverse X
		=\Sigma^{-V^K}\Sigma^\infty_{\AnotherUniverse^K} X^K.
	$$
}
\Def{
	An $H$-universe is called
	{\it complete\/}
	if one can find,
	for any finite dimensional $H$-representation $V$,
	a sub-representation in $\AnotherUniverse$
	isomorphic to $V$.
}
\Rem{
	For $G$ the group of deck transformations
	of the covering $\pi:Y\to Y/G$,
	notice that the universe
	$\CoulombUniverse$
	defined by the Coulomb gauge of $Y$
	is naturally a $G\times\U(1)$-universe
	and its $G$-fixed point space,
	$\CoulombUniverse^G$,
	is the $\U(1)$-universe
	defined by the Coulomb gauge of $Y/G$.
	Let $\Universe'$ denote a
	{\it complete\/}
	$\U(1)\times G$-universe.
	By intertwining with change of universe functors
	defined by an isometry
	$\CoulombUniverse\to\Universe'$,
	as was done in
	\DesuspensionIntertwinedWithChangesOfUniverse,
	one can consider the functor
	$\Sigma^{-W^{(-\mu,0)}}$
	as an endofunctor
	of the category
	$\StabHoSpectra{(\U(1)\times G)}\Universe'$.
}
\Def{
	Define the
	{\it metric dependent
	$G$-equivariant
	Seiberg-Witten Floer spectrum\/}
	as
	$$
		\SWF_G(Y,\pi^*\SpinCStruct_\lambda,\pi^*g)
		:=
		\Sigma^{-W^{(-\mu,0)}}
		\Sigma^\infty_{\Universe'}
		\Conley_{\U(1)\times G}(
			S^\mu_{\pi^*\lambda,r},
			\FDASWFlow_{\pi^*\lambda,r}^\mu
		)
		\in\StabHoSpectra{(\U(1)\times G)}{\Universe'}.
	$$
}
\Rem{
	The author believes it to be possible to
	(de)suspend away the metric dependence
	in an analogous fashion to the non-equivariant case;
	however,
	it seems the details are somewhat subtle,
	so he chose not to pursue that goal in this \Document.
	This would involve considering a localization
	of the representation ring $\RO(\U(1)\times G)$,
	equivariant spectral flow
	and equivariant eta invariants.
}
\Thm{
	\Label\SpectrumFixedPoints
	(\LidmanManolescuCoveringAlt)
	The Seiberg-Witten Floer spectra
	for $Y$ and $Y/G$ are related by
	$$
		\Phi^G\SWF_G(Y,\pi^*\SpinCStruct,\pi^*g)
		=\SWF(Y/G,\SpinCStruct,g).
	$$
}
\Proof{
	Although the statement in {\LidmanManolescuCovering}
	is made in terms of $\U(1)\times G$-spaces
	instead of spectra,
	the properties of the geometric fixed points functor
	mentioned above make clear that
	this is the correct statement for spectra.
}
Now,
the contact invariant construction
shall be considered in this setting.
If necessary,
increase $R$,
$\mu$
and the parameter $r$
used in the Seiberg-Witten equations
so that 
{\ContactConfInCoulombBall},
{\ContactCircNonDegen}
and {\ContactCircRepeller}
be satisfied for
{\it both\/}
$Y$ and $Y/G$.
\par
Notice that
the upstairs contact circle,
$
	\ContactCirc_{\pi^*\lambda}
	\subset(W^\mu)^G
	\subset W^\mu
$,
is naturally identified
with the downstairs
contact circle,
$\ContactCirc_\lambda$.
Beware,
however,
that the dual attractors
to $\ContactCirc_{\pi^*\lambda}$
in $S^\mu_{\pi^*\lambda,r}$
and
in $(S^\mu_{\pi^*\lambda,r})^G$
are,
of course,
{\it not the same.}
\Def{
	\Label\EquivariantContactThomIntro
	Let
	$$
		\ContactThom_G(\pi^*\lambda,g)
		:=\Sigma^{W^{(-\mu,0)}}\Sigma^\infty_{\Universe'}
		\Conley_{\U(1)\times G}(
			\ContactCirc_{\pi^*\lambda},
			\FDASWFlow_{\pi^*\lambda,r}^\mu
		)
		\in\StabHoSpectra{(\U(1)\times G)}{\Universe'}
	$$
	denote the Thom space
	appearing in the codomain of the
	cohomotopical contact invariant
	but now seen as a $(\U(1)\times G)$-spectrum
	(cf.\ \MetricDepContactThomIntro).
}
The attractor-repeller cofibration
$$
	\Conley_{\U(1)\times G}\big(
		(\ContactCirc_{\pi^*\lambda})^*,
		\FDASWFlow_{\pi^*\lambda,r}^\mu
	\big)
	\to
	\Conley_{\U(1)\times G}\big(
		S^\mu_{\pi^*\lambda,r},
		\FDASWFlow_{\pi^*\lambda,r}^\mu
	\big)
	\to
	\Conley_{\U(1)\times G}\big(
		\ContactCirc_{\pi^*\lambda},
		\FDASWFlow_{\pi^*\lambda,r}^\mu
	\big)
$$
can be formed $G$-equivariantly as well.
This is done 
by constructing
a $(\U(1)\times G)$-invariant
index triple $(L,M,N)$
according to
{\IndexTripleExistsCofibWellDef}.
The consequence
is that the contact invariant
gains a $G$-equivariant version.
\Def{
	Define the
	{\it $G$-equivariant
	metric dependent
	contact invariant\/}
	as the map
	$$
		\Psi_G(\pi^*\lambda,\pi^*g):
		\SWF_G(Y,\pi^*\lambda,\pi^*g)
		\to\ContactThom_G(\pi^*\lambda,\pi^*g)
	$$
	given by desuspending the above cofibre map
	by the $(\U(1)\times G)$-representation $W^{(-\mu,0)}$.
}
\Thm{
	The cohomotopical contact invariants
	of $(Y,\pi^*\lambda)$ and $(Y/G,\lambda)$
	are related by
	$$
		\Phi^G\Psi_G(\pi^*\lambda,\pi^*g)
		=\Psi(\lambda,g).
	$$
}
\Proof{
	If $(L,M,N)$
	be a $\U(1)\times G$-invariant index triple
	for $S^\mu_{\pi^*\lambda,r}$,
	then the triple of fixed point sets
	$(L^G,M^G,\allowbreak N^G)$
	is a $\U(1)$-invariant index triple for
	$(S^\mu_{\pi^*\lambda,r})^G=S^\mu_{\lambda,r}$.
	The result follows at once.
}
The extra data available due to the $G$-action
allows one to define
a $G$-equi\-vari\-ant 
(or $(\U(1)\times G)$-equivariant)
version
of the monopole Floer cohomologies
via the $G$-equivariant cohomology of the
$\SWF$ spectrum
since the spectrum has naturally acquired
a $G$-action.
There are three choices of equivariant cohomology theory
which spring to mind in this context:
{\it Borel cohomology,}
{\it equivariant K-theory\/}
and
{\it Bredon cohomology.}
In the present \Document,
the author shall pursue the use of Borel cohomology
as that turned out to be the most readily applicable.

\Sect{Borel Cohomology}
This section shall deal
with the Borel $G$-equivariant cohomology
of the spectrum $\SWF_G$
and consider the resulting equivariant contact class.
Despite the spectrum $\SWF_G$
being a $\U(1)\times G$-spectrum,
in the interests of simplicity,
the author shall ignore the $\U(1)$-action
and only consider cohomologies equivariant
with respect to the $G$-action.
Application of Borel cohomology
to the cohomotopical contact invariant
leads to a
{\it $G$-equivariant cohomological contact invariant\/}
$\psi_G$,
which,
as shall be seen,
holds information
about both the
{\it upstairs\/}
and
{\it downstairs\/}
contact structures of the covering.
Verily,
if one know the value of $\psi_G$,
one can recover the value of both
cohomological contact invariants
of the covering.
\par
Ideally,
what one really wishes
is to infer something about the
{\it upstairs\/}
contact invariant
from knowledge about the
{\it downstairs\/}
contact invariant;
the new
{\it equivariant\/}
contact invariant,
therefore,
may seem not to be leading in that direction.
However,
Borel cohomology,
under appropriate circumstances,
enjoys the very powerful property
that
the cohomology of the fixed points of a $G$-space
conditions significantly
the Borel cohomology of the whole space;
this is the content of the 
{\it localization theorem.}
As a consequence of localization,
one
{\it sometimes\/}
can determine
the equivariant contact invariant
from knowledge of the downstairs contact invariant.
And,
as the equivariant contact invariant
recovers the non-equivariant upstairs contact invariant,
this allows one to infer the upstairs contact invariant
starting only from knowledge of the downstairs contact invariant.
\par
The applicability of localization,
however,
is limited to scenarios where one know precisely
what the $\SWF_G$ spectrum is.
This happens,
for instance,
in the event that there be a unique solution to the Seiberg-Witten equations,
which implies that $\SWF_G$ is an equivariant sphere spectrum.
Nonetheless,
such manifolds are sufficiently abundant
that the results derived
say something quite non-trivial.
A special case of interest
shall be that of {\it elliptic manifolds,}
where a very general theorem may be stated
which aids greatly
in determining whether
the lift of a tight contact structure remains tight
or becomes overtwisted.
\par
For further simplicity,
the author shall avoid the language of $G$-spectra
and work instead at the level of $G$-spaces
as much as possible in this section.
This is similar to what is done in
\LidmanManolescuCovering.
The main reason being that the appropriate form
of the localization theorem
is more difficult to describe for spectra
as questions about the various forms 
of fixed points functors come into play.
\par
Let $Y$, $G$, $\pi$, $g$ and $\lambda$ be as in the previous section.
\Def{
	Define the
	{\it $G$-equivariant Borel metric dependent
	monopole Floer cohomology\/}
	as
	$$
		\BorelHMTilde^*_{G}
		(Y,\pi^*\SpinCStruct,\pi^*g;R)
		:=\BorelHomology^*_G(
			\SWF_G(Y,\pi^*\SpinCStruct,\pi^*g);
			R
		)
	$$
	where $R$ is some commutative ring,
	which will be left implicit in the notation
	henceforth.
}
\Rem{
	In order to make the best out of Borel cohomology,
	it is best to focus on the case $G=\Z/p\Z$ 
	for a prime $p$.
	In this case,
	the coefficient ring $R$
	shall be chosen to be $\Z/p\Z$.
	For the remainder of this section,
	these choices shall be implicit
	lest the notation become overloaded.
}
\Rem{
	Since $\SWF_G(Y,\pi^*\SpinCStruct,\pi^*g)$
	depends on $g$ only up to suspension
	by $G$-represen\-tations,
	it follows that
	$\BorelHMTilde^*_G(Y,\pi^*\SpinCStruct,\pi^*g)$
	depends on $g$ only up to shifts in grading,
	which,
	although cosmetically unpleasant,
	is not a major issue in practice.
}
\Def{
	Given a $G$-representation $V$
	and a ring $R$,
	one says that $V$ is
	{\it $G$-equivariantly $R$-orientable\/}
	if the vector bundle
	$V_G:=(V\times \E G)/G$
	over $\B G$
	be $R$-orientable.
	In which case,
	one writes $e_G(V)\in\BorelHomology^*_G(S^0;R)$
	for its {\it Euler class.}
}
\Rem{
	Since
	the author is using coefficients $R=\Z/p\Z$
	for the group $G=\Z/p\Z$
	and $p$ is a prime,
	it follows that all $G$-representations
	are $G$-equivariantly $R$-orientable.
	In the case $p=2$,
	all vector bundles are $R$-orientable anyway.
	Otherwise,
	if $p$ be an odd prime,
	then any non-trivial representation
	of $G$ is complex and,
	therefore,
	the vector bundle defined over $\B G$
	shall also be complex
	and therefore orientable over any field.
	So,
	in any event,
	an equivariant Euler class always exists
	for the purposes being pursued here.
}
\Rem{
	Note that,
	due to the Thom isomorphism,
	for an orientable $G$-representation $V$,
	the ring
	$\BorelHomology^*_G(V^+)$
	is isomorphic to
	$\BorelHomology^*_G(S^0)$
	with its grading shifted by
	$\dim V$.
	Moreover,
	under this isomorphism,
	notice that the Thom class
	$\ThomClass_G(V_G)$
	of the bundle $V_G\to\B G$
	gets sent to 
	$1\in\RedHomology^0_G(S^0)$.
	In other words,
	one can think of
	$\ThomClass_G(V_G)$
	as a $G$-equivariant fundamental class
	of the $G$-manifold $V^+$.
}
\Thm{
	(Localization Theorem)
	Let $\Gamma$ be a finite $G$-CW-complex
	and $S\subset\RedHomology^*(\B G;\Z/p\Z)$
	consist of those elements
	which be Euler classes of $G$-representations
	having no trivial summand.
	The inclusion of fixed points
	$\Gamma^G\incl \Gamma$
	induces an
	{\it isomorphism between cohomologies
	localized with respect to $S$,}
	$$
		S^{-1}\BorelHomology^*_G(\Gamma;\Z/p\Z)
		\isomorph
		S^{-1}\BorelHomology^*_G(\Gamma^G;\Z/p\Z).
	$$
}
\Proof{
	Vid.\ \TomDieckTransformationAlt,
	Theorem 3.13.
}
\Eg{
	\Label\LocalizationForRepresentationSphereModOddPrime
	Consider the inclusion of fixed points
	$$
		(V^+)^G\incl V^+
	$$
	of a $G$-representation sphere
	coming from a $G$-representation $V$
	potentially having $V^G\neq 0$.
	In this example,
	assume $p>2$.
	Firstly,
	recall that,
	in this case,
	the cohomology of $\B G$
	is the ring
	$$
		\BorelHomology_G^*(S^0)
		=(\Z/p\Z)[u,v]/(v^2)
	$$
	where $\deg u=2$,
	$\deg v=1$.
	It is worth paying close attention to this ring.
	In the following diagram,
	the top row indicates the abelian subgroups at each degree,
	the middle row indicates the generator with the corresponding degree
	and the bottom row indicates the numerical value of the degree.
	$$\matrix{
		\cdots &
		0 &
		0 &
		\Z/p\Z &
		\Z/p\Z &
		\Z/p\Z &
		\Z/p\Z &
		\Z/p\Z &
		\cdots
		\cr
		&
		&
		&
		\langle 1\rangle &
		\langle v\rangle &
		\langle u\rangle &
		\langle uv\rangle &
		\langle u^2\rangle &
		\cr
		&
		(-2) &
		(-1) &
		(0) &
		(1) &
		(2) &
		(3) &
		(4).
		\cr
	}$$
	Now,
	because $p>2$ is an odd prime,
	all representations of $G$ are complex
	and therefore define orientable bundles over $\B G$.
	It follows from the Thom isomorphism that
	$$
		\BorelHomology^*_G(V^+)
		\cong
		\BorelHomology^{*-\dim V}_G(S^0)
	$$
	and likewise for $(V^+)^G$.
	Hence,
	{\it before localization,}
	the inclusion of fixed points
	$(V^+)^G\incl V^+$
	induces a map
	$$
		\BorelHomology^{*-\dim V}_G(S^0)
		\to
		\BorelHomology^{*-\dim V^G}_G(S^0)
	$$
	which,
	in the notation above,
	is depicted as
	$$\matrix{
		\cdots &
		0 &
		0 &
		\cdots &
		0 &
		\Z/p\Z &
		\cdots
		\cr
		&
		\downarrow &
		\downarrow &
		&
		\downarrow &
		\downarrow &
		\cr
		\cdots &
		0 &
		\Z/p\Z &
		\cdots &
		\Z/p\Z &
		\Z/p\Z &
		\cdots,
	}$$
	where all the $\Z/p\Z\to\Z/p\Z$ maps are isomorphisms.
	Now,
	consider the effect of localization.
	The set $S$ with respect to which
	one must perform localization
	is the set $\{u^n\mid n>0\}$;
	this can be seen from looking at
	the representation theory of $G$.
	The {\it localized\/}
	Borel cohomology
	therefore is
	$$
		S^{-1}\BorelHomology^*_G(S^0)
		=(\Z/p\Z)[u,u^{-1},v]/(v^2).
	$$
	Schematically,
	$$\matrix{
		\cdots &
		\Z/p\Z &
		\Z/p\Z &
		\Z/p\Z &
		\Z/p\Z &
		\Z/p\Z &
		\Z/p\Z &
		\Z/p\Z &
		\cdots
		\cr
		&
		\langle u^{-1}\rangle&
		\langle u^{-1}v\rangle &
		\langle 1\rangle &
		\langle v\rangle &
		\langle u\rangle &
		\langle uv\rangle &
		\langle u^2\rangle &
		\cr
		&
		(-2) &
		(-1) &
		(0) &
		(1) &
		(2) &
		(3) &
		(4).
		\cr
	}$$
	After localization,
	the map induced by the inclusion of fixed points
	takes the form
	$$\matrix{
		\cdots &
		\Z/p\Z &
		\Z/p\Z &
		\cdots &
		\Z/p\Z &
		\Z/p\Z &
		\cdots
		\cr
		&
		\downarrow &
		\downarrow &
		&
		\downarrow &
		\downarrow &
		\cr
		\cdots &
		\Z/p\Z &
		\Z/p\Z &
		\cdots &
		\Z/p\Z &
		\Z/p\Z &
		\cdots,
	}$$
	where {\it all\/} maps are isomorphisms.
}
\Eg{
	\Label\LocalizationForRepresentationSphereModTwo
	Now,
	consider the same scenario
	but with $p=2$.
	In this case,
	recall that
	the cohomology of $\B G$
	is the ring
	$$
		\BorelHomology_G^*(S^0)
		=(\Z/2\Z)[u]
	$$
	where $\deg u=1$.
	In particular,
	as an abelian group,
	this is the same as in the case $p>2$;
	only the ring structures differ.
	Again,
	draw the same sort of diagram as before.
	$$\matrix{
		\cdots &
		0 &
		0 &
		\Z/2\Z &
		\Z/2\Z &
		\Z/2\Z &
		\Z/2\Z &
		\Z/2\Z &
		\cdots
		\cr
		&
		&
		&
		\langle 1\rangle &
		\langle u\rangle &
		\langle u^2\rangle &
		\langle u^3\rangle &
		\langle u^4\rangle &
		\cr
		&
		(-2) &
		(-1) &
		(0) &
		(1) &
		(2) &
		(3) &
		(4).
		\cr
	}$$
	Unlike in the $p>2$ case,
	here,
	there are non-trivial real representations of $G$
	which define non-orientable bundles over $\B G$.
	No matter;
	in this case,
	the coefficient ring being used for cohomology theories
	is $\Z/2\Z$,
	and all vector bundles are $\Z/2\Z$-orientable.
	Hence,
	again by the Thom isomorphism,
	$$
		\BorelHomology^*_G(V^+)
		\cong
		\BorelHomology^{*-\dim V}_G(S^0)
	$$
	and likewise for $(V^+)^G$.
	As before,
	the inclusion of fixed points
	$(V^+)^G\incl V^+$
	induces a map
	$$
		\BorelHomology^{*-\dim V}_G(S^0)
		\to
		\BorelHomology^{*-\dim V^G}_G(S^0)
	$$
	which
	one depicts as
	$$\matrix{
		\cdots &
		0 &
		0 &
		\cdots &
		0 &
		\Z/2\Z &
		\cdots
		\cr
		&
		\downarrow &
		\downarrow &
		&
		\downarrow &
		\downarrow &
		\cr
		\cdots &
		0 &
		\Z/2\Z &
		\cdots &
		\Z/2\Z &
		\Z/2\Z &
		\cdots,
	}$$
	where all the $\Z/2\Z\to\Z/2\Z$ maps are $1\mapsto 1$.
	Localization,
	in this case,
	is with respect to the set
	$S=\{u^n\mid n>0\}$.
	Therefore,
	$$
		S^{-1}\BorelHomology^*_G(S^0)
		=(\Z/2\Z)[u,u^{-1}].
	$$
	Schematically,
	$$\matrix{
		\cdots &
		\Z/2\Z &
		\Z/2\Z &
		\Z/2\Z &
		\Z/2\Z &
		\Z/2\Z &
		\Z/2\Z &
		\Z/2\Z &
		\cdots
		\cr
		&
		\langle u^{-2}\rangle&
		\langle u^{-1}\rangle &
		\langle 1\rangle &
		\langle u\rangle &
		\langle u^2\rangle &
		\langle u^3\rangle &
		\langle u^4\rangle &
		\cr
		&
		(-2) &
		(-1) &
		(0) &
		(1) &
		(2) &
		(3) &
		(4).
		\cr
	}$$
	The end result
	is effectively the same as in the $p>2$ case;
	the map induced by the inclusion of fixed points
	is
	$$\matrix{
		\cdots &
		\Z/2\Z &
		\Z/2\Z &
		\cdots &
		\Z/2\Z &
		\Z/2\Z &
		\cdots
		\cr
		&
		\downarrow &
		\downarrow &
		&
		\downarrow &
		\downarrow &
		\cr
		\cdots &
		\Z/2\Z &
		\Z/2\Z &
		\cdots &
		\Z/2\Z &
		\Z/2\Z &
		\cdots,
	}$$
	where {\it all\/} maps are isomorphisms.
}
\Rem{
	\Label\BorelCohomologyThomSpaces
	Recall the {\it unstable normal bundle\/}
	$
		\ContactNormalUnst\mu{\pi^*\lambda}
		\to\ContactCirc_{\pi^*\lambda}
	$
	from {\IntroUnstableBundle}.
	As was seen earlier,
	the $G$-equivariant Conley index of
	$\ContactCirc_{\pi^*\lambda}$
	is the Thom space
	of this bundle,
	$$
		\Conley_G(
			\ContactCirc_{\pi^*\lambda},
			\FDASWFlow^\mu_{\pi^*\lambda,r}
		)
		=\Thom_G(\ContactNormalUnst\mu{\pi^*\lambda}).
	$$
	Let
	$\ContactNVBF\incl\ContactNormal\mu\lambda$
	denote a fibre of the bundle.
	Note that
	$\ContactNVBF$ is a $G$-representation
	and decompose it as
	$\ContactNVBF^G\oplus\ContactFNVBF$
	where $\ContactFNVBF$ is a $G$-representation
	with trivial fixed points,
	$\ContactFNVBF^G=0$.
	Next,
	consider the unstable normal bundle
	$\ContactNormalUnst\mu\lambda\to\ContactCirc_\lambda$.
	The Conley index is again computed as
	$$
		\Conley(
			\ContactCirc_\lambda,
			\FDASWFlow^\mu_{\lambda,r}
		)
		=\Thom(\ContactNormalUnst\mu\lambda).
	$$
	One can ensure that
	the unstable normal bundles be related as
	$$
		\ContactNormalUnst\mu\lambda
		=(\ContactNormalUnst\mu{\pi^*\lambda})^G.
	$$
	Hence,
	the Conley indices
	of $\ContactCirc_\lambda$ in $(W^\mu)^G$
	and of $\ContactCirc_{\pi^*\lambda}$ in $W^\mu$
	are related by
	$$
		\Conley_G(
			\ContactCirc_{\pi^*\lambda},
			\FDASWFlow^\mu_{\pi^*\lambda,r}
		)
		=\Conley(
			\ContactCirc_{\lambda},
			\FDASWFlow^\mu_{\lambda,r}
		)\wedge\ContactFNVBF^+.
	$$
	Since
	$
		\Conley(
			\ContactCirc_{\lambda},
			\FDASWFlow^\mu_{\lambda,r}
		)
	$
	is $G$-trivial,
	one has that
	$$
		\BorelHomology_G^*(
			\Conley(
				\ContactCirc_{\lambda},
				\FDASWFlow^\mu_{\lambda,r}
			)
		)
		\cong
		\RedHomology^*(
			\Conley(
				\ContactCirc_{\lambda},
				\FDASWFlow^\mu_{\lambda,r}
			)
		)
		\mathop\otimes_{\Z/p\Z}
		\BorelHomology_G^*(S^0)
	$$
	and,
	similarly,
	$$
		\BorelHomology_G^*(
			\Conley(
				\ContactCirc_{\lambda},
				\FDASWFlow^\mu_{\lambda,r}
			)
			\wedge
			\ContactFNVBF^+
		)
		\cong
		\RedHomology^*(
			\Conley(
				\ContactCirc_{\lambda},
				\FDASWFlow^\mu_{\lambda,r}
			)
		)
		\mathop\otimes_{\Z/p\Z}
		\BorelHomology_G^*(\ContactFNVBF^+).
	$$
}
\Rem{
	Next,
	consider a {\it non-equivariant\/} Thom class
	$$
		\ThomClass(\ContactNormal\mu\lambda)
		\in
		\RedHomology^{\dim\ContactNVBF^G}(
			\Conley(
				\ContactCirc_{\lambda},
				\FDASWFlow^\mu_{\lambda,r}
			)
		)
	$$
	for the bundle
	$
		\ContactNormalUnst\mu\lambda
		\to\ContactCirc_{\lambda}
	$.
	Bearing in mind the isomorphisms outlined
	in the previous remark,
	define a class
	$$
		\ThomClass_G(\ContactNormalUnst\mu{\pi^*\lambda})
		:=\ThomClass(\ContactNormalUnst\mu\lambda)
		\otimes \ThomClass_G(\ContactFNVBF_G)
		\in
		\RedHomology^*(
			\Thom(\ContactNormalUnst\mu{\lambda})
		)
		\otimes_{\Z/p\Z}
		\BorelHomology_G^*(\ContactFNVBF^+)
		\cong
		\BorelHomology_G^*(
			\Thom_G(\ContactNormalUnst\mu{\pi^*\lambda})
		).
	$$
	It is easy to verify that the class
	$\ThomClass_G(\ContactNormalUnst\mu{\pi^*\lambda})$
	serves as a $G$-equivariant Thom class
	for the $G$-bundle
	$
		\ContactNormalUnst\mu{\pi^*\lambda}
		\to\ContactCirc_{\pi^*\lambda}
	$.
}
\Rem{
	Due to the difficulties of working
	with localization in the context of spectra,
	the author decided to proceed with certain arguments
	applied prior to desuspension.
	For that end,
	it becomes useful to define
	contact invariants dependent
	on the sufficiently large spectral cut-off parameter.
	In what follows,
	let $\mu>0$
	again be large enough to satisfy what is said in
	{\FDCompactness}
	and \FDAContactCircNonDegen.
}
\Def{
	Let the
	{\it cohomological contact invariant
	in finite dimensional approximation\/}
	be the class
	$$
		\psi(\lambda,g,\mu)\in\RedHomology^*(
			\Conley(S^\mu_{\lambda,r},\FDASWFlow^\mu_{\lambda,r})
		)
	$$
	given as the pullback of
	$
		\ThomClass(\ContactNormalUnst\mu\lambda)
	$
	via the cofibre map
	$$
		\Conley(S^\mu_{\lambda,r},\FDASWFlow^\mu_{\lambda,r})
		\to\Thom(\ContactNormalUnst\mu\lambda).
	$$
}
\Def{
	Let the
	{\it equivariant cohomological contact invariant
	in finite dimensional approximation\/}
	be the class
	$$
		\psi_G(\pi^*\lambda,\pi^* g,\mu)
		\in\BorelHomology^*_G(
			\Conley_G(
				S^\mu_{\pi^*\lambda,r},
				\FDASWFlow^\mu_{\pi^*\lambda,r}
			)
		)
	$$
	given as the pullback of
	$
		\ThomClass_G(\ContactNormalUnst\mu{\pi^*\lambda})
	$
	via the cofibre map
	$$
		\Conley_G(
			S^\mu_{\pi^*\lambda,r},
			\FDASWFlow^\mu_{\pi^*\lambda,r}
		)
		\to\Thom_G(\ContactNormalUnst\mu{\pi^*\lambda}).
	$$
}
\Prop{
	\Label\NonDesuspendedEquivariantMapsToDown
	Under the map induced by inclusion of fixed points
	$$
		\BorelHomology_G^*(
			\Conley_G(
				S^\mu_{\pi^*\lambda,r},
				\FDASWFlow^\mu_{\pi^*\lambda,r}
			)
		)
		\to
		\BorelHomology_G^*(
			\Conley(
				S^\mu_{\lambda,r},
				\FDASWFlow^\mu_{\lambda,r}
			)
		)
		\cong
		\RedHomology^*(
			\Conley(
				S^\mu_{\lambda,r},
				\FDASWFlow^\mu_{\lambda,r}
			)
		)
		\otimes
		\BorelHomology_G^*(S^0)
	$$
	the class $\psi_G(\pi^*\lambda,\pi^*g,\mu)$
	gets sent to
	$\psi(\lambda,g,\mu)\otimes e_G(\ContactFNVBF)$.
}
\Proof{
	Consider the commuting diagram
	$$
		\Diagram{
			\BorelHomology_G^*(
				\Conley_G(
					S^\mu_{\pi^*\lambda,r},
					\FDASWFlow^\mu_{
						\pi^*\lambda,
						r
					}
				)
			)
			&
			\Left
			&
			\BorelHomology_G^*(
				\Thom_G(\ContactNormalUnst\mu{\pi^*\lambda})
			)
			\cr
			\Down
			&
			&
			\Down
			\cr
			\BorelHomology_G^*(
				\Conley(
					S^\mu_{\lambda,r},
					\FDASWFlow^\mu_{
						\lambda,
						r
					}
				)
			)
			&
			\Left
			&
			\BorelHomology_G^*(
				\Thom(\ContactNormalUnst\mu\lambda)
			),
			\cr
		}
	$$
	where the horizontal arrows come from the cofibre maps
	and the vertical from the inclusion of fixed points.
	By
	the isomorphisms discussed in
	\BorelCohomologyThomSpaces,
	one can rewrite the right vertical map as 
	$$
		\RedHomology^*(\Thom(\ContactNormalUnst\mu\lambda))
		\otimes
		\BorelHomology_G^*(\ContactFNVBF^+)
		\to
		\RedHomology^*(\Thom(\ContactNormalUnst\mu\lambda))
		\otimes
		\BorelHomology_G^*(S^0).
	$$
	Moreover,
	with respect to these tensor products,
	this map is simply
	$\id\otimes i^*$
	where $i$
	is the inclusion of fixed points
	$i:S^0\incl\ContactFNVBF^+$.
	But note that
	$i^*\ThomClass_G(\ContactFNVBF_G)$
	is,
	by definition,
	the Euler class
	$e_G(\ContactFNVBF)$.
	The result then follows by commutativity
	of the diagram.
}
\Prop{
	\Label\NonDesuspendedEquivariantDeterminesDown
	If $\psi_G(\pi^*\lambda,\pi^*g,\mu)=0$,
	then $\psi(\lambda,g,\mu)=0$.
}
\Proof{
	Since the representation $\ContactFNVBF$
	has no trivial summand,
	the Euler class $e_G(\ContactFNVBF)$ is never zero.
	The result then follows immediately from
	\NonDesuspendedEquivariantMapsToDown.
}
One can now desuspend appropriately
to obtain a result that is not dependent
on the spectral cut-off $\mu$.
\Def{
	Define the
	{\it metric dependent equivariant
	cohomological contact invariant}
	$$
		\psi_G(\pi^*\lambda,\pi^*g)
		\in\BorelHMTilde_G^*(Y,\pi^*\SpinCStruct_\lambda,\pi^*g)
	$$
	as the desuspension of the class
	$\psi_G(\pi^*\lambda,\pi^*g,\mu)$
	by the $G$-representation $W^{(-\mu,0)}$.
}
\Thm{
	\Label\EquivariantDeterminesDown
	If $\psi_G(\pi^*\lambda,\pi^*g)=0$,
	then $\psi(\lambda,g)=0$.
}
\Proof{
	Follows directly from
	{\NonDesuspendedEquivariantDeterminesDown}
	after desuspension.
}
\Rem{
	Recall that there is a natural transformation
	$\BorelHomology_G^*\to\RedHomology^*$
	induced by the inclusion
	of the fibre of the Borel construction.
	That is,
	if $\Gamma$ be a $G$-space,
	$\{p\}\times \Gamma\incl \Gamma_G:=(EG\times \Gamma)/G$.
	In the present context,
	this translates to a forgetful map
	$$
		\BorelHMTilde_G^*(Y,\pi^*\SpinCStruct_\lambda,\pi^*g)
		\to
		\HMTilde^*(Y,\pi^*\SpinCStruct_\lambda,\pi^*g).
	$$
}
\Thm{
	\Label\EquivariantMapsToUp
	Under the map
	$$
		\BorelHMTilde_G^*(Y,\pi^*\SpinCStruct_\lambda,\pi^*g)
		\to
		\HMTilde^*(Y,\pi^*\SpinCStruct_\lambda,\pi^*g),
	$$
	the class
	$\psi_G(\pi^*\lambda,\pi^*g)$
	gets sent to
	$\psi(\pi^*\lambda,\pi^*g)$.
}
\Proof{
	This follows simply from the fact
	that the equivariant Thom class
	$$
		\ThomClass_G(\ContactNormalUnst\mu{\pi^*\lambda})
		\in\BorelHomology^*_G(
			\Conley_G(\ContactCirc_{\pi^*\lambda})
		)
	$$
	is sent to a {\it non-equivariant\/} Thom class
	for the bundle
	$\ContactNormalUnst\mu{\pi^*\lambda}\to\ContactCirc_{\pi^*\lambda}$
	via the map
	$$
		\BorelHomology^*_G(\Conley_G(\ContactCirc_{\pi^*\lambda}))
		\to\RedHomology^*(\Conley(\ContactCirc_{\pi^*\lambda})).
	$$
}
\Cor{
	If the equivariant contact invariant
	$\psi_G(\pi^*\lambda,\pi^*g)$
	vanish,
	then
	{\it the non-equivariant contact invariants
	of both $Y$ and $Y/G$ shall also vanish,}
	$$
		\psi(\pi^*\lambda,\pi^*g)=0,
		\quad
		\psi(\lambda,g)=0.
	$$
}
\Proof{
	Corollary of {\EquivariantDeterminesDown}
	and \EquivariantMapsToUp.
}
In order to seek computable examples,
the author shall make use of the following definition
first used in the work of \LinLipnowskiHyperbolic.
\Def{
	By saying that
	$Y$ is a 
	{\it minimal L-space,}
	one means that the perturbed Seiberg-Witten
	equations admit no irreducible solutions
	for some perturbation of arbitrarily small norm.
}
\Eg{
	Perhaps the best known examples
	of minimal L-spaces
	are the {\it elliptic manifolds\/}
	due to their metrics of positive scalar curvature
	(cf.\ \KronheimerMrowkaMonopolesAndThreeManifolds,
	\S 22.7).
}
\Eg{
	Another well known example
	of a minimal L-space
	is the {\it Hantzsche-Wendt\/} manifold,
	that is,
	the unique {\it flat\/} rational homology $3$-sphere
	(cf.\ \KronheimerMrowkaMonopolesAndThreeManifolds,
	\S 37.4).
}
\Eg{
	By recent work of \LinSolv,
	it is known that all rational homology $3$-spheres
	which have {\it sol geometry\/}
	are minimal L-spaces.
}
\Eg{
	According to the surprising work of \LinLipnowskiHyperbolic,
	it is known that certain well known {\it hyperbolic\/}
	rational homology $3$-spheres are also minimal L-spaces.
}
\Rem{
	Assume that $X$ be a minimal L-space.
	In this case,
	it is easy to see that
	$$
		\SWF(X,\SpinCStruct,g)
		=\Sigma^{-W^{(-\mu,0)}}S
	$$
	where $S\in\Spectra{}{\R^\infty}$
	is the sphere spectrum
	(cf.\ \ManolescuSeibergWittenFloerAlt,
	\S 10,
	example about the Poincar\'e homology sphere).
	Likewise,
	supposing $Y$ to be a minimal L-space
	upon which $G$ freely act,
	it is just as easy to see that
	$$
		\SWF_G(Y,\pi^*\SpinCStruct,\pi^*g)
		=\Sigma^{-W^{(-\mu,0)}}S_G,
	$$
	where
	$S_G\in\Spectra{G}{\Universe'}$
	is the $G$-equivariant sphere spectrum
	for the universe $\Universe'$
	seen as a complete $G$-universe;
	this is because,
	there being a unique solution
	to the Seiberg-Witten equations,
	$G$ can only act trivially on it
	and,
	therefore,
	the action on the Conley index
	will come entirely from the linear action
	on $W^\mu$.
}
\Rem{
	In order to apply the localization theorem,
	it is easiest to work at the level of spaces
	instead of spectra
	and consider the inclusion of fixed points
	(cf.\ {\SpectrumFixedPoints})
	$$
		\Conley(
			S^\mu_{\lambda,r},
			\FDASWFlow^\mu_{\lambda,r}
		)
		\to
		\Conley_G(
			S^\mu_{\pi^*\lambda,r},
			\FDASWFlow^\mu_{\pi^*\lambda,r}
		).
	$$
	This is simply {\it the inclusion
	of fixed points
	of a $G$-representation sphere,}
	exactly
	as was studied in
	{\LocalizationForRepresentationSphereModOddPrime}
	and {\LocalizationForRepresentationSphereModTwo}.
	The crux of the matter
	is to deduce what this says
	about the contact invariants at play.
	The first observation
	is the following.
}
\Prop{
	\Label\NonDesuspendedLocalizedEquivariantMapsToDown
	Let $S\subset \Homology^*(\B G;\Z/p\Z)$
	be the set of Euler classes of $G$-representations
	having no trivial summands.
	The map
	$$
		S^{-1}\BorelHomology_G^*(
			\Conley_G(
				S^\mu_{\pi^*\lambda,r},
				\FDASWFlow^\mu_{\pi^*\lambda,r}
			)
		)
		\to
		S^{-1}\BorelHomology_G^*(
			\Conley(
				S^\mu_{\lambda,r},
				\FDASWFlow^\mu_{\lambda,r}
			)
		)
	$$
	on localized Borel cohomology
	induced by the inclusion of fixed points
	maps
	$\psi_G(\pi^*\lambda,\pi^*g,\mu)$
	to $\psi(\lambda,g,\mu)\otimes e_G(\ContactFNVBF)$,
	where these classes are being seen now
	as classes in localized Borel cohomology.
}
\Proof{
	Follows from
	{\NonDesuspendedEquivariantMapsToDown}
	after localizing.
}
\Thm{
	\Label\MinimalLSpaceDownDeterminesEquivariant
	Suppose $\pi:Y\to Y/G$
	be a regular $p$-fold covering
	of minimal L-spaces
	where $p$ be prime
	and let $\lambda$ be a contact form on $Y/G$.
	It follows that
	$\psi_G(\pi^*\lambda,\pi^*g)=0$
	if and only if
	$\psi(\lambda,g)=0$.
}
\Proof{
	It suffices to work with the invariants
	$\psi_G(\pi^*\lambda,\pi^*g,\mu)$
	and $\psi(\lambda,g,\mu)$
	dependent on the spectral cut-off parameter
	as $\psi_G(\pi^*\lambda,\pi^*g)$
	and $\psi(\lambda,g)$
	are just grading-shifted versions of those.
	The localization theorem asserts that
	$$
		S^{-1}\BorelHomology_G^*(
			\Conley_G(
				S^\mu_{\pi^*\lambda,r},
				\FDASWFlow^\mu_{\pi^*\lambda,r}
			)
		)
		\to
		S^{-1}\BorelHomology_G^*(
			\Conley(
				S^\mu_{\lambda,r},
				\FDASWFlow^\mu_{\lambda,r}
			)
		)
	$$
	is an isomorphism.
	Since
	$\Conley_G(S^\mu_{\pi^*\lambda,r}, \FDASWFlow^\mu_{\pi^*\lambda,r})$
	is a $G$-representation sphere,
	its Borel cohomology
	is a free one-dimensional
	$\BorelHomology_G^*(S^0)$-module;
	hence,
	$\psi_G(\pi^*\lambda,\pi^*g,\mu)$
	cannot be a torsion element
	(with respect to the
	$\BorelHomology^*_G(S^0)$-module
	structure)
	of
	$$
		\BorelHomology_G^*(
			\Conley_G(
				S^\mu_{\pi^*\lambda,r},
				\FDASWFlow^\mu_{\pi^*\lambda,r}
			)
		).
	$$
	The result then follows from 
	{\NonDesuspendedLocalizedEquivariantMapsToDown}
	and,
	again,
	the observation that $e_G(\ContactFNVBF)\neq0$
	since $\ContactFNVBF^G=0$.
}
\Thm{
	Suppose $\pi:Y\to Y/G$
	be a regular $p$-fold covering
	of minimal L-spaces
	where $p$ be prime
	and let $\lambda$ be a contact form on $Y/G$.
	If the contact invariant 
	$\psi(\lambda)$
	with $\Z/p\Z$ coefficients
	on $Y/G$
	vanish,
	the contact invariant
	$\psi(\pi^*\lambda)$
	with $\Z/p\Z$ coefficients
	on $Y$
	shall also vanish.
}
\Proof{
	Noting,
	again,
	that it suffices to work with
	the analogous invariants
	dependent on the metric
	and the spectral cut-off parameter,
	combine 
	{\MinimalLSpaceDownDeterminesEquivariant}
	and \EquivariantMapsToUp.
}
\Def{
	For a rational homology $3$-sphere $Y$
	with $\SpinC$ structure
	$\SpinCStruct$,
	the {\it Fr\o yshov invariant,}
	$h(Y,\SpinCStruct)$,
	is defined as follows.
	Recall firstly that,
	as a $\Z[U]$-module,
	$\HMFrom^*(Y,\SpinCStruct)$
	splits as a rank one
	free summand
	and a torsion summand.
	Let $h(Y,\SpinCStruct)$
	be minus one half of the minimal degree,
	with respect to the absolute $\Q$-grading,
	in which the free summand
	be non-zero.
}
\Rem{
	There is another well known invariant
	extracted from the $\Q$-grading
	of a Floer theory;
	in this case,
	from the Heegaard Floer theory.
	This one is denoted $d(Y,\SpinCStruct)$
	and was first defined by
	{\OzsvathSzaboAbsolutely}
	in a similar fashion
	but referencing $\HF^+_*(Y,\SpinCStruct)$
	instead of $\HMFrom^*(Y,\SpinCStruct)$.
	According to the corpus
	that has equated Heegaard Floer homology,
	monopole Floer cohomology
	and embedded contact homology,
	the invariants $h$ and $d$,
	in fact,
	hold the same information
	but,
	for historical reasons,
	they are related by
	$-2d(Y,\SpinCStruct)=h(Y,\SpinCStruct)$.
	This is achieved
	via the isomorphisms
	that preserve the
	absolute grading of the Floer theories
	by {\it hyperplane fields\/}
	(vid.\ \CristofaroGardinerAbsoluteAlt,
	\RamosAbsoluteAlt)
	and the observation that the $\Q$-grading,
	in each case,
	is recovered through Gompf's
	$d_3$ invariant of the hyperplane fields,
	perhaps plus one half
	depending on one's conventions.
	The advantage of working with the
	Ozsv\'ath-Szab\'o invariant
	is that it is readily computable
	for many important classes of manifolds.
}
\Thm{
	\Label\CoveringMinimalLSpacesPsiUpNonZero
	Suppose $\pi:Y\to Y/G$
	be a regular $p$-fold covering
	of minimal L-spaces
	where $p$ be prime
	and let $\lambda$ be a contact form on $Y/G$.
	If $\psi(\lambda)\neq 0$
	and
	$
		\deg\psi(\pi^*\lambda)
		=-2h(Y,\pi^*\SpinCStruct_\lambda)
	$,
	it follows that $\psi(\pi^*\lambda)\neq 0$.
}
\Proof{
	As before,
	work with the metric
	and spectral cut-off
	dependent versions.
	Then,
	$\psi(\lambda,g,\mu)\neq 0$
	implies
	$\psi_G(\pi^*\lambda,\pi^*g,\mu)\neq 0$
	by localization.
	Meanwhile,
	the map
	$$
		\BorelHomology_G^*(
			\Conley_G(
				S^\mu_{\pi^*\lambda,r},
				\FDASWFlow^\mu_{\lambda,r}
			)
		)
		\to
		\RedHomology^*(
			\Conley(
				S^\mu_{\pi^*\lambda,r},
				\FDASWFlow^\mu_{\lambda,r}
			)
		),
	$$
	which occurs in {\EquivariantMapsToUp},
	can be written in terms of each degree as
	$$\matrix{
		\cdots &
		0 &
		\Z/p\Z &
		\Z/p\Z &
		\cdots
		\cr
		&
		\downarrow &
		\downarrow &
		\downarrow &
		\cdots
		\cr
		\cdots &
		0 &
		\Z/p\Z &
		0 &
		\cdots,
	}$$
	where the central map $\Z/p\Z\to\Z/p\Z$
	is an isomorphism.
	Hence,
	the invariant $\psi(\pi^*\lambda,\pi^*g,\mu)$
	is non-zero
	precisely
	when
	$$
		\RedHomology^{
			\deg\psi(\pi^*\lambda,\pi^*g,\mu)
		}
		(
			\Conley(
				S^\mu_{\pi^*\lambda,r},
				\FDASWFlow^\mu_{\lambda,r}
			)
		)
		=\Z/p\Z.
	$$
	Or,
	equivalently,
	after desuspensions,
	$$
		\HMTilde^{\deg\psi(\pi^*\lambda)}(
			Y,\pi^*\SpinCStruct
		)=\Z/p\Z.
	$$
	\par
	Now,
	to relate this to the Fr\o yshov invariant,
	recall that,
	for an L-space,
	the value of $-2h(Y,\pi^*\SpinCStruct_\lambda)$
	is the grading of the {\it only\/} degree
	in which
	$\HMTilde^*(Y,\pi^*\SpinCStruct_\lambda)$
	is non-trivial.
}
\Rem{
	Note that
	$\deg\psi(\pi^*\lambda)=-2h(Y,\pi^*\SpinCStruct_\lambda)$
	is also an obvious {\it necessary\/} condition for
	$\psi(\pi^*\lambda)\neq 0$.
	Therefore,
	the above theorem can be seen as strengthening this
	to a necessary and
	{\it sufficient\/}
	condition.
}
\Thm{
	Suppose $\pi:Y\to Y/G$
	be a regular $p$-fold covering
	of minimal L-spaces
	where $p$ be prime
	and let $\lambda$ be a contact form on $Y/G$.
	If $\deg\psi(\pi^*\lambda)<-2h(Y,\pi^*\SpinCStruct_\lambda)$,
	then $\psi(\lambda)=0$.
}
\Proof{
	Note
	$
		\deg\psi_G(\pi^*\lambda,\pi^*g)
		=\deg\psi(\pi^*\lambda,\pi^*g)
	$,
	so,
	if $\deg\psi(\pi^*\lambda)<-2h(Y,\pi^*\SpinCStruct\lambda)$,
	then $\psi_G(\pi^*\lambda,\pi^*g)=0$
	because
	$
		\BorelHMTilde_G^{\deg\psi_G(\pi^*\lambda,\pi^*g)}(
			Y,\pi^*\SpinCStruct_\lambda,\pi^*g
		)=0
	$.
	The result now follows from {\EquivariantDeterminesDown}.
}
\Rem{
	Recall that the grading of the contact invariant
	is given simply by the
	$3$-dimensional obstruction theoretic invariant
	of hyperplane fields as
	$$
		\deg \psi(\lambda)
		=d_3(\Ker\lambda)
		+{1\over2}.
	$$
	The reader can find a precise definition of $d_3$,
	originally due to \GompfHandlebody,
	in the following section,
	which shall de dedicated to understanding its behaviour
	under coverings.
}
\Thm{
	(\MatkovicClassificationAlt;
	\GhigginiTightAlt;
	\GhigginiLiscaStipsiczClassificationAlt;
	\WuLegendrianAlt)
	The tight contact structures
	on a small Seifert fibred L-space
	all have non-vanishing contact invariant
	and no pair of non-isotopic tight contact structures
	share the same $\SpinC$ structure.
}
\Rem{
	Elliptic manifolds are
	small Seifert fibred L-spaces,
	and this theorem 
	allows one to state a stronger version
	of the above results.
}
\Cor{
	Suppose $\pi:Y\to Y/G$
	be a regular $p$-fold covering
	of {\it elliptic manifolds\/}
	where $p$ be prime
	and let $\lambda$ be a {\it tight\/} contact form on $Y/G$.
	It follows that $\pi^*\lambda$ is tight
	if and only if
	$d_3(\pi^*\lambda)+{1\over2}=-2h(Y,\pi^*\SpinCStruct_\lambda)$.
}
\Rem{
	One can avoid computing the Fr\o yshov invariant here
	if one know precisely the list
	of tight contact structures on $Y$
	together with their corresponding
	$\SpinC$ structures and $d_3$ invariants.
	In that event,
	it suffices for one to compute
	$\pi^*\SpinCStruct_\lambda$
	and $d_3(\pi^*\Ker\lambda)$
	and compare with the entries of that list.
	If there be a match,
	$\pi^*\lambda$ {\it must\/} be tight,
	if not, $\pi^*\lambda$ must be overtwisted.
	This may be useful in some cases
	where computing Fr\o yshov invariants is difficult
	and the classification of tight contact structures
	is known by other methods,
	but it should be noted that,
	in many cases of interest,
	the calculations required
	to classify tight contact structures
	on small Seifert fibred L-spaces
	require comprehensive knowledge of the Fr\o yshov invariants
	(cf.\ \MatkovicClassificationAlt).
	In summary,
	this can be phrased as follows.
}
\Cor{
	\Label\EllipticLiftHomotopicTight
	Suppose $\pi:Y\to Y/G$
	be a regular $p$-fold covering
	of {\it elliptic manifolds\/}
	where $p$ be prime
	and let $\lambda$ be a {\it tight\/} contact form on $Y/G$.
	It follows that $\pi^*\lambda$ is
	{\it isotopic\/}
	to a given tight contact structure
	-- and therefore is also tight --
	if and only if it is
	{\it homotopic\/}
	to said tight contact structure.
}
Use $\HantzscheWendt$
to denote the
{\it Hantzsche-Wendt\/}
manifold,
that is,
the unique rational homology $3$-sphere
with Euclidean geometry.
A curious fact about $\HantzscheWendt$
is that there are
cyclic {\it self-coverings\/}
$\HantzscheWendt\to\HantzscheWendt$
of any odd degree
(vid.\ \ChelnokovMednykhHantzscheWendtAlt).
As $\HantzscheWendt$ is also known
to be a minimal L-space,
one could hope to apply the above methods here.
Unfortunately,
this fails because these self-coverings
are never {\it regular\/}
due to the first homology of $\HantzscheWendt$
having even order
whilst the coverings odd order.
\par
The case of the rational homology spheres
admitting {\it Sol geometry\/}
may be one where the methods developed here can work.
The problem,
however,
lies in the limited knowledge available
about the contact topology of such manifolds;
there is no known classification
of the isotopy classes of tight contact structures.
Hence,
before being able to apply
{\CoveringMinimalLSpacesPsiUpNonZero}
in this case,
one would have to produce tight contact structures
and develop a way to compute the respective
Fr\o yshov invariants.
The author intends to pursue that route
in later work.
\par
A similar story can be told
of the few hyperbolic manifolds
known to be minimal L-spaces;
although the methods above make
non-trivial statements
about lifting tight contact structures,
one hardly has the ability to apply them
due to lack of information about the contact topology
of those manifolds.
One case where the methods developed here
{\it can be\/} applied
is the double covering 
of the hyperbolic manifold
${\tt m007(3,1)}$
by ${\tt m036(-3,2)}$,
where the notation being used
is that of the Hodgeson-Weeks census
of small volume closed hyperbolic $3$-manifolds
from the well known {\it SnapPy\/} software
of Culler,
Dunfield,
Goerner
\& Weeks.
That both these manifolds are minimal L-spaces
is the consequence of the work of
\LinLipnowskiHyperbolic,
Theorem 1.
\par
This leaves the elliptic manifolds
currently as the best place
in which to seek to perform concrete computations.
In the next few sections,
the author shall proceed to develop
certain topological techniques
which shall be required for these calculations.
The goal shall be to apply
{\EllipticLiftHomotopicTight}
to find certain
tight contact structures
which lift to tight contact structures
by showing simply
that the lift have the same homotopy
theoretic invariants
as a known tight contact structure.
In order to do that,
one must address the problem
of computing the
obstruction theoretic invariants
of the lifted contact structure.
There are two obstruction theoretic invariants:
the $\SpinC$ structure
and the $d_3$ invariant of Gompf.
Lifting either of them
provides a challenge of its own
that must be overcome.
The next two sections shall deal
with these two problems.

\Sect{The d3 Invariant and Finite Coverings}
Consider the following problem.
Say $\lambda$ be a contact form on a rational homology $3$-sphere $Y/G$
regularly and finitely covered by $\pi:Y\to Y/G$.
Given the value of $d_3(\Ker\lambda)$
is it possible to say something
about the value of $d_3(\pi^*\Ker\lambda)$?
This problem was studied by {\KhuzamLifting}
via use of the
{\it $G$-signature theorem.}
In summary,
one can solve the problem by seeking
an almost complex $4$-manifold-with-boundary
having the contact manifold $(Y,\pi^*\lambda)$
as its almost complex boundary
and extending the $G$-action into its interior.
The action need not remain free inside the
$4$-manifold-with-boundary;
it need only have a properly embedded closed surface
as its $G$-fixed points.
Having such an {\it equivariant almost complex filling,}
one can infer the lifting behaviour of the $d_3$ invariant
of the contact form $\lambda$.
\par
Producing such a $4$-manifold-with-boundary
is typically a difficult task.
The goal of this section
is to describe a method
that the author devised
by use of Kirby calculus
to construct such a filling.
The method consists,
in essence,
of the following.
One starts with a Kirby diagram of $Y/G$,
which also defines
a $4$-dimensional $2$-handlebody
having $Y/G$ as its boundary.
Then,
one removes
a set of $2$-handles
judiciously
so that the resulting handlebody
evidently have a {\it branched\/} covering
with branching locus a properly embedded surface
with boundary at the curves
along which the removed $2$-handles were originally attached.
Next,
one computes this branched covering,
and proceeds to {\it equivariantly\/} attach $2$-handles
in the hope that,
in the quotient,
these $2$-handles be precisely the $2$-handles
that were removed initially.
This procedure is not always easy to carry out,
but,
in simple cases such as when $Y/G$ be a surgery on a knot,
it can be readily done.
One must also be careful with the almost complex structure
on the filling.
To achieve that,
one considers {\it Legendrian\/} Kirby diagrams
for the contact manifold $(Y/G,\lambda)$.
\par
As an application of the method,
the author shall calculate the $d_3$ invariant
of the lift of a tight contact structure
of the $(-8)$-surgery on the left-handed trefoil
via a double covering.
This is an example of a {\it prism\/} manifold,
therefore elliptic,
so one can apply the results of the preceding section
to try and determine if the lifted contact structure
remains tight.
Recall that,
for the lift to be tight,
one needs the lift to have an appropriate pair
of $d_3$ invariant {\it and\/} $\SpinC$ structure.
It shall be shown that it does indeed have the appropriate $d_3$ invariant,
whereas,
understanding the behaviour of the $\SpinC$ structures under coverings
is the content of the next section.
\par
To begin,
the author shall recall the relevant definitions
and the standard results that shall be needed.
Let $Y$, $G$, $\pi$, $g$ and $\lambda$ be as in the previous section.
\Def{
	A contact $3$-manifold $(Y,\lambda)$
	is called the
	{\it almost complex boundary\/}
	of an almost complex $4$-manifold-with-boundary $(M,J)$
	if $\partial M=Y$
	and $\Ker\lambda = \Tan Y\cap J\Tan Y$.
}
\Def{
	(\GompfHandlebodyAlt)
	Given a contact structure $\Ker\lambda$
	on a rational homology $3$-sphere $Y$,
	the $d_3$ invariant is defined as follows.
	Choose some almost complex
	$4$-manifold-with-boundary $(M,J)$
	such that $Y$ be its almost complex boundary.
	Then,
	define
	$$
		d_3(\Ker\lambda)
		:={1\over4}\left(
			c_1(M,J)^2
			-2\chi(M)
			-3\sigma(M)
		\right)
	$$
	where the notation
	$c^2\in\Q$,
	for a class $c\in\Homology^2(X;\Z)$,
	is defined in
	\SquareInSecondCohomologyOfFourManiWithQHSThreeBoundary.
}
\Rem{
	The author shall sometimes write
	$d_3(\lambda)$
	instead of
	$d_3(\Ker\lambda)$
	as that is more convenient
	in the present context.
}
\Rem{
	Gompf originally called this invariant $\theta$
	and the factor of $1/4$ was not present
	in his definition.
	There is another convention
	sometimes found in the literature,
	particularly
	in the context of Heegaard Floer theory,
	wherein the value of $d_3$ has $1/2$ added to
	what was defined above.
	The appeal of doing so is that the grading
	of the contact invariant becomes $d_3$
	instead of $d_3+1/2$.
}
Now,
consider again the case
of finite coverings.
The following result of Khuzam
summarizes what can be said by means
of the $G$-signature theorem
about the lifting behaviour of $d_3$.
\Thm{\Label\KhuzamThm
	(\KhuzamLiftingAlt, Theorem 2)
	Suppose $\pi:Y\to Y/G$
	be an $m$-fold regular cyclic covering.
	Let $\lambda$ be a contact form on $Y/G$
	and $\pi^*\lambda$ its lift to $Y$.
	Suppose further that one have
	$\Pi:M\to M/G$
	an $m$-fold cyclic {\it branched\/} covering of
	almost complex $4$-manifolds-with-boundary
	where the branching locus be a closed surface
	$S\subset M$
	satisfying $S\cap \partial M=\emptyset$
	and where $(Y,\pi^*\lambda)$ be the almost complex boundary
	of $M$.
	Then,
	the $d_3$ invariants of $\lambda$ and $\pi^*\lambda$
	are related by
	$$
		d_3(\pi^*\lambda)
		=md_3(\lambda)
		+{3\over4}m\sigma(M/G)
		-{3\over4}\sigma(M)
		-{3\over4}\sum_{k=1}^{m-1}
			(S\cdot S)
			\csc^2(\gamma_k/2),
	$$
	where $\gamma_k$
	denotes the angle of rotation
	given by the action of the element
	$k\in\Z/m\Z$
	on the normal planes to the surface $S$
	when $X$ be equipped with a
	$\Z/m\Z$-invariant metric.
	In fact,
	this metric can always be chosen
	so that $\gamma_k=2\pi k/m$.
}
\Rem{
	In the present \Document,
	the author follows the convention
	that the lens space $L(p,1)$
	is the one given by $(-p)$-surgery
	on the unknot in $S^3$.
}
\Eg{
	\Label\ExampleOfLensSpacesBoundingDiskBundles
	In \KhuzamLifting,
	the only example studied is that of the universal covering
	$S^3\to L(p,1)$
	of the lens space $L(p,1)$.
	As a first example here,
	the author shall extend that to the lens covering
	$L(p,1)\to L(mp,1)$.
	One can construct the branched covering
	$\Pi:M\to M/G$
	as follows.
	Let $M$ be the disk bundle over $S^2$
	of Euler class $-p$
	so that $\partial M= L(p,1)$.
	Note,
	$G:=\Z/m\Z$ acts on $M$
	by rotating the disk fibres.
	The action is free away from the zero section,
	which is fixed;
	therefore, $S:=M^G\cong S^2$.
	The manifold $M/G$
	is no other than the disk bundle of Euler class $-mp$ over $S^2$,
	so its boundary is $L(mp,1)$.
	Both $M$ and $M/G$ have almost complex structures
	coming from their disk bundle structures
	and $\Pi$ is pseudoholomorphic
	with respect to these.
	Clearly,
	this is precisely the scenario of {\KhuzamThm}
	and it becomes possible to compute
	$d_3(\pi^*\lambda)$
	from 
	the value of $d_3(\lambda)$.
	The lens space $L(mp,1)$
	admits precisely $mp-1$ tight contact forms 
	$\lambda_1,\ldots,\lambda_{mp-1}$
	(vid.\ {\HondaClassificationAlt} or {\GirouxStructuresAlt})
	satisfying
	$$
		d_3(\lambda_k)
		={1\over4}\left(
			-mp - 1 + 4k - {4k^2\over mp}
		\right).
	$$
	Noting that $S\cdot S=-p$,
	one computes
	$$
		\eqalign{
			&d_3(\pi^*\lambda_k)
			\cr
			=&
			{m\over4}
			\left(
				-mp - 1 + 4k - {4k^2\over mp}
			\right)
			+{3\over4}m(-1)
			-{3\over4}(-1)
			-{3\over4}
			\sum_{j=1}^{m-1}(-p)\csc^2(j\pi/m)
			\cr
			=&
			-m
			+km - {k^2\over p}
			+{3\over4}
			-{p\over4}
			\cr
			=&
			{1\over4}
			\left(
				-p
				+3
				+4m(k-1)
				-{4k^2\over p}
			\right),
		}
	$$
	where the author used the well known identity
	$\sum_{j=1}^{m-1}\csc^2(j\pi/m)=(m^2-1)/3$,
	cf.\ \CauchyCours, Note VIII.
	Comparing this with the possible values
	of the $d_3$
	invariants for the tight contact structures
	on $L(p,1)$,
	one easily sees that the only
	values of $k$
	for which
	$d_3(\pi^*\lambda_k)$
	matches the value of the $d_3$ invariant
	of some tight contact structure on $L(p,1$)
	are $k=1$ and $k=mp-1$.
	These are the
	{\it universally\/}
	tight contact structures
	of $L(mp,1)$
	and,
	in any event,
	one knows that they lift to the universally tight
	contact structures on $L(p,1)$,
	so there is nothing interesting happening
	in this family of examples;
	that is,
	all virtually overtwisted contact structures
	immediately lift to overtwisted contact structures.
}
\Rem{
	As the reader may have foreseen,
	finding the equivariant almost complex filling
	of the contact $G$-manifold $(Y,\pi^*\lambda)$
	in order to apply {\KhuzamThm},
	is not typically an easy task.
	The remainder of this section,
	shall introduce a method
	for producing examples
	of such fillings
	by use of Kirby calculus.
	The rationale is the following.
	Firstly,
	one constructs a {\it branched covering\/}
	of $4$-manifolds-with-boundary
	where the branching locus is allowed
	to intersect the boundary;
	secondly,
	a set of $2$-handles is {\it equivariantly\/} glued
	in order to cap off the branching locus
	thereby making the boundary freely acted
	by the deck transformations.
}
\Rem{
	Such matters shall require working with knot diagrams.
	The author prefers to draw what are called
	{\it grid diagrams.}
	These are knot diagrams
	where only vertical and horizontal lines occur,
	and apparent corners should be understood
	as being arcs with of a very small radius.
	This style of diagram
	has the advantage of bringing
	isotopies and Reidemeister moves
	more evidently into the realm of combinatorics.
	The following figure
	shows the example of the left-handed trefoil
	as a grid diagram.
}
\Fig{
	\Label\LeftTrefoil
	\vbox{%
	\baselineskip = 0pt%
	\hbox{\hskip105bp}%
	\vbox{\vskip52.5bp}%
	\vtop{\vskip52.5bp}%
	\includegraphics{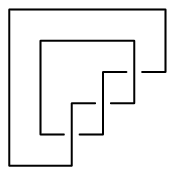}%
}%

}
\Rem{
	Now,
	follows a series of standard definitions
	and propositions
	which shall be needed
	in pursuing the goal of the remainder of this section.
}
\Def{
	By a {\it Legendrian knot or link\/}
	in a contact $3$-manifold $(Y,\lambda)$,
	one means a knot or link everywhere tangent
	to the contact structure $\Ker\lambda$.
}
\Rem{
	It is standard to represent Legendrian knots
	in $S^3=(\R^3)^+$,
	with its standard contact form
	$\lambda=\d z + x\d y$,
	via their
	{\it front projection diagrams.}
	Those are the knot diagrams one obtains
	when projecting a knot
	to the $yz$-plane.
	This leads to a diagram
	that possesses 
	no tangencies parallel to the $y$-axis
	but may have {\it cusps\/} pointing
	in directions parallel to the $z$-axis
	and whose transverse crossings
	must always
	have the strand
	passing {\it underneath\/}
	be the one for which
	the slope $\d z/\d y$ be {\it highest.}
	Conversely,
	any such diagram
	where all crossings be transverse
	is the front projection
	of a Legendrian knot.
	The grid knot diagrams,
	which the author favours,
	can readily be used to represent
	front projection diagrams
	if the following convention be agreed:
	the $y$-axis
	is understood to point in the
	southeast direction
	and the $z$-axis
	in the northeast direction.
	Here,
	the corners of the curve pointing
	in the
	{\it southwest\/}
	and 
	{\it northeast\/}
	directions
	are understood to be
	{\it smoothed\/}
	as before
	whereas those pointing
	in the
	{\it northwest\/}
	and
	{\it southeast\/}
	are understood to be
	{\it cusps.}
	Note that the diagram in {\LeftTrefoil}
	can be reinterpreted
	as a front projection diagram
	as it satisfies the required conditions.
	The following figure shows
	another {\it inequivalent\/} example
	of a Legendrian left-handed trefoil.
	When the author desire
	for the reader to interpret a particular
	diagram as Legendrian,
	he shall make it clear from the context;
	otherwise,
	the diagrams are to be understood simply
	as representing smooth links.
}
\Fig{
	\Label\LegendrianLeftTrefoilTBMinusSevenRotZero
	\vbox{%
	\baselineskip = 0pt%
	\hbox{\hskip90bp}%
	\vbox{\vskip45bp}%
	\vtop{\vskip45bp}%
	\includegraphics{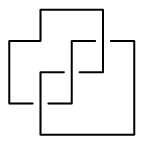}%
}%

}
\Def{
	The {\it contact framing\/}
	of a Legendrian knot $K$
	in $(Y,\lambda)$
	is the framing induced by
	the contact structure $\Ker\lambda$.
}
\Def{
	The {\it Thurston-Bennequin invariant,}
	$\TB(K)\in\Z$,
	of a Legendrian knot $K\subset S^3$
	is the value of the contact framing
	relative to the $0$-framing;
	that is,
	relative to the framing induced
	by a Seifert surface.
}
\Prop{
	For a Legendrian knot $K\subset S^3$,
	if $w(K)$ denote its {\it writhe\/}
	and $c(K)$ denote the number of {\it cusps\/}
	in one of its front projection diagrams,
	then
	$$
		\TB(K)=w(K)-{1\over2}c(K).
	$$
}
\Proof{
	Vid.\ e.g.\ \GeigesIntroduction, Proposition 3.5.9.
}
\Eg{
	\Label\LegendrianLeftTrefoilTBMinusSevenRotZeroHasTBMinusSeven
	The Legendrian left-handed trefoil $K$
	in {\LegendrianLeftTrefoilTBMinusSevenRotZero}
	has $\TB(K)=-7$.
}
\Def{
	The {\it rotation number,}
	$r(K)$,
	of an {\it oriented\/} Legendrian knot in $S^3$
	is the degree of the map
	$f:S^1\to S^1$ defined as follows.
	Take some parametrization
	$\gamma:S^1\to S^3$ of $K$.
	Let $\Sigma\subset S^3$
	denote a Seifert surface for $K$.
	Trivialize the standard contact structure of $S^3$
	over $\Sigma$
	yielding a bundle isomorphic to $\Sigma\times\R^2$.
	Since $K$ is Legendrian,
	the derivative $\gamma'$
	together with this trivialization
	defines a map $S^1\to\R^2\setminus\{0\}$.
	Given a deformation retraction
	$\R^2\setminus\{0\}\to S^1$,
	this defines the desired map $f:S^1\to S^1$.
}
\Def{
	For an oriented Legendrian knot
	$K\subset S^3$,
	denote,
	respectively,
	by $c_+(K)$ and $c_-(K)$
	the number upwardly and downwardly oriented cusps
	of a front projection diagram of $K$.
	In a Legendrian grid diagram
	of the sort being used by the author,
	$c_+(K)$
	is the number of corners
	which
	{\it start by pointing east
	and finish by pointing north\/}
	as one traverses it along the orientation of $K$,
	while $c_-(K)$
	is the number of corners
	which
	{\it start by pointing south
	and finish by pointing west.}
}
\Prop{
	For an oriented Legendrian knot $K\subset S^3$,
	the rotation number can be computed by
	$$
		r(K)={1\over2}\left(c_-(K)-c_+(K)\right).
	$$
}
\Proof{
	Vid.\ e.g.\ \GeigesIntroduction, Proposition 3.5.19.
}
\Eg{
	\Label\LegendrianLeftTrefoilTBMinusSevenRotZeroHasRotZero
	The Legendrian left-handed trefoil $K$
	in {\LegendrianLeftTrefoilTBMinusSevenRotZero}
	has $r(K)=0$.
}
\Def{
	A $3$-manifold $Y$
	given as surgery on a Legendrian link in $S^3$
	where the framing of a link component $K$
	is precisely {\it the contact framing minus 1\/}
	is called a {\it Legendrian surgery.}
}
A Legendrian surgery naturally inherits
a {\it tight\/} contact structure from $S^3$.
Moreover,
the $4$-manifold-with-boundary
obtained from the corresponding Kirby diagram
has a natural {\it complex structure;}
indeed,
it is a {\it Stein surface.}
For a reference,
vid.\ \GompfHandlebody.
This means that $d_3$ invariants
can be read off this Kirby diagram
as follows.
\Prop{(\GompfHandlebodyAlt)
	\Label\DThreeInvariantOfContactSurgery
	Let $(Y,\lambda)$ be a contact
	rational homology $3$-sphere
	given as Legendrian surgery on a Legendrian link
	$\bigsqcup_{i=1}^n K_i$
	in $S^3$.
	Use $L$ to denote the linking matrix
	of the link
	and $\sigma(L)$ its signature
	as a symmetric bilinear form.
	Then,
	$$
		d_3(\lambda)
		=
		{1\over4}
		\left(
			\sum_{i,j=1}^n
			r(K_i)
			(L^{-1})_{ij}
			r(K_j)
			-2(n+1)
			-3\sigma(L)
		\right).
	$$
}
\Rem{
	\Label\KirbyCoveringConstructionOne
	Let $(X,\lambda)$ be a contact rational homology $3$-sphere.
	Suppose one have a Kirby diagram
	for an almost complex $4$-manifold-with-boundary $(N,J)$
	consisting entirely of $2$-handles
	with its almost complex boundary
	being $\partial(N,J)=(X,\lambda)$.
	Label the $2$-handles of $N$ as
	$\{h_i\mid i\in I\}$
	and their respective attaching knots in $S^3$
	as $\{K_i\mid i\in I\}$.
	Now,
	consider a
	{\it subset\/}
	of these $2$-handles,
	denote it $\{h_i\mid i\in I'\}$,
	and
	{\it remove\/}
	them from the Kirby diagram.
	This yields
	{\it another\/}
	almost complex $4$-manifold-with-boundary,
	call it $(N',J')$,
	having as its $2$-handles the set
	$\{h_i\mid i\in I\setminus I'\}$;
	denote its almost complex boundary by
	$(X',\lambda'):=\partial(N',J')$.
	Consider the link
	$L:=\sqcup_{i\in I'}K_i$
	as a link in $X'$.
	Fix some
	{\it pseudoholomorphic curve\/}
	$S'\subset N'$
	with $\partial S'=L$
	and $\Interior S'\subset\Interior N'$.
	Now,
	the $p$-fold
	{\it branched covering\/}
	of $N'$ with branching locus $S'$ exists
	precisely when the homology class
	$[S']\in\Homology_2(N',\partial N';\Z)$
	be divisible by $p$.
	Assume that this be indeed the case.
	Then,
	in fact,
	the branched cover,
	call it $M'$,
	can be made to come 
	with an almost complex structure
	and the covering map
	$\Pi':M'\to N'$
	can be made to respect the almost complex structures.
	Write $Y':=\partial M'$;
	this is a contact manifold.
	The restriction to the boundaries
	gives a branched covering of $3$-manifolds
	$\pi':Y'\to X'$
	having branching locus the link $L\subset X'$.
	Required now is the concept
	of {\it equivariant handle attachment.}
}
\Def{
	Consider a $2$-handle $h$
	as a copy of $D^2\times D^2$
	to be attached along $S^1\times D^2$.
	Regard $D^2\times D^2$
	as having the action of $G:=\Z/p\Z$
	defined by
	$e^{2\pi i/p}\cdot(z_1,z_2):=(z_1,e^{2\pi i/p}z_2)$.
	This action defines a $g$-fold branched covering
	$D^2\times D^2\to D^2\times D^2$
	with branching locus the core disk $D^2\times\{0\}$.
	Think of this as the {\it $G$-$2$-handle.}
	When one be in possession
	of a $G$-$4$-manifold-with-boundary $M$
	with a $G$-fixed surface $S$
	having $\Interior S\subset\Interior M$
	and $\partial S\subset\partial M$,
	define the
	{\it equivariant handle attachment\/}
	of $h$ along a component $K$ of $\partial S$
	to be the $G$-$4$-manifold-with-boundary
	$M\cup_\varphi h$
	given by gluing $h$ to $M$
	via an {\it equivariant framing\/}
	$\varphi$ of $K$;
	that is,
	an equivariant diffeomorphism onto its image
	$\varphi:S^1\times D^2\to M$
	such that its image be a tubular neighbourhood of $K$.
}
\Rem{
	Notice that the resulting manifold $M\cup_\varphi h$
	has one of the boundary components
	of the $G$-fixed surface $S$
	{\it capped\/} by the $G$-fixed core disk of $h$.
	Hence,
	$M\cup_\varphi h$
	has one less $G$-fixed circle on its boundary
	compared to $M$.
}
\Rem{
	\Label\CappingBranchingLocus
	Returning now to the context of
	{\KirbyCoveringConstructionOne},
	one desires to equivariantly attach $2$-handles
	to the manifold $M'$
	in order to convert the {\it branched\/} covering
	of $3$-manifolds
	$\pi:Y'\to X'$
	to a {\it genuine\/} covering
	$\pi:Y\to X$,
	where $X$ be the $3$-manifold
	with which one has started.
	Whether this is possible or not
	is a matter about the framings
	of the handles $\{h_i\mid i\in I'\}$
	which were initially removed from $N$
	resulting in the manifold $N'$.
	Consider a fixed $i\in I'$.
	Now,
	equivariantly attach a $2$-handle
	$\tilde h_i$
	along the knot $(\Pi')^{-1}(K_i)$
	with framing determined by an integer
	$m_i\in\Z$.
	One hopes to be able to choose
	$m_i$ so that the corresponding
	{\it quotient\/} handle $\tilde h_i/G$
	have the same framing as the handle $h_i$ of $N$;
	thereby allowing one to identify them.
	The main issue is that determining the behaviour
	of framings under the quotient
	turns out to be a subtle matter.
}
Instead of considering this general case,
simplify the scenario as follows.
Consider instead a manifold $N$
as above
but consisting of a {\it single\/}
$2$-handle $h$
attached along a knot $K\subset S^3$
with framing $n$.
Then,
one forms $N'$
by removing $h$
to obtain,
of course,
the disk $D^4$.
The manifold $M'$,
therefore,
is a familiar sort of manifold;
it is the branched covering of $D^4$
branched over the Seifert surface of $K$
with its interior pushed into the interior of $D^4$.
In this case,
it is an easy matter to determine
the behaviour of framings under the quotient.
If one equivariantly attach
a $2$-handle $\tilde h$
to $M'$
along $(\Pi')^{-1}(K)$
with framing $m$,
the framing of the quotient handle $\tilde h/G$
shall be $mp$.
Hence,
for the procedure to work,
one needed $n$ to be a multiple of $p$.
\par
The process of computing branched coverings
in Kirby calculus
is outlined in \GompfStipsiczFourManifolds,
\S 6.3.
In general,
it can be a difficult procedure
where one obtains a Kirby diagram
for the complement of a tubular neighbourhood
of the branching surface,
computes the genuine covering of the resulting manifold
and then glues handles
to fill in the hole left by the removal
of the branching surface.
There is a simpler approach,
discussed in
\RolfsenKnots,
which has the disadvantage
of requiring blow-ups to be performed.
By a blow-up
the author means
the formation of
the connected sum with a copy of
$\CPBar^2$ or $\CP^2$.
In the case of $\CPBar^2$,
almost complex complex structures
can be naturally carried to the blow-up;
in the case of $\CP^2$,
this is not the case.
Hence,
for the purposes being pursued
in this section,
one may only perform blow-ups
of the type where one takes
the connected sum
with a copy $\CPBar^2$.
The procedure is as follows.
Starting with the Kirby diagram of $N$
consisting of a single $2$-handle $h$
attached along a knot $K$,
perform a sequence of blow-ups
in order to untie $K$
without changing the framing $n$ of the handle $h$.
Now,
the branched covering that one needs
to compute is simply branched over a disk
since $K$ has become the unknot.
This is a significantly easier task
than the general case.
\Fig{
	\Label\MinusEightLeftTrefoilBlowUp
	\vbox{%
	\baselineskip = 0pt%
	\hbox{\hskip105bp}%
	\vbox{\vskip52.5bp}%
	\vtop{\vskip52.5bp}%
	\includegraphics{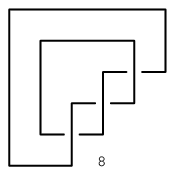}%
}%

	$\longrightarrow$
	\vbox{%
	\baselineskip = 0pt%
	\hbox{\hskip135bp}%
	\vbox{\vskip67.5bp}%
	\vtop{\vskip67.5bp}%
	\includegraphics{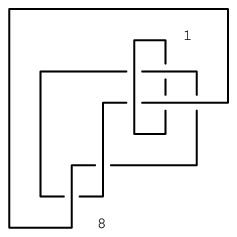}%
}%

}
Consider now a concrete example
in which to apply this procedure.
Start with the manifold $X$
given as $(-8)$-surgery
on the left handed trefoil knot $K$.
This manifold admits a
double covering,
which is the one that shall be studied here.
One can perform a single blow-up
to untie $K$
as depicted in
{\MinusEightLeftTrefoilBlowUp}.
Hence,
the manifold $N$ shall be the $2$-handlebody
defined by the right hand side Kirby diagram
of {\MinusEightLeftTrefoilBlowUp}.
After a series of isotopies,
one can achieve the form
on the left hand side of the following figure.
\Fig{
	\Label\MinusEightLeftTrefoilReadyForCover
	\vbox{%
	\baselineskip = 0pt%
	\hbox{\hskip105bp}%
	\vbox{\vskip67.5bp}%
	\vtop{\vskip67.5bp}%
	\includegraphics{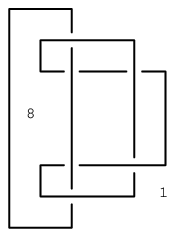}%
}%

	$\longrightarrow$
	\vbox{%
	\baselineskip = 0pt%
	\hbox{\hskip135bp}%
	\vbox{\vskip67.5bp}%
	\vtop{\vskip67.5bp}%
	\includegraphics{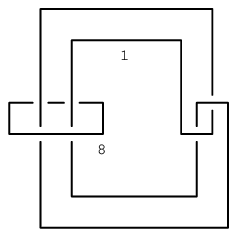}%
}%

}
After another isotopy,
one obtains the diagram
in right hand side of
{\MinusEightLeftTrefoilReadyForCover}.
What one desires to do now
is to {\it remove\/}
the $(-8)$-framed handle
in order to form the manifold $N'$
as depicted in the left hand side
of the following figure.
Once that be done,
one can read off the double branched covering
branched over the obvious Seifert disk
for the knot $K$.
The resulting manifold,
$M'$,
is depicted on the right hand side
of the following figure.
The more lightly stroked curve
indicates the handle which has been
{\it removed,}
which is,
therefore,
also the boundary of the branching locus.
\Fig{
	\Label\LeftTrefoilDoubleBranchedCover
	\vbox{%
	\baselineskip = 0pt%
	\hbox{\hskip135bp}%
	\vbox{\vskip67.5bp}%
	\vtop{\vskip67.5bp}%
	\includegraphics{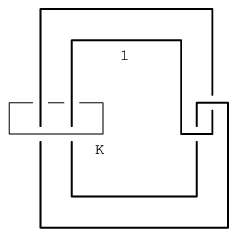}%
}%

	\vbox{
		\baselineskip=10pt
		\hbox to 1.5cm {\hfill $2:1$\hfill}
		\hbox to 1.5cm {\leftarrowfill}
	}
	\vbox{%
	\baselineskip = 0pt%
	\hbox{\hskip135bp}%
	\vbox{\vskip67.5bp}%
	\vtop{\vskip67.5bp}%
	\includegraphics{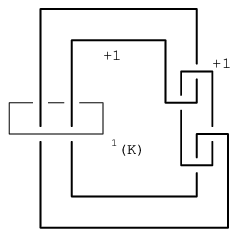}%
}%

}
\Rem{
	Notice that the framings
	behave under the branched covering
	in such a manner that the
	{\it blackboard framing\/}
	be preserved
	irrespective of the change in the writhes.
	This is why the $(-1)$-framed handle
	lifts to a pair of $(+1)$-framed handles.
}
Next,
one equivariantly attaches
a $2$-handle along 
$(\Pi')^{-1}(K)$
in order to cap off
the branching surface $S$
leading to the branched covering
of $4$-manifolds-with-boundary
seen in the following figure
where the branching locus no longer
intersects the boundary.
\par
\Fig{
	\Label\MinusEightLeftTrefoilDoubleCover
	\vbox{%
	\baselineskip = 0pt%
	\hbox{\hskip135bp}%
	\vbox{\vskip67.5bp}%
	\vtop{\vskip67.5bp}%
	\includegraphics{figs/minus_eight_left_trefoil_ready_for_cover.ps}%
}%

	\vbox{
		\baselineskip=10pt
		\hbox to 1.5cm {\hfill $2:1$\hfill}
		\hbox to 1.5cm {\leftarrowfill}
	}
	\vbox{%
	\baselineskip = 0pt%
	\hbox{\hskip135bp}%
	\vbox{\vskip67.5bp}%
	\vtop{\vskip67.5bp}%
	\includegraphics{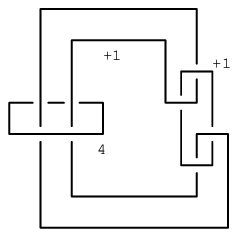}%
}%

}
\Prop{
	\Label\ContactPrismExampleDThreeOneQuarter
	The manifold $X$
	defined by $(-8)$-surgery on the left handed trefoil
	admits a contact form $\lambda$
	which lifts via a double covering
	$Y\to X$
	to a contact form $\pi^*\lambda$
	having $d_3(\pi^*\lambda)=1/4$.
}
\Proof{
	Define $N$ and $M$
	to be the $4$-manifolds-with-boundary
	on the left and right hand side
	of {\MinusEightLeftTrefoilDoubleCover}
	respectively.
	Here,
	the branching locus no longer intersects
	the boundaries 
	$X=\partial N$
	and
	$Y:=\partial M$;
	therefore,
	one can apply {\KhuzamThm}
	but only to the subset
	-- and this subset may be proper --
	of the contact structures on $X$
	which occur as almost complex boundaries
	of $N$.
	The author shall focus
	on a particular contact structure.
	Let $\lambda$ be the contact structure on $X$
	defined as the {\it Legendrian surgery\/}
	according to the Legendrian representative of $K$
	seen in
	{\LegendrianLeftTrefoilTBMinusSevenRotZero}.
	Recall,
	from
	{\LegendrianLeftTrefoilTBMinusSevenRotZeroHasTBMinusSeven},
	that the Thurston-Bennequin invariant
	of this Legendrian knot
	is $-7$;
	hence,
	the Legendrian surgery is indeed the
	topological $(-8)$-surgery.
	Since,
	in defining the $4$-manifold-with-boundary $N$,
	only a $(-1)$-blow-up was performed,
	this contact structure is still
	an almost complex boundary of $N$.
	From
	{\LegendrianLeftTrefoilTBMinusSevenRotZeroHasRotZero},
	one knows that the rotation number is $0$.
	With these data,
	one applies
	{\DThreeInvariantOfContactSurgery}
	to find that
	$$
		d_3(\pi^*\lambda)
		=
		{1\over4}
		\left(
			0
			-4
			+3
		\right)
		=-{1\over4}.
	$$
	Now,
	using {\KhuzamThm},
	one computes
	$$
		\eqalign{
			&d_3(\pi^*\lambda)
			\cr
			=&
			-{1\over2}
			+{3\over4}\left(
				+2\sigma\pmatrix{
					-8 &  0 \cr
					 0 & -1 \cr
				}
				-\sigma\pmatrix{
					-4 &  0 &  0 \cr
					 0 &  1 & -2 \cr
					 0 & -2 &  1 \cr
				}
				-\sum_{k=1}^2(-4)
				\csc^2\left({
					\pi k/m
				}\right)
			\right)
			\cr
			=&
			{1\over4},
			\cr
		}
	$$
	where the matrices
	are read off from {\MinusEightLeftTrefoilDoubleCover}.
}
\Prop{
	The manifold $X$
	defined by $(-8)$-surgery
	on the left handed trefoil
	is the {\it prism manifold\/}
	having Seifert fibration
	$(S^2;(1,-1),(3,2),(2,1),(2,1))$.
}
\Proof{
	The left-handed trefoil is a torus knot;
	hence,
	one can establish the Seifert invariants
	of $X$ from a well known theorem
	of \MoserElementary.
}
\Cor{
	The double cover $Y$ of $X$
	is the lens space $L(12,7)$.
}
\Proof{
	One can compute coverings
	directly from the Seifert invariants.
	The Seifert manifold
	$(S^2;(1,-1),(3,2),\allowbreak(2,1),(2,1))$
	has a horizontal
	double covering
	(meaning it lifts Seifert fibres
	to Seifert fibres
	$1$-to-$1$)
	with invariants
	$$
		(S^2;(1,-1),(1,-1),(3,2),(3,2),(1,1),(1,1))
	$$
	and this is one of the Seifert structures
	on the lens space $L(12,7)$.
}
\Rem{
	The lens space $L(12,7)$
	can be obtained as surgery
	on a chain of three unknots
	with framings respectively $-2$, $-4$ and $-2$.
	According to the classification
	of tight contact structures on lens spaces
	due,
	independently,
	to {\HondaClassification}
	and {\GirouxStructures},
	the manifold $Y$
	has precisely $3$ isotopy classes
	of tight contact structures.
	Each of these contact structures
	comes from Legendrian surgery
	on the three possible Legendrian stabilizations
	of this chain of unknots having
	the appropriate Thurston-Bennequin invariants.
	The three Legendrian surgery diagrams
	are depicted in the following figure.
}
\Fig{
	\Label\LensTwelveSevenContact
	\vbox{%
	\baselineskip = 0pt%
	\hbox{\hskip105bp}%
	\vbox{\vskip52.5bp}%
	\vtop{\vskip52.5bp}%
	\includegraphics{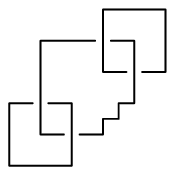}%
}%

	\vbox{%
	\baselineskip = 0pt%
	\hbox{\hskip97.5bp}%
	\vbox{\vskip48.75bp}%
	\vtop{\vskip48.75bp}%
	\includegraphics{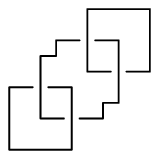}%
}%

	\vbox{%
	\baselineskip = 0pt%
	\hbox{\hskip105bp}%
	\vbox{\vskip52.5bp}%
	\vtop{\vskip52.5bp}%
	\includegraphics{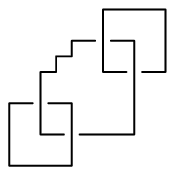}%
}%

}
\Prop{
	The three tight contact structures on $Y=L(12,7)$
	depicted in {\LensTwelveSevenContact}
	have $d_3$ invariants,
	respectively,
	equal to $-1/12$, $1/4$ and $-1/12$.
}
\Proof{
	Compute using 
	\DThreeInvariantOfContactSurgery.
}
\Rem{
	\Label\ContactPrismExampleDThreeLiftMatch
	In light of
	\ContactPrismExampleDThreeOneQuarter,
	the contact structure
	depicted in the middle of {\LensTwelveSevenContact}
	is the one that shall receive special attention
	as its $d_3$ invariant
	matches that of the lift
	of the tight contact structure
	on $X$
	via the double cover.
	In order to apply
	{\EllipticLiftHomotopicTight}
	and see that these two contact structure are
	one and the same,
	what remains to be shown
	is whether the $\SpinC$ structures
	also agree.
	This problem shall be studied in the next section.
}

\Sect{The Spin-C Structure and Finite Coverings}
This section shall deal with a fairly elementary problem
but one that turns out to be quite difficult
to solve in practice.
The problem is that of lifting $\SpinC$ structures
across finite coverings.
There are various methods
which one might try to use
in approaching this problem.
For instance,
in the case of Seifert manifolds,
there is a certain canonical $\SpinC$ structure
defined by the Seifert fibration.
This canonical $\SpinC$ structure
is pulled back naturally via coverings
and,
therefore,
it suffices to have a good description of it
upstairs and downstairs in the covering
in order to understand how
every other $\SpinC$ structure lifts.
This is so
because the set of $\SpinC$ structures on a manifold $Y$
is a $\Homology^2(Y;\Z)$-torsor,
and,
under the map
induced by a covering,
this torsor structure is preserved;
hence,
it suffices to understand
the induced map in the second cohomology
together with how one single $\SpinC$ structure lifts
in order to deduce the action on every other.
Another approach is to use contact topology.
Suppose one know of a certain
{\it universally\/} tight contact structure $\lambda$ on $Y$,
then,
its lift must always be universally tight
via any covering.
If there not be many of these upstairs,
this reveals information about what the lift
of the $\SpinC$ structure associated with $\lambda$ is.
This strategy is particularly useful in the case of lens spaces.
Another method is to turn to {\it Spin structures.}
To any Spin structure is associated
a $\SpinC$ structure in a natural way.
Hence,
it is suffices to know how a Spin structure
lifts via a covering
to deduce how the $\SpinC$ structures lift as well.
\par
This last approach shall be the one pursued in the present section.
Spin structures can be studied
with the aid of the Kirby calculus
and this fits well into the picture
of lifting $d_3$ invariants of the preceding section.
Indeed,
{\GompfStipsiczFourManifolds}
describe not only how to express Spin structures
on a $3$-manifold given by a Kirby diagram
by ascribing extra decorations to it,
but also how this description changes under Kirby moves.
The main feat of the current section
is to describe the manner in which this description of Spin structures
behaves under finite coverings.
\par
The double covering of the $(-8)$-surgery on the left-handed trefoil
considered in the previous section
shall continue to furnish examples in this section.
By the end,
it shall be established that the tight contact structure
whose lift was postulated to be tight
indeed has its $\SpinC$ structure lift to the correct one,
so that the theorems concerning the equivariant contact invariant
assert that the lift be tight.
\Def{
	Let $M$ be an $n$-manifold
	with $n\geq3$
	and use $M_k$
	to denote the $k$-skeleton of $M$
	in some cellular decomposition of $M$.
	By a {\it Spin structure\/} on $M$
	one means a trivialization of
	$\Tan M|_{M_1}$
	which extend to $\Tan M|_{M_2}$.
}
According to \GompfStipsiczFourManifolds,
Spin structures can be understood
in terms of Kirby calculus
in a rather convenient fashion
which shall be recalled next.
Consider a $3$-manifold $Y$
given as the boundary
of a $2$-handlebody $M$.
Denote by $\{K_i\subset S^3\mid i\in I\}$
the set of knots onto which the $2$-handles of $M$
are attached
and by $\{h_i\mid i\in I\}$
the corresponding $2$-handles.
\Def{
	Given a Spin structure $\SpinStruct$ on $Y$,
	define the class
	$w_2(M,\SpinStruct)\in\Homology^2(M,Y;\Z/2\Z)$
	as the obstruction
	to extending $\SpinStruct$
	from $Y$ to $M$.
}
\Rem{
	By standard obstruction theory,
	this class can be characterized
	in terms of the its evaluations
	on each of the classes
	$[D_i,\partial D_i]\in\Homology_2(M,Y;\Z/2\Z)$
	associated to the cocores
	$\{D_i\}$
	of the $2$-handles
	$\{h_i\}$.
	One can ask whether $\SpinStruct$ extends
	across $h_i$
	to a Spin structure on $h_i\cong D^2\times D^2$.
	It follows that
	$w_2(M,\SpinStruct)$
	evaluates on the class $[D_i,\partial D_i]$
	as $0$
	when $\SpinStruct$
	extends across $h_i$
	and as $1$
	when it does not.
}
\Def{
	An element $w$ of $\Homology^2(M,Y;\Z/2\Z)$
	is called {\it characteristic\/}
	when it gets mapped to the second Stiefel-Whitney class
	$w_2(M)\in\Homology^2(M;\Z/2\Z)$
	via the map
	$\Homology^2(M,Y;\Z/2\Z)\to\Homology^2(M;\Z/2\Z)$.
}
For each $i$,
suppose $F_i\subset M$ to be a closed surface
constructed by taking a Seifert surface
for the knot in $S^3$
along which the $2$-handle $h_i$ was attached,
pushing the interior of this Seifert surface
into the interior of $D^4$
and capping it
by gluing it to the core disk of the handle $h_i$.
Orient the closed surface $F_i$
in a manner compatible to the orientation of the
knot along which $h_i$ was attached.
Notice that the classes $[F_i]$
span $\Homology_2(M;\Z)$ as a $\Z$-module.
\Prop{
	A class $w\in\Homology^2(M,Y;\Z/2\Z)$
	is characteristic
	if and only if,
	for all $i\in I$,
	$$
		w([D_i,\partial D_i])
		=
		[F_i]\cdot [F_i]
		\quad \hbox{mod}\ 2.
	$$
}
\Proof{
	Follows from Wu's formula.
	Vid.\ \GompfStipsiczFourManifolds,
	Exercise 5.7.3.
}
\Prop{
	(\GompfStipsiczFourManifoldsAlt)
	The class
	$w_2(M,\SpinStruct)\in\Homology^2(M,Y;\allowbreak\Z/2\Z)$
	completely characterizes the Spin structure
	$\SpinStruct$ on $Y$.
	Conversely,
	given any characteristic element
	$w\in\Homology^2(M,Y;\Z/2\Z)$
	there is some Spin structure $\SpinStruct$ on $Y$
	such that $w_2(M,\SpinStruct)=w$.
}
\Rem{
	In a Kirby diagram,
	the author shall
	denote a Spin structure $\SpinStruct$
	by writing
	a $0$ or a $1$
	following the framing coefficient
	separated from it by a comma
	next to each $2$-handle.
	This signifies the value of the evaluation
	of $w_2(Y,\SpinStruct)$
	on the class in
	$H_2(M,Y;\Z/2\Z)$
	corresponding to the respective $2$-handle.
	The following figure shows
	this notation in the main example
	from the preceding section.
}
\Fig{
	\Label\MinusEightLeftTrefoilSpin
	\vbox{%
	\baselineskip = 0pt%
	\hbox{\hskip105bp}%
	\vbox{\vskip52.5bp}%
	\vtop{\vskip52.5bp}%
	\includegraphics{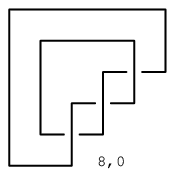}%
}%

}
\Rem{
	It is then possible to keep track
	of the Spin structure
	during the performance of Kirby moves.
	The rule when performing a handle slide
	of $h_i$ over $h_j$,
	where the respective components
	of $w_2(M,\SpinStruct)$
	be $n_i$ and $n_j$,
	is that $n_j$ changes
	precisely when $n_i=1$,
	whereas $n_i$ always stays unchanged.
	For a blow-up,
	the component of $w_2(M,\SpinStruct)$
	associated to the newly attached $(\pm 1)$-framed
	$2$-handle is always $1$.
	The next figure illustrates
	the result of both these operations.
}
\Fig{
	\Label\MinusEightLeftTrefoilBlowUpSpin
	\vbox{%
	\baselineskip = 0pt%
	\hbox{\hskip135bp}%
	\vbox{\vskip67.5bp}%
	\vtop{\vskip67.5bp}%
	\includegraphics{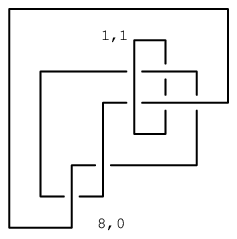}%
}%

}
The matter is now to consider how a Spin structure
lifts via a finite covering
in the Kirby calculus setting.
As in the previous section,
the procedure shall be carried out in two parts.
Firstly,
one performs a {\it branched\/} covering
over a $4$-manifold-with-boundary
obtained by judiciously removing $2$-handles
and then one equivariantly glues $2$-handles
in the branched covering to cap the branching locus.
In the first part,
one needs to take care of how a Spin structure
behaves near the {\it free\/} $2$-handles;
in the second part,
one needs to take care of how it behaves near
the {\it non-free\/} $2$-handles.
\par
Consider a $3$-manifold $X$,
the boundary
of a $4$-manifold-with-boundary
$N$ having only a $0$-handle
and a set of $2$-handles
$\{h_i\mid i\in I\}$.
Let $\SpinStruct$ be a Spin structure on $X$.
Denote by $\{K_i\subset S^3\}$
the set of attaching knots of the $2$-handles
and by $\{D_i\}$ the cocores.
Given a subset of $2$-handles
$\{h_i\mid i\in I'\}$,
denote by $N'$
the $4$-manifold-with-boundary
given by removing this set of $2$-handles from $N$;
denote by $X'$ its boundary.
Now,
{\it assume\/}
that there exist
the $p$-fold branched covering
$\Pi':M'\to N'$
branched at a Seifert surface $S'$
for the link
$L:=\bigsqcup_{i\in I'}K_i$
with its interior pushed into the interior
of $D^4$.
Let $Y':=\partial M'$
and $\pi':Y'\to X'$
be the restriction of $\Pi'$.
For each $2$-handle $h_i$ of $N'$,
the lift $(\Pi')^{-1}(h_i)$
consists of $p$ $2$-handles
of $M'$ freely permuted
by the deck transformations.
\Prop{
	\Label\LiftSpinFreeHandle
	For each handle $h$
	in $(\Pi')^{-1}(h_i)$
	with cocore $D$,
	$$
		w_2(Y',(\pi')^*\SpinStruct)([D,\partial D])
		=w_2(X',\SpinStruct)([D_i,\partial D_i]).
	$$
}
\Proof{
	Recall that $w_2(X',\SpinStruct)([D_i,\partial D_i])$
	is characterized by whether
	the Spin structure $\SpinStruct$
	extends to the whole interior of the $2$-handle $h_i$.
	That extension existing,
	the lifted Spin structure
	$(\pi')^*\SpinStruct$
	shall also extend over
	each of the preimage $2$-handles
	$h$ of $h_i$.
	Conversely,
	should the extension not exist downstairs,
	it cannot extend upstairs over any $h$ either.
}
\Eg{
	The next figure
	exemplifies {\LiftSpinFreeHandle}
	in the already familiar setting
	of the branched double cover
	over the left-handed trefoil $K$
	from the preceding section
	(cf.\ {\LeftTrefoilDoubleBranchedCover}).
}
\Fig{
	\Label\LeftTrefoilBranchedDoubleCoverSpin
	\vbox{%
	\baselineskip = 0pt%
	\hbox{\hskip135bp}%
	\vbox{\vskip67.5bp}%
	\vtop{\vskip67.5bp}%
	\includegraphics{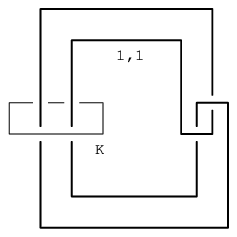}%
}%

	\vbox{
		\baselineskip=10pt
		\hbox to 1.5cm {\hfill $2:1$\hfill}
		\hbox to 1.5cm {\leftarrowfill}
	}
	\vbox{%
	\baselineskip = 0pt%
	\hbox{\hskip135bp}%
	\vbox{\vskip67.5bp}%
	\vtop{\vskip67.5bp}%
	\includegraphics{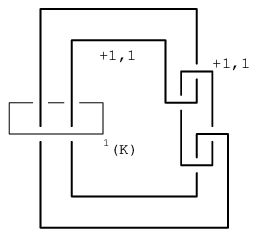}%
}%

}
Now,
one needs to understand
what happens to the non-free $2$-handles.
Proceeding as in
{\CappingBranchingLocus},
one equivariantly attaches
a $2$-handle $\tilde h_i$
along $(\Pi')^{-1}(K_i)$
for each $i\in I'$
in order to cap the branching locus
thereby defining the manifolds
$M:=M'\cup\bigcup_{i\in I'} \tilde h_i$
and $N:=M/G$.
Let $m_i$ denote the framing coefficient
of the equivariant handle $\tilde h_i$.
As before,
one hopes to be able to choose $m_i$
so that the original manifold $Y$
be $\partial N$.
Assume that this be the case.
The question then becomes:
given the value of $w_2(X,\SpinStruct)([D_i,\partial D_i])$
for an $i\in I'$,
what is the value of
$w_2(Y,\pi^*\SpinStruct)([\tilde D_i,\partial\tilde D_i])$
where $\tilde D_i$ denotes
the cocore of $\tilde h_i$?
\Prop{
	For each $i\in I'$,
	if $w_2(X,\SpinStruct)([D_i,\partial D_i])=1$,
	then $w_2(Y,\allowbreak\pi^*\SpinStruct)([\tilde D_i,\tilde\partial D_i])$
	is always $1$ as well.
	Conversely,
	if $w_2(X,\SpinStruct)([D_i,\partial D_i])=0$,
	then $w_2(Y,\pi^*\SpinStruct)([\tilde D_i,\partial\tilde D_i])$
	is $1$ if the covering multiplicity $p$
	be even
	and $0$ otherwise.
}
\Proof{
	Firstly,
	consider the $2$-handle $D^2\times D^2$
	having the Spin structure $\SpinStruct_0$ defined
	along $D^2\times S^1$
	which does not extend across $D^2\times D^2$.
	Consider the $p$-fold branched cover
	$D^2\times D^2\to D^2\times D^2$
	given by rotating the second $D^2$ factor
	and keeping the first one fixed.
	One can trivialize the bundle
	$\Tan(D^2\times D^2)|_{D^2\times S^1}$
	by splitting it as
	$\R^3\oplus\Tan S^1$.
	It is clear now that this trivialization
	is pulled back to itself along the restriction
	of the covering to $D^2\times S^1$.
	This means that the Spin structure $\SpinStruct_0$
	is pulled back to itself
	and the lift also does {\it not\/}
	extend across $D^2\times D^2$.
	Secondly,
	consider instead the Spin structure
	$\SpinStruct$
	which extend across
	$D^2\times D^2$.
	Then,
	to understand the lifting behaviour,
	make use of the other Spin structure already considered,
	that is $\SpinStruct_0$,
	which does {\it not\/} extend across the covering.
	The obstruction class
	to a homotopy between
	$\SpinStruct$ and $\SpinStruct_0$
	can be seen as a non-trivial
	class in the cohomology
	$\Homology^1(S^1;\pi_1(\SO(4)))\cong\pi_1(\SO(4))$.
	One now easily sees that,
	under the $p$-fold covering,
	this class in $\pi_1(\SO(4))$
	gets traversed $p$ times.
	Since $\pi_1(\SO(4))$
	is the cyclic group in two elements,
	the result follows.
}
\Eg{
	The next figure continues the line of examples
	coming from the $(-8)$-surgery on the left-handed trefoil
	and shows the perhaps
	slightly unexpected behaviour in this case:
	the Spin structure does extend
	across the $(-8)$-framed handle downstairs
	but,
	upstairs,
	its lift does not extend across
	the lifted $2$-handle
	because the cover is a {\it double\/} cover.
}
\Fig{
	\Label\MinusEightLeftTrefoilDoubleCoverSpin
	\vbox{%
	\baselineskip = 0pt%
	\hbox{\hskip135bp}%
	\vbox{\vskip67.5bp}%
	\vtop{\vskip67.5bp}%
	\includegraphics{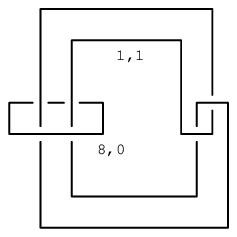}%
}%

	\vbox{
		\baselineskip=10pt
		\hbox to 1.5cm {\hfill $2:1$\hfill}
		\hbox to 1.5cm {\leftarrowfill}
	}
	\vbox{%
	\baselineskip = 0pt%
	\hbox{\hskip135bp}%
	\vbox{\vskip67.5bp}%
	\vtop{\vskip67.5bp}%
	\includegraphics{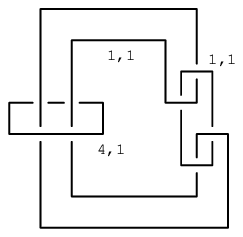}%
}%

}
\Prop{
	\Label\MinusEightLeftTrefoilDoubleCoverSpinLift
	Let $X$ be
	the $(-8)$-surgery
	on the left-handed trefoil.
	Let $\SpinStruct$
	be the Spin structure on $X$
	as defined by the Kirby diagram in
	{\MinusEightLeftTrefoilSpin}.
	Let $Y$ be the lens space $L(12,7)$.
	Denote by
	$\pi:Y\to X$
	the double covering.
	Let $M$ be the standard
	linearly plumbed $4$-manifold-with-boundary
	whose boundary is $L(12,7)$.
	Then,
	the Spin structure $\pi^*\SpinStruct$ on $Y$
	is the one for which $w_2(M,\pi^*\SpinStruct)=0$.
}
\Proof{
	Using the techniques developed above,
	this becomes a Kirby calculus affair.
	Starting from what is seen
	on the right hand side of
	\MinusEightLeftTrefoilDoubleCoverSpin,
	one computes,
	as shown in the following
	sequence of diagrams,
	(i) by blowing down of one of the $+1$-framed unknots,
	(ii) blowing up twice by $+1$-framed unknots
	and (iii) blowing down the $-1$-framed unknot.
	\par
	\hbox to \hsize{
		\hfill
		\vbox{
			\baselineskip=0pt
			\halign{
				\hfill#\hfill
				&
				\hfill#\hfill
				&
				\hfill#\hfill
				\cr
				\vbox{%
	\baselineskip = 0pt%
	\hbox{\hskip135bp}%
	\vbox{\vskip67.5bp}%
	\vtop{\vskip67.5bp}%
	\includegraphics{figs/minus_eight_left_trefoil_double_cover_spin.ps}%
}%

				&
				$\longrightarrow$
				&
				\vbox{%
	\baselineskip = 0pt%
	\hbox{\hskip135bp}%
	\vbox{\vskip52.5bp}%
	\vtop{\vskip52.5bp}%
	\includegraphics{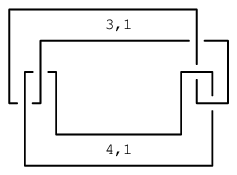}%
}%

				\cr
				&
				&
				$\Big\downarrow$
				\cr
				\vbox{%
	\baselineskip = 0pt%
	\hbox{\hskip150bp}%
	\vbox{\vskip56.25bp}%
	\vtop{\vskip56.25bp}%
	\includegraphics{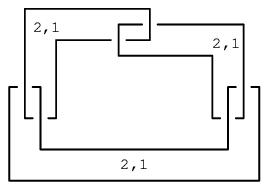}%
}%

				&
				$\longleftarrow$
				&
				\vbox{%
	\baselineskip = 0pt%
	\hbox{\hskip135bp}%
	\vbox{\vskip67.5bp}%
	\vtop{\vskip67.5bp}%
	\includegraphics{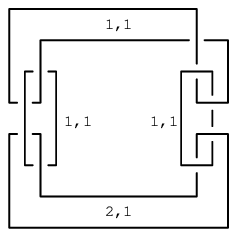}%
}%

				\cr
			}
		}
		\hfill
	}
	\par
	Now,
	if one pay careful attention
	at the signs of the crossings
	on the last diagram,
	one sees that one can
	perform a handle slide
	of any one of the $2$-framed handles
	over the $(-2)$-framed handle
	in order to obtain a linear chain of unknots.
	Thence,
	it is a straightforward matter
	to perform blow-ups and blow-downs
	to reach the standard plumbing diagram
	of the lens space $L(12,7)$.
	The procedure is outlined in the following diagrams.
	\par
	\hbox to \hsize{
		\hfill
		\vbox{
			\baselineskip=0pt
			\halign{
				\hfill#\hfill
				\cr
				\vbox{
					\halign{
						\hfill#\hfill
						&
						\hfill#\hfill
						\cr
						\vbox{%
	\baselineskip = 0pt%
	\hbox{\hskip135bp}%
	\vbox{\vskip45bp}%
	\vtop{\vskip45bp}%
	\includegraphics{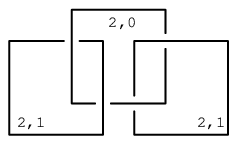}%
}%

						&
						\vbox{%
	\baselineskip = 0pt%
	\hbox{\hskip135bp}%
	\vbox{\vskip45bp}%
	\vtop{\vskip45bp}%
	\includegraphics{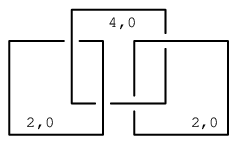}%
}%

						\cr
						$\searrow$
						&
						$\nearrow$
						\cr
					}
				}
				\cr
				\vbox{%
	\baselineskip = 0pt%
	\hbox{\hskip195bp}%
	\vbox{\vskip45bp}%
	\vtop{\vskip45bp}%
	\includegraphics{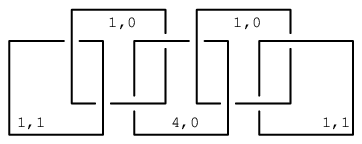}%
}%

				\cr
			}
		}
		\hfill
	}
}
\Prop{
	A $\SpinC$ structure on a $3$-manifold $Y$
	consists of precisely the same data
	as a complex trivialization of $\Tan Y\oplus\R$
	over the $2$-skeleton of $Y$
	which extend across its $3$-skeleton.
}
\Proof{
	Follows from obstruction theory.
	Cf. \GompfStipsiczFourManifolds,
	Remark 5.6.9(a),
	for the case of Spin structures.
	Also vid.\ \GompfHandlebody.
}
\Rem{
	From this characterization
	of $\SpinC$ structures on $Y$,
	it is easy to see that a Spin structure $\SpinStruct$
	defines a $\SpinC$ structure
	by taking the trivialization of
	$\Tan Y|_{Y_2}$ defined by $\SpinStruct$
	and picking the trivial complex structure
	on $\Tan Y|_{Y_2}\oplus\R$.
	One can then check that this trivialization
	extends to $\Tan Y|_{Y_3}$.
}
\Def{
	Given a Spin structure $\SpinStruct$,
	denote its induced $\SpinC$ structure
	by $\SpinStruct^\C$.
}
\Prop{
	(\GompfHandlebodyAlt, Theorem 4.12)
	\Label\GompfSpinCSpinDiffFormula
	Let $(Y,\lambda)$
	be a contact $3$-manifold
	given as Legendrian surgery on a Legendrian link
	$L:=\bigsqcup_{i\in I} K_i$ in $S^3$.
	Let $M$ be the $4$-manifold-with-boundary
	defined by the same Kirby diagram.
	Suppose $\SpinStruct$ be a Spin structure on $Y$
	and use $L'=\bigsqcup_{i\in I'}K_i\subset L$
	to denote the sublink of $L$
	consisting of those components $K_i$
	for which $w_2(Y,\SpinStruct)([D_i,\partial D_i])\neq 0$,
	where $D_i$ denotes the cocore of the handle attached along $K_i$.
	Recall that the author uses
	$\SpinCStruct_\lambda$
	to denote the $\SpinC$ structure
	defined by the contact form $\lambda$.
	Then,
	the difference class
	$$
		\SpinCStruct_\lambda-\SpinStruct^\C
		\in
		\Homology^2(Y;\Z)
	$$
	of $\SpinC$ structures
	is determined by the restriction to $Y$
	of a class $\rho\in\Homology^2(M;\Z)$
	defined by its evaluations
	on the classes $[F_i]\in\Homology_2(M;\Z)$,
	according to the formula
	$$
		(\SpinCStruct_\lambda
		-\SpinStruct^\C)([F_i])
		=
		{1\over2}\left(
			r(K_i)
			+
			\sum_{j\in I'}
			[F_i]\cdot[F_j]
		\right).
	$$
}
\Eg{
	\Label\MinusEightLeftTrefoilGompfGammaLift
	In the case of the Spin $3$-manifold $X$
	defined by {\MinusEightLeftTrefoilSpin}
	and studied in {\MinusEightLeftTrefoilDoubleCoverSpinLift},
	consider the $\SpinC$ structure $\SpinCStruct_\lambda$
	defined by the contact structure $\Ker\lambda$
	produced by Legendrian surgery on
	the Legendrian left-trefoil
	depicted in
	{\LegendrianLeftTrefoilTBMinusSevenRotZero},
	that is,
	having Thurston-Bennequin invariant $-7$
	and rotation number zero.
	Computing using
	\GompfSpinCSpinDiffFormula,
	one readily finds that
	$\SpinCStruct-\SpinStruct^\C=0$.
}
\Rem{
	Using {\GompfSpinCSpinDiffFormula}
	and the Kirby calculus techniques developed above,
	one can compute the lift of
	a $\SpinC$-structure $\SpinCStruct$
	via a finite covering
	by computing the lift of a Spin structure $\SpinStruct$
	and then computing the lift of
	the second cohomology class
	$\SpinCStruct-\SpinStruct^\C$.
}
\Thm{
	Consider the manifold $X$
	given by $(-8)$-surgery on the left-handed trefoil
	with tight contact structure $\Ker\lambda$ given
	by Legendrian surgery
	according to {\LegendrianLeftTrefoilTBMinusSevenRotZero}.
	The lift of $\Ker\lambda$
	via the double cover $\pi:Y\to X$
	is the tight contact structure $\tilde\lambda$
	on $Y\cong L(12,7)$
	given by Legendrian surgery
	on the middle diagram of
	{\LensTwelveSevenContact}.
}
\Proof{
	By {\ContactPrismExampleDThreeLiftMatch},
	$d_3(\pi^*\lambda)=d_3(\tilde\lambda)$.
	Meanwhile,
	{\MinusEightLeftTrefoilDoubleCoverSpinLift}
	and {\MinusEightLeftTrefoilGompfGammaLift}
	combine to assert that the $\SpinC$ structures
	also match;
	that is,
	$\pi^*\SpinCStruct_\lambda\cong\SpinCStruct_{\tilde\lambda}$.
	Since the $d_3$ invariant
	and the $\SpinC$ structure
	completely characterize the homotopy class
	of a hyperplane field,
	$\pi^*\lambda$ and $\tilde\lambda$ are homotopic.
	Now,
	according to \EllipticLiftHomotopicTight,
	this is a sufficient condition
	for $\pi^*\lambda$ and $\tilde\lambda$
	to be isotopic.
}
\Rem{
	To the best of the author's knowledge,
	there is no other way to obtain this result
	short of working in coordinates,
	which would probably prove itself to be untenable.
}

\Sect{Appendix}
This appendix shall deal with the proof of
{\NoTrajToContactConf},
where it was claimed that
there are no non-trivial Seiberg-Witten trajectories
with positive infinity limit the contact configuration.
A very similar assertion is made by
\TaubesWeinsteinII,
Proposition 5.15.
The key difference is that
Taubes deals with the case of a non-torsion
$\det\SpinCStruct_\lambda$,
which is not the case here.
The fact that 
$\det\SpinCStruct_\lambda$
is non-torsion in
\TaubesWeinsteinII,
Proposition 5.15,
also causes its conclusion
to be somewhat weaker than what is proved here.
Another minor difference is that,
in the present {\Document},
the author must not assume that
the negative infinity limit of the trajectory
be non-degenerate.
These differences require
only significant modifications
to the final part of the proof;
therefore,
the proof presented here
shall be very similar
to what one can find in
{\TaubesWeinsteinII}
and other works of Taubes
that have drawn significantly from it.
However,
Taubes' proof is difficult to follow
for someone who has not read
{\TaubesWeinsteinI} and {\TaubesWeinsteinII}
in their entirety;
therefore,
the author felt that it was necessary
to include some of the details here
for the reader's convenience.
The reader can also find a complete account
with self-contained proofs for {\it all\/} lemmata
in \RosoThesis, \S12.
\par
The overall strategy shall consist of the following.
On the one hand,
it shall be shown that only the contact configuration
and its gauge equivalent configurations
have the property that the spinor component
be bounded
{\it away from zero\/}
in a certain way.
On the other,
it shall be shown that the trajectory
necessarily satisfies the same sort of bound
on its spinor component
for {\it all time.}
As a consequence,
both endpoints shall have to be gauge equivalent.
\par
Henceforth,
the author shall deal
with Seiberg-Witten trajectories.
Suppose $\gamma:\R\to\Conf$
to be a finite type trajectory
satisfying
$$
	\dbyd t\gamma(t)=
	-\SW_{\lambda,r}\big(\gamma(t)\big).
$$
Assume that both limits
$\lim_{t\to\pm\infty}\gamma(t)$
be well defined
in the $\Sobolev5,2$
topology.
Note that,
by Sobolev embedding,
the limits are also well defined in $\DiffClass 2$.
The author shall not {\it yet\/} be assuming
that the positive time limit
be the contact configuration $C_\lambda$;
that assumption shall be added later.
Write
$
	\gamma(t)=(A(t),\psi(t))
	\in\Conf(Y,\SpinCStruct_\lambda)
$
for its connexion and spinor components
respectively.
Write $A(t):=A_\lambda(t)+a(t)$
where $a(t)$ is a purely imaginary $1$-form on $Y$.
Decompose the spinor as
$\psi(t)=(\alpha(t),\beta(t))$
according to the direct sum
$
	\Spinor_\lambda
	=\Exterior^{0,0}\xi^*
	\oplus\Exterior^{0,1}\xi^*
$.
\Lem{
	\Label\SecondOrdEqnForSpinor
	\LabelNum\SecondOrdEqnForSpinorNum
	The spinor $\psi$ satisfies the second order equation
	$$
		-\pdbydo t2\psi
		+\nabla_A^*\nabla_A\psi
		-\cliff\left(
			{1\over2}\Hodge F_{A_\lambda}
			+r\tau(\psi)
			-{ir\over2}\lambda
		\right)
		\psi
		+{s\over4}\psi
		=0,
	$$
	where $s:Y\to\R$
	denotes the scalar curvature of $Y$.
}
\Proof{
	The Dirac equation
	$$
		\pdbyd t\psi=-\Dirac_A \psi
	$$
	implies the second order equation
	$$
		\pdbydo t2\psi
		=-\pdbyd t\Dirac_A\psi
		=-\Dirac_A\pdbyd t\psi
		-\cliff\left(\pdbyd t A\right)\psi
		=\Dirac_A^2\psi
		-\cliff\left(\pdbyd t A\right)\psi.
	$$
	Now,
	apply the
	the well known Weitzenb\"ock formula
	(cf.\ \NicolaescuNotesAlt, (1.3.11)),
	$$
		\Dirac_A^*\Dirac_A \psi
		=\nabla_A^*\nabla_A \psi
		-{1\over2}\cliff(\Hodge F_A)\psi
		+{s\over4}\psi,
	$$
	to find that
	$$
		-\pdbydo t2\psi
		+\nabla_A^*\nabla_A\psi
		-\cliff\left(
			\pdbyd t A
			+{1\over2}\Hodge F_A
		\right)
		\psi
		+{s\over4}\psi
		=0.
	$$
	The result then follows by applying
	the curvature component
	of the Seiberg-Witten trajectory equations.
}
\Rem{
	\Label\TaubesVersionOfSWEqnsWithT
	A distinction worth mentioning
	between what is done here
	and what is done in {\TaubesWeinsteinII}
	is that
	some (but not all)
	of the results in {\TaubesWeinsteinII}
	concern a different version
	of the Seiberg-Witten trajectory equations.
	To be precise,
	this different version
	says
	$$
		\pdbyd t A
		=-{1\over2}\Hodge\left(
			F_A
			-F_{A_\lambda}
		\right)
		-re^{T(t)}\left(
			\tau(\psi)
			-{i\over2}\lambda
		\right),
		\quad
		\pdbyd t\psi=-\Dirac_A\psi,
	$$
	for a fixed function $T:\R\to[0,\infty)$
	satisfying $T(t)=0$
	for all $t\leq1$
	and $T(t)=t$
	for all $t\geq2$.
	The equations at use in the present {\Document} 
	are equivalent to setting $T=0$ everywhere.
	The use of such a $T$
	is certainly important in {\TaubesWeinsteinII}
	and brings about some complications
	which do not occur in the present context.
	Some of the results proved in {\TaubesWeinsteinII}
	for this version of the Seiberg-Witten equations
	are easily adapted to the present context
	as is remarked in the proof of
	\TaubesWeinsteinII, Proposition 5.15.
	In particular,
	for the next two lemmata,
	the author shall refer the reader to
	\TaubesWeinsteinII.
}
\Rem{
	In what follows,
	a sequence of bounds
	on various quantities
	shall be claimed.
	These shall involve constants
	$\Const_n$
	labelled by $n$
	the number of the lemma in which they appear.
	These quantities may depend on the metric of $Y$,
	and the contact form $\lambda$,
	but
	they shall not depend on the particular solution
	$(A,\psi)$ to the Seiberg-Witten equations
	{\it nor shall they shall not depend on the value
	of the parameter $r>0$.}
}
\Lem{
	\Label\PointwisePsiBound
	\LabelNum\PointwisePsiBoundNum
	There is a constant
	$\Const_\PointwisePsiBoundNum>0$
	such that
	$$
		|\alpha|^2+|\beta|^2\leq
		1+\Const_\PointwisePsiBoundNum r^{-1}.
	$$
}
\Proof{
	Analogous to the first statement
	of
	\TaubesWeinsteinII,
	Lemma 5.1
	with the function $T$
	from {\TaubesVersionOfSWEqnsWithT} set to zero.
	Alternatively,
	vid.\ \RosoThesis, Lemma 12.5.
}
\Lem{
	\Label\PointwiseBetaBound
	\LabelNum\PointwiseBetaBoundNum
	There are constants
	$\Const_\PointwiseBetaBoundNum>0$
	and $r_\PointwiseBetaBoundNum>0$
	such that
	for any $r>r_\PointwiseBetaBoundNum$
	the following inequality hold.
	$$
		|\beta|^2\leq
		\Const_\PointwiseBetaBoundNum\left(
			r^{-1}(1-|\alpha|^2)
			+r^{-2}
		\right).
	$$
}
\Proof{
	Analogous to the second statement
	of
	\TaubesWeinsteinII,
	Lemma 5.1
	with the function $T$
	from {\TaubesVersionOfSWEqnsWithT} set to zero.
	Alternatively,
	vid.\ \RosoThesis, Lemma 12.10.
}
\Lem{
	\Label\PointwiseCurvBound
	\LabelNum\PointwiseCurvBoundNum
	There is a constant
	$\Const_\PointwiseCurvBoundNum>0$
	such that
	$$
		\modulus{\pdbyd tA}
		+\modulus{F_A}
		\leq \Const_\PointwiseCurvBoundNum r.
	$$
}
\Proof{
	Analogous to
	\TaubesWeinsteinII,
	Proposition 5.2
	with the function $T$
	from {\TaubesVersionOfSWEqnsWithT} set to zero.
	Alternatively,
	vid.\ \RosoThesis, Lemma 12.25.
}
Henceforth,
add the assumption that
$$
	\lim_{t\to+\infty}\gamma=C_\lambda.
$$
The goal shall be bound $\modulus\alpha^2$
{\it away\/}
from zero
in the same manner as is the case for $C_\lambda$
but for the entirety of the trajectory.
\par
For that end,
start by introducing the symplectic form
$
	\omega:=
	e^{2t}(\d t\wedge\lambda+\Hodge\lambda),
$
on $\R\times Y$
and fixing some almost complex structure $J$
compatible with $\omega$.
Also fix a non-decreasing $\DiffClass\infty$ function
$\sigma:[0,\infty)\to[0,1]$
satisfying $\sigma|_{[0,1/2]}=0$
and $\sigma|_{[1,\infty)}=1$.
Use $\sigma'$ to denote its derivative.
Assume further that $\sigma':\R\to[0,3]$.
\Def{
	Given $\delta>0$,
	denote $\sigma_\delta:\R\times Y\to[0,1]$
	the function
	$$
		(t,y)\mapsto
		\sigma\left(
			\delta^{-1}
			\left(1-|\alpha|^2\right)
		\right).
	$$
}
\Def{
	Given $\delta>0$,
	denote $\sigma'_\delta:\R\times Y\to[0,1]$
	the function
	$$
		(t,y)\mapsto
		\sigma'\left(
			\delta^{-1}
			\left(1-|\alpha|^2\right)
		\right).
	$$
}
Now,
introduce $\nabla_A$ to denote the $\SpinC$ covariant derivative
on $\Spinor_\lambda$
induced by the connexion $A$ on $\det\SpinCStruct_\lambda$.
This is to say that Clifford multiplication
is covariantly constant for $\nabla_A$.
It shall also be necessary to work
with covariant derivatives defined on each of the summands
$
	\Spinor_\lambda
	\cong
	\Exterior^{0,0}\xi^*
	\oplus
	\Exterior^{0,1}\xi^*
$.
Both these derivatives shall be denoted
by the same symbol as
$$
	\nabla_A':
	\Gamma(\Exterior^{0,k}\xi^*)
	\to
	\Gamma(\Exterior^{0,k}\xi^*\otimes\Cotan Y)
$$
and shall be defined by projecting $\nabla_A$
to the respective summand of $\Spinor_\lambda$.
This can be conveniently expressed in terms of
Clifford multiplication by the formula
$$
	\nabla_A'=
	{1\over2}
	\left(
		1+i(-1)^k\cliff(\lambda)
	\right)
	\nabla_A.
$$
\par
It shall prove advantageous to work
with objects over $\R\times Y$
rather than time dependent objects over $Y$.
For that end,
introduce
$\hat\d$
and $\hat{\Hodge}$
to denote,
respectively,
the exterior derivative
and the Hodge operator
over
$\R\times Y$.
The $\SpinC$ structure
$\SpinCStruct_\lambda$
over $Y$
also defines a $\SpinC$ structure
over $\R\times Y$,
which shall be denoted $\hat\SpinCStruct_\lambda$.
The time dependent connexion $A$ over
$\det\SpinCStruct_\lambda$
also defines a connexion 
on the bundle $\det\hat\SpinCStruct_\lambda$;
denote this connexion by $\hat A$
and note that it is characterized by requiring
its induced covariant derivative to be
$$
	\nabla_{\hat A}:=\d t\otimes\pdbyd t+\nabla_A.
$$
Hence,
one sees that its curvature is given by
$$
	F_{\hat A}
	=\d t\wedge\pdbyd tA+F_A.
$$
The covariant derivative $\nabla'_A$
induces a covariant derivative
on the trivial complex line bundle over $\R\times Y$
defined by
$$
	\nabla_{\hat A}':=\d t\otimes{\pdbyd t} + \nabla_A'.
$$
This covariant derivative can then be extended
to a covariant exterior derivative on complex differential forms,
which shall be denoted
$$
	\hat\d_{\hat A}:
	\Omega^k(\R\times Y)\otimes\C
	\to\Omega^{k+1}(\R\times Y)\otimes\C.
$$
\Def{
	Given $\delta>0$,
	define a $2$-form
	$\wp_\delta\in\Omega^2(\R\times Y)$
	by the formula
	$$
		\wp_\delta:=
		{1\over\delta}\sigma'_\delta
		\cdot
		\hat\d_{\hat A}\alpha
		\wedge
		\hat\d_{\hat A}\conj\alpha
		+\sigma_\delta\cdot
		F_{\hat A}.
	$$
}
\Rem{
	Notice that the form $\wp_\delta$ is closed.
}
\Lem{
	\Label\WPSuppBoundedAbove
	For all
	$\delta\in(0,\delta_\SpinorContactConfBoundedbelowNum)$,
	there exists
	$s_\delta\in\R$,
	such that,
	the $2$-form $\wp_\delta$
	vanishes identically on
	$[s_\delta,\infty)\times Y\subset \R\times Y$.
}
\Proof{
	By {\SpinorContactConfBoundedbelow}
	and the fact that
	$\lim_{t\to\infty}\gamma=\ContactConf_\lambda$
	in the $\DiffClass 0$ topology,
	it follows that there exists
	$s_\delta\in\R$
	such that,
	for $t>s_\delta$,
	$|\psi(t)|\geq 1-\delta_\SpinorContactConfBoundedbelowNum$.
	Since the function
	$\sigma$ is supported on $[1/2,\infty)$,
	the claimed result follows.
}
\Lem{
	For $r>r_\SpinorContactConfBoundedbelowNum$
	and
	$\delta\in(0,\delta_\SpinorContactConfBoundedbelowNum)$,
	$\omega\wedge\wp_\delta$
	is integrable over all of $\R\times Y$.
}
\Proof{
	Notice that
	there is some $K>0$,
	which may well depend on $r$,
	$A$ or $\psi$ here,
	such that,
	for sufficiently negative $s<0$,
	$|\omega\wedge\wp_\delta|$
	restricted to $(-\infty,s)$
	is no greater than 
	$Ke^s$.
	This is a consequence of the fact that
	the limit
	$\lim_{t\to-\infty}\gamma(t)$
	is well defined in the Sobolev norm $\Sobolev5,2$
	and,
	therefore,
	also in $C^2$
	by Sobolev embedding.
	Together with
	{\WPSuppBoundedAbove},
	this implies integrability.
}
\Lem{
	\Label\OmegaWPIntZero
	Provided $r>r_\SpinorContactConfBoundedbelowNum$
	and 
	$\delta\in(0,\delta_\SpinorContactConfBoundedbelowNum)$,
	it follows
	$\int_{\R\times Y}\omega\wedge\wp_\delta=0$.
}
\Proof{
	Note that $\omega$ is exact
	and $\wp_\delta$ is closed
	and integrate by parts.
}
Introduce the notation
$$
	\hd_{\hat A}:
	\Omega^{i,j}(\R\times Y)
	\to
	\Omega^{i+1,j}(\R\times Y),
	\quad
	\ahd_{\hat A}:
	\Omega^{i,j}(\R\times Y),
	\to
	\Omega^{i,j+1}(\R\times Y),
$$
for the covariant Cauchy-Riemann operators associated
to the connexion $\hat A$
and the almost complex structure $J$.
It can be verified that,
for $\eta\in\Omega^{0,0}(\R\times Y)$,
they take the form
$$
	\hd_{\hat A}\eta:=
	{1\over2}\left(
		\hat\d_{\hat A}\eta
		+i e^{-2t}\hat{\Hodge}(\omega
		\wedge
		\hat\d_{\hat A}\eta)
	\right),\quad
	\ahd_{\hat A}\eta:=
	{1\over2}\left(
		\hat\d_{\hat A}\eta
		-i e^{-2t}\hat{\Hodge}(\omega
		\wedge
		\hat\d_{\hat A}\eta)
	\right).
$$
\Lem{
	\Label\OmegaWPIntZeroExpanded
	Vanishing of the integral
	in {\OmegaWPIntZero}
	amounts to saying
	$$
		\int_{\R\times Y}
			\hat{\Hodge}e^{2t}\left(
			\delta^{-1}
			\sigma_\delta'\left(
				|\hd_{\hat A}\alpha|^2
				-|\ahd_{\hat A}\alpha|^2
			\right)
			+
			r\sigma_\delta\left(
				1-|\alpha|^2+|\beta|^2
			\right)
		\right)
		=0.
	$$
}
\Proof{
	Can be checked directly
	by expanding with the expressions given above
	for $\hd_{\hat A}$ and $\ahd_{\hat A}$
	acting on $0$-forms
	and by using the four-dimensional Seiberg-Witten curvature equation.
}
\def\GaussChart{\phi}
\Lem{(cf.\ \TaubesWeinsteinIIAlt, Lemma 5.13)
	\Label\AHDAlphaBound
	\LabelNum\AHDAlphaBoundNum
	For a given
	$\delta\in(0,\delta_\SpinorContactConfBoundedbelowNum)$,
	there exist $r_\AHDAlphaBoundNum>0$
	and $\Const_\AHDAlphaBoundNum>0$
	such that,
	for all $r>r_\AHDAlphaBoundNum$,
	$$
		|\ahd_{\hat A}\alpha|\leq\Const_\AHDAlphaBoundNum.
	$$
}
\Proof{
	It is well known
	(cf.\ \NicolaescuNotesAlt, \S1.4.3),
	that,
	in terms of the notation above,
	the Dirac equation,
	$$
		\pdbyd t \psi
		=-\Dirac_A\psi
	$$
	takes the more familiar form
	$$
		\ahd_{\hat A}\alpha + \ahd_{\hat A}^*\beta = f_1(\psi),
	$$
	where $f_1$ is some bundle homomorphism
	independent of $r$.
	Hence,
	a bound on $\ahd_{\hat A}\alpha$
	of the variety required
	follows from a bound
	on $\pdbyd t\beta$,
	$\nabla_A\beta$
	and $\modulus\psi$.
	The bound on $\modulus\psi$
	was already provided by \PointwisePsiBound.
	The rest of the proof
	shall concern itself with the other two bounds.
	\par
	Let $r_\AHDAlphaBoundNum$
	be large enough so that
	the injectivity radius of
	$Y$
	be strictly
	larger than $(r_\AHDAlphaBoundNum)^{-1/2}$
	and so that the bundle
	$\Exterior^{0,2}\Cotan(\R\times Y)$
	be trivial when restricted
	to any ball of radius
	$(r_\AHDAlphaBoundNum)^{-1/2}$
	in $\R\times Y$.
	Then,
	let $r>r_\AHDAlphaBoundNum$,
	and fix a point $p\in \R\times Y$.
	Introduce
	$\GaussChart_{p,r}$
	to denote the Gaussian chart
	centred at $p$
	rescaled so that the ball of radius
	$1$ of $\R^4$
	be mapped to the geodesic ball of radius
	$r^{-1/2}$.
	Now,
	let
	$\tilde\beta:\Ball(\R^4,1)\to\C$
	denote the pullback of $\beta$
	via this chart  
	seen as a complex valued function
	by trivialising the bundle
	$\Exterior^{0,2}\Cotan(\R\times Y)$
	on this chart
	via the parallel transport map
	of the connexion $\hat A$.
	One can check that
	{\SecondOrdEqnForSpinor}
	implies a certain second order equation
	for $\tilde\beta$
	of the form
	$$
		\Delta\tilde\beta
		+\sum_{j=1}^4 f_{j+1}\pdbyd{x_j}\tilde\beta
		+f_6\tilde\beta
		+r^{-1}f_7
		=0,
	$$
	where $f_2,\ldots,f_7$
	are complex valued functions
	with norms bounded above by some constant
	$K_1>0$
	independent of $r$ and the point $p$.
	Recall {\PointwiseBetaBound};
	it implies that,
	for some constant $K_2\geq K_1$,
	also independent of $r$ and $p$,
	that $|\beta|\leq Kr^{-1/2}$.
	Standard elliptic theory then provides a bound
	of the form
	$$
		\left|\pdbyd x_j\tilde\beta(0)\right|
		\leq K_3 r^{-1/2}
	$$
	where $K_3\geq K_2$
	is some potentially larger constant
	independent of $r$ and $p$.
	Since the chart at hand was scaled
	so that the ball of radius $1$
	be mapped to the geodesic ball of radius $r^{-1/2}$,
	and $\Exterior^{0,2}\Cotan(\R\times Y)$
	was trivialized by the parallel transport
	of the connexion $\hat A$,
	whose curvature satisfies the bound of
	{\PointwiseCurvBound},
	if the reader care to check,
	it follows that,
	for some constant $K_4\geq K_3$
	independent of $r$ and $p$,
	$$
		\left|\pdbyd t\beta\right|
		+|\nabla_A\beta|
		\leq K_4.
	$$
}
The next few results needed
shall require mention
of a variant of the
{\it vortex equations\/}
on $\R^4=\C^2$.
These are defined next.
In what follows,
use $\omega_0$
to denote the standard symplectic form
on $\R^4$.
\Def{
	\Label\VortexEqns
	Consider a pair
	$(a_0, \alpha_0)$,
	where $a_0\in i\Omega^1(\R^4)$
	and $\alpha_0\in \Omega^{0,0}(\R^4)$.
	Then,
	$(a_0,\alpha_0)$
	satisfy the
	{\it vortex equations
	with bound $K>0$\/}
	when
	$$
		\ahd\alpha_0=0,\quad
		|\alpha_0|\leq 1,\quad
		\d^+ a={1\over2}(1-|\alpha_0|^2)\omega_0,\quad
		|\d^-a|\leq K.
	$$
}
\Lem{
	(\TaubesWeinsteinIIAlt, Lemma 5.14)
	\Label\VortexFormBoundVolumes
	\LabelNum\VortexFormBoundVolumesNum
	For any $K>0$
	and $\delta\in(0,1)$,
	there exist
	$R_\VortexFormBoundVolumesNum>2$
	and $\Const_\VortexFormBoundVolumesNum>1$
	such that,
	for any vortex $(a_0,\alpha_0)$
	with bound $K$,
	if one use $V$ to denote the volume of the set
	$$
		\SetComp{
			x\in\Ball(
				\R^4,
				R_\VortexFormBoundVolumesNum
			)
		}{
			(1-|\alpha_0(x)|^2)>\delta
		}
	$$
	and $V'$ denote the volume of the set
	$$
		\SetComp{
			x\in\Ball(
				\R^4,
				\Half
				R_\VortexFormBoundVolumesNum
			)
		}{
			\delta
			>(1-|\alpha_0|^2)
			\geq\Half\delta
		},
	$$
	then it follows that
	$V'\leq\Const_\VortexFormBoundVolumesNum V$.
}
\Rem{
	The ability to apply this lemma
	in the present context comes from the following.
}
\Lem{
	\Label\VortexCloseMonopole
	\LabelNum\VortexCloseMonopoleNum
	There exists
	$\Const_\VortexCloseMonopoleNum>0$
	such that,
	for any
	$R\geq 1$
	and
	$\epsilon>0$,
	there exists
	$r_\VortexCloseMonopoleNum>0$
	such that,
	for all
	$p\in \R\times Y$
	and $r>r_\VortexCloseMonopoleNum$,
	it follows that
	there exists a vortex $(a_0,\alpha_0)$
	bounded by $\Const_\VortexCloseMonopoleNum$
	and a gauge transformation $u$
	such that
	$$
		\norm{
			(a_0,(\alpha_0,0))
			-(\GaussChart_{p,r,R})^*u\cdot(a,(\alpha,\beta))
		}_{\DiffClass0(\Disk(\R^4,R))}
		<\epsilon,
	$$
	where $\GaussChart_{p,r,R}$
	denotes the rescaled Gaussian chart centred at $p$
	so that the geodesic ball of radius $Rr^{-1/2}$
	be mapped to the ball of radius $R$ of $\R^4$.
}
\Proof{
	According to \TaubesWeinsteinII,
	the proof is an adaptation of
	a similar statement made in \TaubesSWToGr.
	The idea is as follows.
	Assume the contrary.
	Then,
	for any
	$\Const_\VortexCloseMonopoleNum>0$,
	there exists $R\geq 1$,
	$\epsilon>0$,
	an increasing unbounded sequence $\{r_n\}$,
	a sequence of points $p_n\in\R\times Y$
	and a sequence of Seiberg-Witten trajectories
	$(A_\lambda+a_n,(\alpha_n,\beta_n))$
	with parameter $r=r_n$
	such that,
	if one define
	$(\tilde a_n,(\tilde\alpha_n,\tilde\beta_n))$
	to mean the pullback of
	$(a_n,(\alpha_n,\beta_n))$
	via the chart $\GaussChart_{p_n,r_n,R}$
	seen as
	functions of the ball of radius $R$ of $\R^4$,
	then $(\tilde a_n,(\tilde\alpha_n,\tilde\beta_n))$
	do not lie in a ball of radius $\epsilon$
	around any vortex
	bounded by $\Const_\VortexCloseMonopoleNum$
	in the $\DiffClass0(\Disk(\R^4,R))$ norm.
	Now,
	the Seiberg-Witten equations imply that,
	after redefining
	$(\tilde a_n,(\tilde\alpha_n,\tilde\beta_n))$
	by applying some gauge transformation $u_n$,
	these functions obey {\it elliptic\/}
	equations of the form
	$$
		\d^+ \tilde a_n=f_1(\tilde\psi_n,\tilde\psi_n)+f_2,
		\quad

		\d^*\tilde a_n=0,
		\quad

		\ahd\tilde\alpha_n
		+\ahd^*\tilde\beta_n
		=f_3(\tilde a_n,\tilde\psi_n),
	$$
	for certain polynomial functions
	$f_1$, $f_2$ and $f_3$.
	On the other hand,
	by using {\PointwisePsiBound},
	{\PointwiseBetaBound},
	{\PointwiseCurvBound}
	and standard elliptic theory arguments,
	one can conclude that
	there must be a subsequence
	uniformly convergent on the disk of radius $R$.
	Let the limit of this convergent subsequence
	be denoted
	$(\tilde a_\infty, (\tilde\alpha_\infty,\tilde\beta_\infty))$.
	The point now is that,
	by {\PointwiseBetaBound},
	one can see further that
	$\beta_\infty=0$.
	Hence,
	if one care to check,
	it follows,
	by examining carefully the terms
	$f_1$,
	$f_2$
	and $f_3$,
	that
	the equations above reduce to the vortex equations
	in the limit $n\to\infty$.
	Consequently,
	$(a_\infty, \alpha_\infty)$
	satisfies the vortex equations.
	But then note that,
	if
	$\Const_\VortexCloseMonopoleNum$
	be large enough
	relative to the constant
	$\Const_\PointwiseCurvBoundNum$
	provided by
	{\PointwiseCurvBound},
	$(a_\infty, \alpha_\infty)$
	is a vortex bounded by
	$\Const_\VortexCloseMonopoleNum$,
	which is a contradiction.
}
Set the following notation
$$
	\Omega_\delta
	:=\int_{\R\times Y}
		\hat{\Hodge}e^{2t}\sigma_\delta,
	\quad
	\Omega'_\delta
	:=\int_{\R\times Y}
		\hat{\Hodge}e^{2t}\sigma'_\delta.
$$
\Lem{
	(cf.\ \TaubesWeinsteinIIAlt, Lemma 5.17)
	\Label\OmegaDashOmegaBound
	\LabelNum\OmegaDashOmegaBoundNum
	Given
	$\delta\in(0,\delta_\SpinorContactConfBoundedbelowNum)$,
	there exist $r_\OmegaDashOmegaBoundNum>0$
	and $\Const_\OmegaDashOmegaBoundNum>0$
	such that,
	for all $r>r_\OmegaDashOmegaBoundNum$,
	$$
		\Omega_\delta'
		\leq
		\Const_\OmegaDashOmegaBoundNum\Omega_\delta.
	$$
}
\Proof{
	The analogous statement in {\TaubesWeinsteinII}
	lacks a proof as it is claimed to be similar
	to a previous lemma;
	therefore,
	the author shall furnish the details here.
	Firstly,
	notice that the function $\sigma_\delta$
	is non-zero only at points where
	$(1-|\alpha|^2)\geq\delta$;
	meanwhile,
	$\sigma'_\delta$
	is non-zero only at points where
	$\delta\geq(1-|\alpha|^2)\geq\Half\delta$.
	Now,
	If $g$ be the product metric of $\R\times Y$
	consider instead the metric $e^{2t}g$.
	Since the set of points where
	$(1-|\alpha|^2)\geq\delta$
	has $t$ coordinate bounded above,
	the volume of this set
	in the metric $e^{2t}g$
	is,
	in fact,
	finite.
	Hence,
	note that,
	in order to prove the claimed result,
	it suffices to bound the volume of points
	where
	$\delta\geq(1-|\alpha|^2)\geq\Half\delta$
	by some constant times
	the volume of the set of points where 
	$(1-|\alpha|^2)\geq\delta$;
	both volumes being with respect to the metric
	$e^{2t}g$.
	For the rest of this proof,
	the metric at use shall be $e^{2t}g$
	whenever the author talk of volumes or geodesy.
	For brevity,
	use $R:=R_\VortexFormBoundVolumesNum$
	to denote the constant from
	{\VortexFormBoundVolumes}.
	Consider a set $\Lambda$
	of disjoint geodesic balls in $Y$
	centred at points $p$ where
	$\delta\geq(1-|\alpha(p)|^2)\geq\Half\delta$
	and all having radius
	$\Frac14Rr^{-1/2}$.
	Due to compactness of $Y$,
	one can also assume
	$\Lambda$ to be maximal
	with respect to inclusion.
	For each ball $B\in\Lambda$,
	let $B''\supset B'\supset B$ denote
	the concentric geodesic balls
	having radii $Rr^{-1/2}$
	and $\Half Rr^{-1/2}$
	respectively.
	Now,
	suppose that
	some point $p\in \R\times Y$
	at which
	$\delta\geq(1-|\alpha(p)|^2)\geq\Half\delta$
	not be in the set
	$\bigcup_{B\in\Lambda}B'$.
	Then,
	if $r$ be sufficiently large
	compared to $R$,
	say larger than some
	$r_\OmegaDashOmegaBoundNum>0$,
	the ball of radius $\Frac14Rr^{-1/2}$
	centred at $p$
	would not intersect any of the balls
	in the set $\Lambda$;
	so,
	if this were the case,
	$\Lambda$ could not be maximal.
	Therefore,
	$\bigcup_{B\in\Lambda}B'$
	covers the set of points where
	$\delta\geq(1-|\alpha|^2)\geq\Half\delta$.
	Moreover,
	perhaps after an increase to 
	$r_\OmegaDashOmegaBoundNum>0$,
	Riemannian geometry provides
	an upper bound
	for the maximal number $n$ such that
	there be a set of balls $\{B_i\}\subset\Lambda$
	satisfying
	$B''_1\cap\cdots\cap B''_n\neq\emptyset$;
	this upper bound
	is independent of $\delta$ and $r$
	except that one must ensure
	$r>r_\OmegaDashOmegaBoundNum$.
	The consequence of all of this is that
	it suffices to provide the desired sort of bound
	for each of the balls $B$.
	For that end,
	given a ball $B\in\Lambda$,
	denote by $V_B$
	the volume of the subset of $B''$
	where
	$\delta\geq(1-|\alpha(p)|^2)\geq\Half\delta$;
	likewise,
	denote by $V'_B$
	the volume of the subset of $B'$
	where
	$(1-|\alpha|^2)\geq\delta$.
	Now,
	apply
	{\VortexFormBoundVolumes}
	in combination with {\VortexCloseMonopole}
	in order to obtain a constant
	$\Const_\OmegaDashOmegaBoundNum>0$
	independent of $B$,
	satisfying
	$V'_B\leq\Const_\OmegaDashOmegaBoundNum V_B$.
}
\Lem{
	\Label\OmegaDeltaZero
	\LabelNum\OmegaDeltaZeroNum
	Given
	$\delta\in(0,\delta_\SpinorContactConfBoundedbelowNum)$,
	there exists $r_\OmegaDeltaZeroNum$
	such that
	for all $r>r_\OmegaDeltaZeroNum$,
	$\Omega_\delta=0$.
}
\Proof{
	To start,
	let
	$
		r_\OmegaDeltaZeroNum
		=\max\{
			r_\OmegaDashOmegaBoundNum,
			r_\AHDAlphaBoundNum
		\}
	$.
	Recall that
	{\OmegaWPIntZeroExpanded}
	asserted that
	$$
		\int_{\R\times Y}
			\hat{\Hodge}e^{2t}\left(
			\delta^{-1}
			\sigma_\delta'\left(
				|\hd_{\hat A}\alpha|^2
				-|\ahd_{\hat A}\alpha|^2
			\right)
			+
			r\sigma_\delta\left(
				1-|\alpha|^2+|\beta|^2
			\right)
		\right)
		=0.
	$$
	Neglecting some positive terms yields the inequality
	$$
		\int_{\R\times Y}
			\hat{\Hodge}e^{2t}\left(
			-\delta^{-1}
			\sigma_\delta'
			|\ahd_{\hat A}\alpha|^2
			+
			r\sigma_\delta
			(1-|\alpha|^2)
		\right)
		\leq 0.
	$$
	Focusing on the second term of the integrand,
	note that,
	at a point where $\sigma_\delta\neq0$,
	it is necessarily the case that
	$(1-|\alpha|^2)\geq\delta$;
	hence,
	$$
		\int_{\R\times Y} \hat{\Hodge}e^{2t}\left(
			-\delta^{-1}
			\sigma_\delta'
			|\ahd_{\hat A}\alpha|^2
			+
			r\sigma_\delta\delta
		\right)
		\leq 0.
	$$
	By applying
	{\OmegaDashOmegaBound}
	and {\AHDAlphaBound},
	one finds
	$$
		\big(
			-\Const_\OmegaDashOmegaBoundNum
			\Const_\AHDAlphaBoundNum^2
			\delta^{-1}
			+r\delta
		\big)
		\Omega_\delta
		\leq 0.
	$$
	But $\Omega_\delta\geq0$.
	Hence,
	perhaps after increasing 
	$r_\OmegaDeltaZeroNum$,
	it follows that,
	for $r>r_\OmegaDeltaZeroNum$,
	$\Omega_\delta=0$.
}
\Thm{
	\Label\AlphaBoundedFromBelow
	Given
	$\delta\in(0,\delta_\SpinorContactConfBoundedbelowNum)$,
	for all $r>r_\OmegaDeltaZeroNum$,
	it follows that,
	pointwise on all of $\R\times Y$,
	$$
		1-|\alpha|^2\leq\delta.
	$$
}
\Proof{
	$\Omega_\delta$ 
	is the integral of
	$e^{2t}\sigma_\delta$,
	which is non-negative
	and strictly positive
	wherever $(1-|\alpha|^2)>\delta$.
}
\Cor{
	$C=\lim_{t\to-\infty}\gamma$ is gauge equivalent to $\ContactConf_\lambda$.
}
\Proof{
	$C$ is a solution to the Seiberg-Witten equations on $Y$.
	Its $\Exterior^{0,0}\xi^*$ component
	is $\lim_{t\to-\infty}\alpha$.
	Therefore,
	it also satisfies the bound in {\AlphaBoundedFromBelow}.
	{\SpinorContactConfBoundedbelow}
	then guarantees that $C$
	must be gauge equivalent to $\ContactConf_\lambda$.
}
\Cor{
	$\gamma$ is the constant trajectory at $\ContactConf_\lambda$.
}
\Proof{
	Assume the contrary.
	Because the Seiberg-Witten vector field
	$\SW_{\lambda,r}$
	is minus the gradient of the functional
	$\CSD_{\lambda,r}$,
	the value of $\CSD_{\lambda,r}$
	never increases along the trajectory $\gamma$.
	Moreover,
	the contact configuration $\ContactConf_\lambda$
	is a non-degenerate fixed point
	of the downward gradient flow of
	$\CSD_{\lambda,r}$
	by {\ContactConfNonDegen},
	which means that
	the value of $\CSD_{\lambda,r}$
	certainly decreases along $\gamma$
	for sufficiently large time
	as one approaches $\ContactConf_\lambda$.
	But
	since both endpoints of $\gamma$
	are gauge equivalent
	and $b_1=0$,
	the values of $\CSD_{\lambda,r}$
	are the same
	for gauge equivalent configurations.
	This cannot be.
}
\Rem{
	Hereby,
	the author concludes the proof of
	{\NoTrajToContactConf}.
}

\vfill
\eject
\Sect{Bibliography}
\Bib
\bye